\documentclass[11pt,psamsfonts]{amsart}
\usepackage{graphicx}
\usepackage{url}
\usepackage{amssymb}
\newtheorem{thm}{Theorem}[section]
\newtheorem{lem}{Lemma}[section]
\newtheorem*{conj}{Conjecture}
\newcommand{\dom}{\mathop{\mathrm{dom}}}

\newcommand{\VS}{\mathcal{VS}}
\newcommand{\SG}{\mathcal{SG}}
\newcommand{\E}{\mathcal{E}}
\numberwithin{equation}{section}

\title[The Laplacian and Spectral Operators on the Vicsek Set]{Analysis of the Laplacian and Spectral Operators on the Vicsek Set}
\author{Sarah~Constantin}
\address[Sarah~Constantin]
        {Department of Mathematics\\
        Fine Hall, Washington Road\\
        Princeton University\\
        Princeton NJ  08544-1000}
\email{sconstan@princeton.edu}

\author{Robert~S.~Strichartz}
\address[Robert~S.~Strichartz]
        {Mathematics Department, Malott Hall\\
        Cornell University\\
        Ithaca, NY  14853}
\email{str@math.cornell.edu}

\author{Miles~Wheeler}
\address[Miles~Wheeler]
        {Department of Mathematics\\
        Box 1917\\
        Brown University\\
        Providence, RI  02912}
\email{mhw33@cornell.edu}

\thanks{The research of the first and third authors was supported by
the National Science Foundation through the Research Experiences for
Undergraduates Program at Cornell University.}
\thanks{The research of the second author was supported in part by the
National Science Foundation, grant DMS-0652440.}

\begin{document}

\maketitle
\begin{abstract}
  We study the spectral decomposition of the Laplacian on a family of
  fractals $\mathcal{VS}_n$ that includes the Vicsek set for $n=2$,
  extending earlier research on the Sierpinski Gasket.  We implement
  an algorithm \cite{zhou} for spectral decimation of eigenfunctions
  of the Laplacian, and explicitly compute these eigenfunctions and
  some of their properties.  We give an algorithm for computing inner
  products of eigenfunctions.  We explicitly compute solutions to the
  heat equation and wave equation for Neumann boundary conditions.  We
  study gaps in the ratios of eigenvalues and eigenvalue clusters. We
  give an explicit formula for the Green's function on
  $\mathcal{VS}_n$. Finally, we explain how the spectrum of the
  Laplacian on $\mathcal{VS}_n$ converges as $n \to \infty$ to the
  spectrum of the Laplacian on two crossed lines (the limit of the
  sets $\mathcal{VS}_n$.) 
\end{abstract}

\tableofcontents

\section{Introduction}\label{sec:intro}
Kigami \cite{kigami} has developed a theory of Laplacians on a class
of fractals called \emph{pcf self-similar fractals}.  One example, the
Sierpinski gasket $\SG$ has become the ``poster child'' for
this theory \cite{fractals} in the belief that it is the simplest
nontrivial example.  As a result, a lot of very concrete results have
been obtained for $\SG$. This paper extends some of these
concepts and results to a different family of finitely ramified
self-similar fractals, the Vicsek sets $\VS_n$, with $n=2$
corresponding to the Vicsek set $\VS$.  We also obtain results
for $\VS$ that have no analogs on $\SG$.

To review the standard theory, a pcf self-similar fractal $V$ will be
a compact set in the plane, defined as the limit of a sequence of
graphs $\Gamma_0, \Gamma_1, \ldots$ with vertices $V_0 \subset V_1
\subset \cdots$.  The property of self-similarity takes the form of a
family of mappings from $V$ to itself, $\{ F_i \}$ which are
contractive similarities and have the property that $V = F_0(V) \cup
F_1(V) \cup \cdots \cup F_k(V)$.  For example, the Sierpinski Gasket
is defined by three similarities, each of which sends the entire set
$\SG$ to one of its three smaller triangular component copies.  We
refer to the graph at stage $m$ of the approximation as the
\emph{$m$th level graph approximation}.  The Vicsek set (specifically
the second order Vicsek set $\VS_2$, but sometimes simply called the
Vicsek set) is the fractal defined by the similarities $F_i\colon
\VS_2 \to \VS_2$
 \begin{align*}
   F_1(x)  & = x/3 & 
   F_4(x)  & = x/3 + 2/3 (0, 1) \\
   F_2(x)  & = x/3 + 2/3 (1, 0) & 
   F_0(x)  & = x/3 + 2/3 (1/2, 1/2) \\
   F_3(x)  & = x/3 + 2/3 (1, 1) & 
 \end{align*}
The first graph approximation $\Gamma_0$ is the complete graph on four
vertices (that is, the vertices of a unit square and an edge
connecting every pair of vertices).  The next approximation $\Gamma_1$
consists of five miniature copies of $\Gamma_0$ arranged in an X shape
with branches of length 2 (hence $\VS_2$).  Further graph
approximations likewise consist of five copies of the previous level;
they display finer levels of branching.  
\begin{figure}
  \begin{minipage}{1.8in}
    \centering
    \includegraphics[width=1.8in]{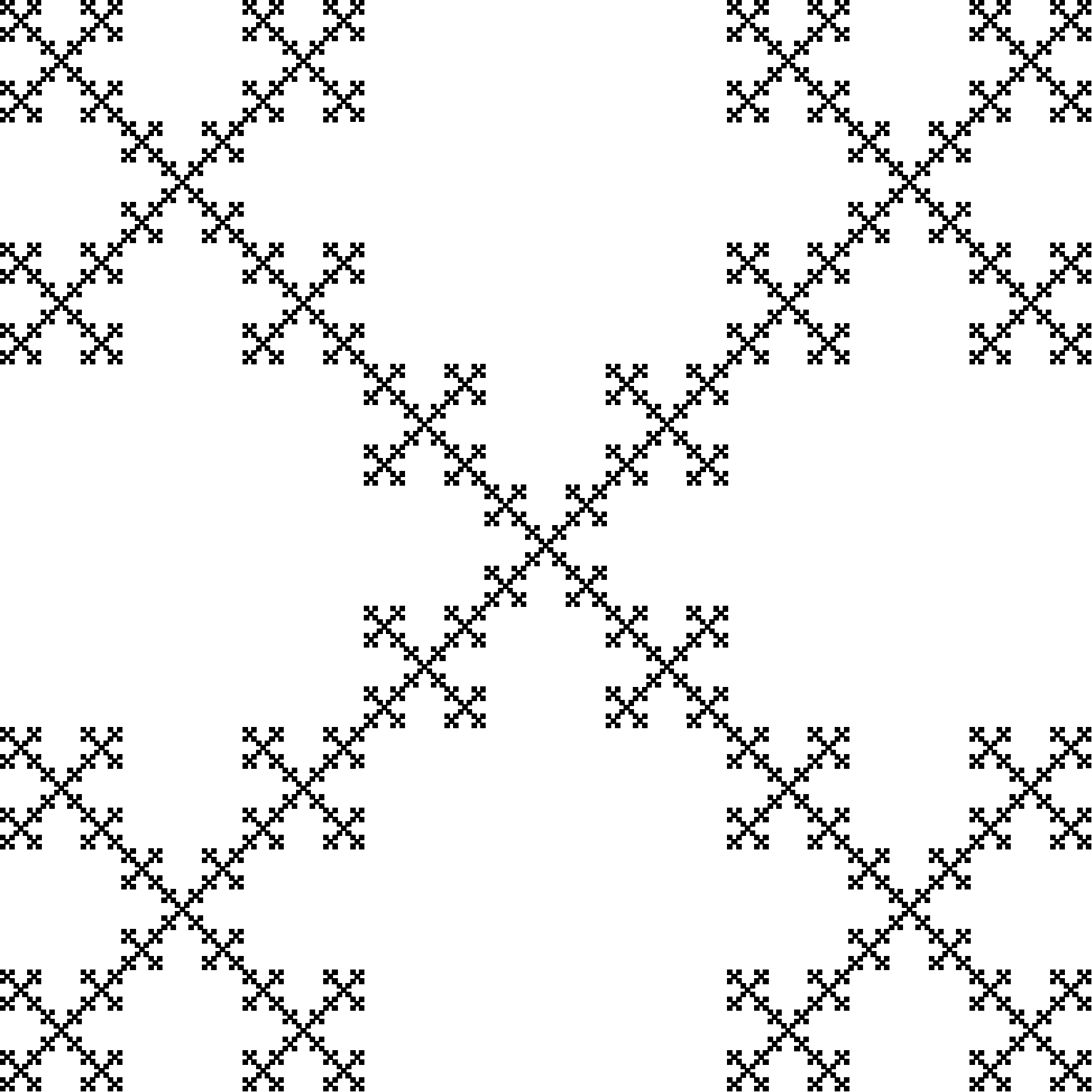}
    \caption{$\VS_2$}
  \end{minipage}
  \qquad
  \begin{minipage}{1.8in}
    \centering
    \includegraphics[width=1.8in]{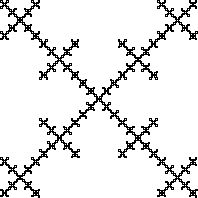}
    \caption{$\VS_3$}
  \end{minipage}
\end{figure}
Higher order Vicsek sets $\VS_n$ are similar, except that $\Gamma_1$
is an X-shaped graph consisting not of five but $4n -3$ copies of
$\Gamma_0$, with arms of length $n$.  Instead of five similarities, we
have $4n-3$ similarities.  

It is intuitive from the picture and also easy to demonstrate that as
$n \to \infty$, $\VS_n$ approaches the pair of crossed line segments
between $(0,0)$ and $(1, 1)$ and $(1,0)$ and $(0,1)$.  (That is, the
maximum Euclidean distance of any point in $\VS_n$ from the crossed
lines approaches zero.)  This is important to note because it suggests
a connection between fractal analysis on the Vicsek sets and classical
analysis on the line; later in this paper we show that the spectrum of
the (Neumann) Laplacian as defined on the Vicsek sets does, in fact,
approach the spectrum for the classical Neumann Laplacian on the
cross.  

On $\VS_n$, we can define a standard self-similar probability measure
as follows: for each graph approximation, let $\nu_m$ be the
probability measure which weights each vertex by its degree.  Then the
standard measure $\mu$ on $\VS_n$ is defined by
\begin{align*}
  \int_{\VS_2} f d\mu = \lim_{m \to \infty} \int_{\Gamma_m} f d\nu_m.
\end{align*}

We define the unrenormalized energy of a function on 
$\Gamma_m$ by 
\begin{align*}
  E_m(u) = \sum_{x \sim y} |u(x)-u(y)|^2.
\end{align*}
The renormalization factor for $\VS_n$ is $2n-1$, so the renormalized
graph energy on $\Gamma_m$ is
\begin{align*}
  \E_m(u) = (2n-1)^{-m} E_m(u),
\end{align*}
and we can define the fractal energy $\E(u) = \lim_{m \to \infty}
\E_m(u)$.  We define $\dom \E$ as the space of continuous functions
with finite energy.

Now we have the tools to define a fractal Laplacian.  In $\dom \E$,
$\E$ extends by the polarization formula to a bilinear form $\E(u, v)$
which defines an inner product in this space.  If $\mu$ is the
standard measure, we can define the Laplacian with a weak formulation:
$\Delta u = f$ if $f$ is continuous, $u \in \dom \E$, and 
\begin{align*}
  \E(u, v) = -\int f v\, d\mu
  \qquad \forall v \in \textstyle\dom_0 \E
\end{align*}
where $\dom_0 \E = \{ v \in \E : v|_{\rm bdry} = 0\}$.  There is also
a pointwise formula (which is proven to be equivalent in
\cite{fractals}) which, for nonboundary points in $\VS_n$ computes
\begin{align*}
  \Delta u(x) = \lim_{m \to \infty} K (4n-3)^m (2n-1)^m \Delta_m u(x), 
\end{align*}
with $K$ a constant, and where $\Delta_m$ is a discrete Laplacian
associated with the graph $\Gamma_m$, defined by
\begin{align*}
  \Delta_m u(x) = \frac 1{\mathop{\text{deg}} x}\sum_{y \sim x}
  (u(y)-u(x)),
  \qquad \text{for $x$ not on the boundary}.
\end{align*}
The Laplacian satisfies the scaling property
\begin{align*}
  \Delta(u \circ F_i) = (4n-3)(2n-1)(\Delta u) \circ F_i 
\end{align*}
and by iteration
\begin{align*}
  \Delta(u \circ F_w) = \big((4n-3)(2n-1)\big)^m(\Delta u) \circ F_w 
\end{align*}
for $F_w = F_{w_1} \circ F_{w_2} \circ \cdots \circ F_{w_m}$.

In this paper, we restrict attention to the Laplacian defined with
Neumann boundary conditions.  The Neumann boundary conditions are
``natural'', in the sense that the weak formulation need only be
modified to allow all $v \in \dom \E$, and the pointwise formulation
is also valid at boundary points. It is also possible to define a
normal derivative $\partial_n u(q_i)$ at boundary point, and the Neumann
condition is $\partial_n u(q_i) = 0$.  Moreover, there are infinitely
many points in $\VS_n$ that have neighborhoods isometric to
neighborhoods of boundary points; the Neumann boundary conditions
treat the boundary points no differently from these equivalent points.
(Note that this is not true on $\SG$.) These are ample reasons to
prefer Neumann to Dirichlet boundary conditions. An additional benefit
is that the theory is considerably simpler.

The Laplacian on a fractal such as $\SG$ or $\VS_n$ has a discrete
spectrum of positive eigenvalues $\lambda_0 < \lambda_1 < \lambda_2 <
\cdots$, which can be computed explicitly by the method of spectral
decimation developed by Fukushima and Shima, and applied to the Vicsek
set in \cite{zhou}.  Spectral decimation is a method of relating
eigenfunctions and eigenvalues from one graph approximation to a finer
one.  In Section~\ref{sec:decimation}, we describe the method and
explicitly compute an algorithm for spectral decimation on $\VS_2$,
which allows us to numerically calculate eigenfunctions on the Vicsek
set, and observe patterns in the data. 

Let $\{\lambda_j\}$ denote the spectrum of the Laplacian, and let
$\{u_j\}$ denote an orthonormal basis of eigenfunctions.  Then for any
bounded function $f$, we can define the spectral operator $f(-\Delta)$
on $L^2(\VS_n)$ by
\begin{align*}
  f(-\Delta) u = \sum_{j = 1}^\infty f(\lambda_j) \langle u, u_j
  \rangle u_j.
\end{align*}
These operators include the fundamental solutions to the heat and wave
equations, and solutions for other space-time equations.   Because of
the importance of spectral operators to classical analysis,
understanding spectral operators and the Laplacian on $\VS$ is a key
goal in the development of analysis on fractals.

In computing a spectral operator, we can group terms in the sum
corresponding to the same eigenvalue, and write
\begin{align*}
  f(-\Delta) u(x) = \sum_{\lambda_j} f(\lambda) \int P_\lambda(x, y)
  u(y) d\mu(y) 
\end{align*}
where, at a given point $x$
\begin{align*}
  P_\lambda(x, y) = \sum_j u_j(x) u_j(y), 
\end{align*}
$\{u_j\} $ being an orthonormal basis of the $\lambda$-eigenspace
$E_\lambda$.  In Section~\ref{sec:boundary} we show how, for certain
special points $x$, we can simplify this sum to a single term.
Fixing a point $x$ on the boundary, or at the center, and letting
$E_\lambda^x$ denote the subspace of $E_\lambda$ of functions
vanishing at $x$, we can choose the orthonormal basis so that the
first element $u_1$ is in $(E_\lambda^x)^\perp$ and the rest belong to
$E_\lambda^x$.  Then, 
\begin{align*}
  P_\lambda(x, y) = u_1(x) u_1(y).
\end{align*}

Additionally, in Section~\ref{sec:boundary} we prove a formula for the
inner product of two eigenfunctions on a graph approximation, and show
that it converges in the limit to the inner product on the Vicsek set.
This ensures that functions which are orthogonal on
graph approximations remain orthogonal on the Vicsek set, and makes it
possible to compute $P_\lambda$ when $x$ is a point on the boundary or
at the center.  Here we follow some of the ideas in
\cite{spectralops}.

In Section~\ref{sec:data}, we give some numerical data using our
MATLAB algorithms for the eigenvalues and eigenfunctions of the
Laplacian on $\VS_2$ and $\VS_n$.  We also give data on the eigenvalue
counting function $N(x)$ and the Weyl ratio $N(x)/x^\alpha$, for the
appropriate power $\alpha$.

In Section~\ref{sec:spectral}, we give numerical results for the heat
kernel, the propagator for the wave equation, and the spectral
projections onto the 0-series.

In Section~\ref{sec:diag}, we show that each 0-series eigenfunction is
determined by its restriction to the diagonal of the Vicsek set. 

In Section~\ref{sec:gaps}, we prove, following \cite{products} the
existence of a ratio gap in the spectrum of the Laplacian.  A ratio
gap is an interval $(a, b)$ such that the ratio of any two eigenvalues
must fall outside the interval; this is a measure of the sparseness of
the spectrum. Related results have been obtained in \cite{ratiogaps}.

In Section~\ref{sec:clusters}, we show the existence of eigenvalue
clusters; that is, arbitrarily many distinct eigenvalues in an
arbitrarily small interval.

In Section~\ref{sec:green}, we calculate an explicit Green's function
for the Laplacian on the Vicsek set.

In Section~\ref{sec:higher}, we examine the convergence of
eigenfunctions and eigenvalues of the Laplacian on $\VS_n$ as $n \to
\infty$ and show that they approach the corresponding values for the
Laplacian on the cross.  

In Section~\ref{sec:weyl} we establish some properties of the Weyl
ratio on $\VS_n$ that begin to explain the curious apparent
convergence to a function that is unrelated to the Weyl ratio on the
cross.

For more data and programs, refer to \url{www.math.cornell.edu/~mhw33}
(\cite{website}).

It is possible to describe $\VS_n$ as the closure of a countable union
of straight line segments; start with the two diagonals, and take all
images under all iterates of $\{ F_i \}$.  (Some images will be proper
subsets of other line segments and should be deleted to eliminate
redundancy.)  We call this the \emph{skeleton} of $\VS_n$, $SK(\VS_n)
= \cup_{j = 1}^\infty I_j$, where the line segments $I_j$ intersect
only at points.  Since the skeleton is dense, any continuous function
is uniquely determined by its restriction to the skeleton, but the
skeleton is not all of the Vicsek set, since it has $\mu$-measure
zero.

Each line segment $I_j$ has a simple one-dimensional energy
\begin{align*}
  \E_j (u, v) = \int_{-b_j}^{b_j} u'(s_j(t)) v'(s_j(t)) dt
\end{align*}
where $s_j: [-b_j, b_j] \to I_j$ is the linear parametrization.  It
is not difficult to see that
\begin{align*}
  \E(u, v) = c \sum_{j = 1}^\infty \E_j (u|_{I_j}, v|_{I_j})
\end{align*}
for the appropriate constant $c$.  From this point of view, the energy
form on $\VS_n$ is trivial.  Because we combine the trivial energy
with the unrelated measure $\mu$, we obtain a nontrivial Laplacian.

On the other hand, there is a natural measure on the skeleton: just
take the sum of Lebesgue measure on each $I_j$.  By the embedding of
the skeleton in $\VS_n$ we may also regard this as a measure $\nu$ on
$\VS_n$.  Of course it is not a finite measure, as the sum of the
lengths of the line segments $I_j$ diverges.  It satisfies the
self-similar identity
\begin{align*}
  \nu = \frac 1{2n-1} \sum_i \nu \circ F_i^{-1}
\end{align*}
in contrast to the self-similar identity
\begin{align*}
  \mu = \frac{1}{4n-3} \sum_i \mu \circ F_i^{-1}
\end{align*}
for $\mu$.  There is good reason to consider $\nu$ as the
\emph{universal energy measure} on $\VS_n$.  If $f \in \dom \E$ then
we may define an associated energy measure $\nu_f$ with $\E(f, f) =
\nu_f (\VS_n)$ and roughly speaking $\nu_f(A)$ is the contribution to
$\E(f, f)$ coming from the set $A$, for any simple set $A$ (for
example, a finite union of cells.)  For each $I_j$ consider the
function $f_j$ defined by $f_j(s_j(t)) = t$ on $I_j$, which is
constant on every other interval that intersects $I_j$.  Then $f_j$ is
harmonic at every point except the endpoints of $I_j$, and $\nu_{f_j}$
is exactly Lebesgue measure on $I_j$.  So 
\begin{align*}
  \nu = \sum_{j = 1}^\infty \nu_{f_j}.
\end{align*}
We can also see that $f = \sum_{j = 1}^\infty f_j$ is a finite sum on
each $I_j$ and $\nu = \nu_f$, although $f$ does not have finite
energy.  One can also show that $\nu_f \ll \nu$ for every function $f
\in \dom \E$.  This is the ``universal'' property of $\nu$.  

On $\SG$ one can define the \emph{Kusuoka measure } $\nu = \nu_{h_1} +
\nu_{h_2}$ where $\{ h_1, h_2 \}$ is an orthonormal basis of global
harmonic functions (modulo constants) in the energy norm, and this
serves as a universal energy measure.  A similar approach would not
work on $\VS_n$, since it would produce a measure supported on the two
diagonals alone.  

It is possible to define an \emph{energy Laplacian} on $\VS_n$ using
the energy $\E$ and the energy measure $\nu$ in place of $\mu$,
although there are some technical problems because $\nu$ is not
finite.  Such a Laplacian would be rather ``trivial'', since it would
amount to the second derivative along each line segment $I_j$,
together with matching conditions on first derivatives at points of
intersection.  We will not consider this Laplacian further in this
paper.  

We hope that this paper makes a strong case that the Vicsek
sets deserve to be considered the simplest nontrivial examples of pcf
self-similar fractals.  There are two sides to this statement. The
first is that the analysis is nontrivial.  Indeed, if you just
restrict attention to harmonic functions on $\VS_n$, the theory is
basically trivial: these are just linear functions on each of the arms
of $\VS_n$ that are constant on all trees that attach to an arm.  But
the graphs we have obtained for eigenfunctions of the Laplacian reveal
that these are nontrivial functions.

The other side of our assertion is that $\VS_n$ is simpler than $\SG$.
The expression for the Green's function and the numerical data for
solutions of the wave equation are good a posteriori evidence for
this.  We can also point to two structural features that can be
considered a priori evidence.  The first is topological: $\VS_n$ is
contractible while $\SG$ has infinite dimensional homology.  Indeed,
the cycles in $\SG$ play a role in the description of the structure of
some of the eigenspaces of the Laplacian (the 5-series in the
terminology of \cite{fractals}.)  The second relates to symmetry:
while $\SG$ only has a 6-element symmetry group, $\VS_n$ has an
infinite symmetry group.  Indeed this group is a semidirect product of
one copy of $S_4$ and infinitely many copies of $S_3$ and $S_2$.
($S_k$ denotes the permutation group on $k$ letters.)  The $S_4$
symmetries are the permutations of the 4 arms, which fix the center
point $q_0$.  For any cell, $F_w V$ with center point $F_w q_0$, with
$w_m \neq 0$, there will be either $S_2$ or $S_3$ symmetries permuting
2 or 3 of the arms of the cell, depending on whether the cell $F_w V$
has 2 or 1 neighboring cells (the permutable arms are the ones with no
neighbors.)

\section{Spectral Decimation}\label{sec:decimation}
The method of spectral decimation was invented by Fukushima and Shima
\cite{fukushima} for $\SG$ to relate eigenfunctions and eigenvalues on
the graph approximations to each other and the eigenfunctions and
eigenvalues on $\SG$. In essence, an eigenfunction on $\Gamma_m$ with
eigenvalue $\lambda_m$ can be extended to an eigenfunction on
$\Gamma_{m+1}$ with eigenvalue $\lambda_{m+1}$, where $\lambda_m =
R(\lambda_{m+1})$ for an explicit functions $R$, except for certain
specified forbidden eigenvalues, and all eigenfunctions on $\SG$ arise
as limits of this process starting at some level $m$. This is true
regardless of the boundary conditions, but if we specify Dirichlet or
Neumann boundary conditions we can describe explicitly all eigenspaces
and their multiplicities. This method was extended to the Vicsek sets
by Zhou \cite{zhou}.

We describe the procedure briefly here.  First, there is a local
extension algorithm that shows how to uniquely extend an eigenfunction
$u$ defined on $V_m$ to a function defined on $V_{m+1}$ such that the
$\lambda$-eigenvalue equations hold on all points of $V_{m + 1}
\setminus V_m$.  Then there is a rational function $R(\lambda)$ such
that if $u$ satisfies a $\lambda_m$-eigenvalue equation on $V_m$, then
the extended function will satisfy the $\lambda_{m+1}$-eigenvalue
equation on $V_{m+1}$ if $\lambda_m = R(\lambda_{m+1})$ and
$\lambda_m$ is not a \emph{forbidden eigenvalue}.  (Forbidden
eigenvalues are singularities of the spectral decimation function $R$.
It is ``forbidden'' to decimate to a forbidden eigenvalue. Because
forbidden eigenvalues have no predecessor --- there is no
$\lambda_{m-1}$ corresponding to $\lambda_m$ --- we speak of forbidden
eigenvalues being ``born'' at a level of approximation $m$.)

We have the following theorem from \cite{zhou}:
\begin{thm}
Define
\begin{align*}
  f_n(\lambda) &= T_n(3 \lambda -1) - 3T_{n -1} (3 \lambda -1)\\
  g_n(\lambda) &= U_{n -1}(3 \lambda -1) - U_{n -2} (3 \lambda -1)\\
  h_n(\lambda) &= U_{n -1}(3\lambda -1) - 3 U_{n-2} (3 \lambda -1)
\end{align*}
where $T_n$ and $U_n$ are the Chebyshev polynomials of the first and
second kind.  Then the spectral decimation function $R$ is
\begin{align*}
  R(\lambda_m) = \lambda_{m-1} = \lambda_m g_m(\lambda_m) h_m(\lambda_m).
\end{align*}
Moreover, the forbidden eigenvalues are 4/3 and the zeroes of $f_n$
and $g_n$.  
\end{thm}
We also have a matrix equation for the eigenfunction extension
formula: If $u|_{V_0}$ is a vector of the values of $u$ on $V_0$ and
$u|_{V_1\setminus V_0}$ is defined analogously, then
\begin{align*}
  u|_{V_1\setminus V_0} = -(X + \lambda_1 M)^{-1} J u_{V_0}.
\end{align*}
where $J$ is the $V_0 \times (V_1 \setminus V_0)$ adjacency matrix, $X$ is the
adjacency matrix for $V_1\setminus V_0$, with the degrees of each
vertex as its diagonal entries, and $M$ is a diagonal matrix
with $M_{ii} = -X_{ii}$. Multiplying this matrix by the values of $u$
on any $k$-cell (with $\lambda_1$ replaced by $\lambda_{k+1}$), we
similarly get the values of $u$ on the $(k+1)$-cells contained in that
$k$-cell.

\begin{figure}
  \centering
  \includegraphics[width=2.5in]{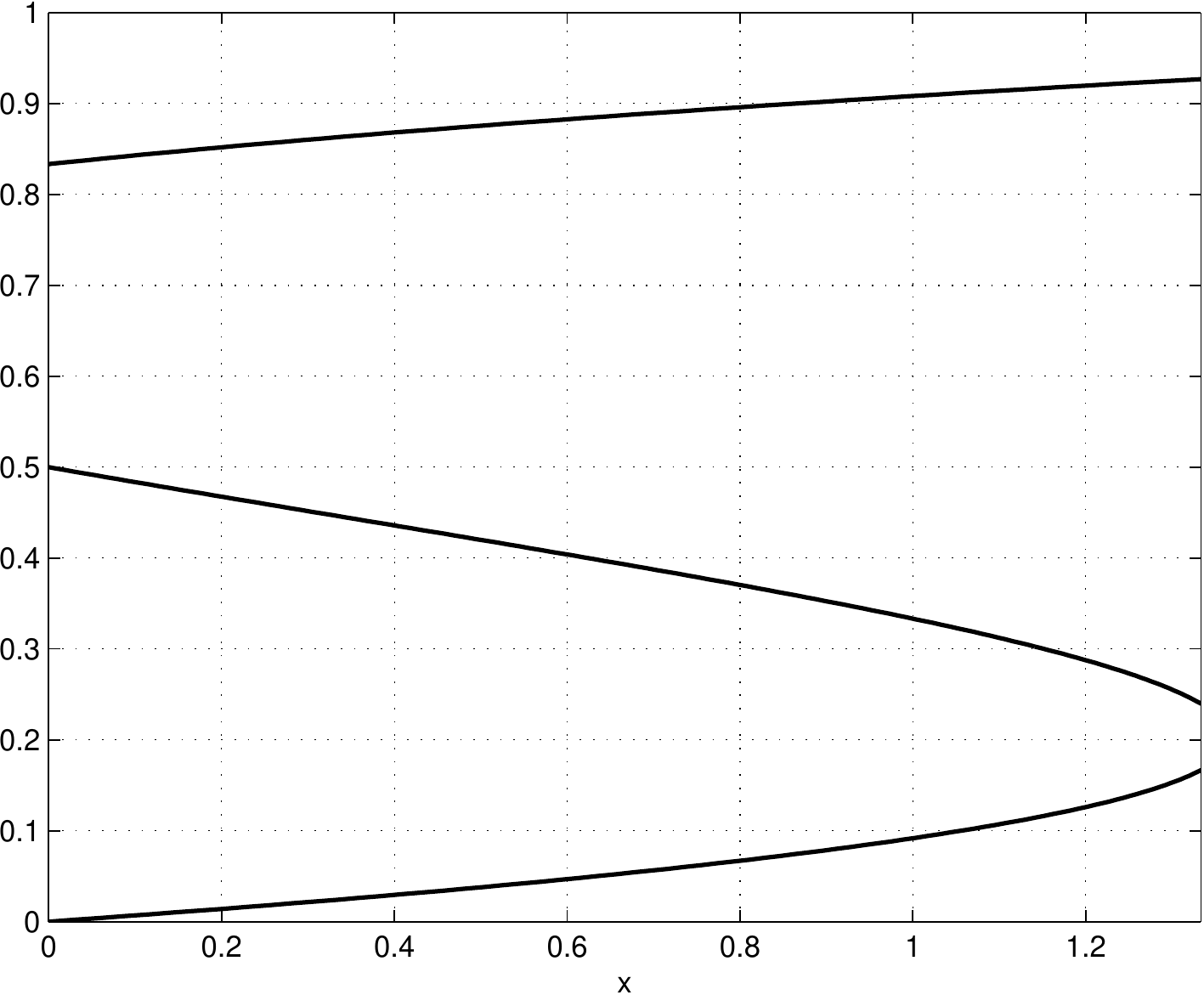}
  \caption{Inverse functions on $\VS_2$. \label{fig:inverses}}
\end{figure}

\begin{figure}
  \centering
  \includegraphics[width=1.5in]{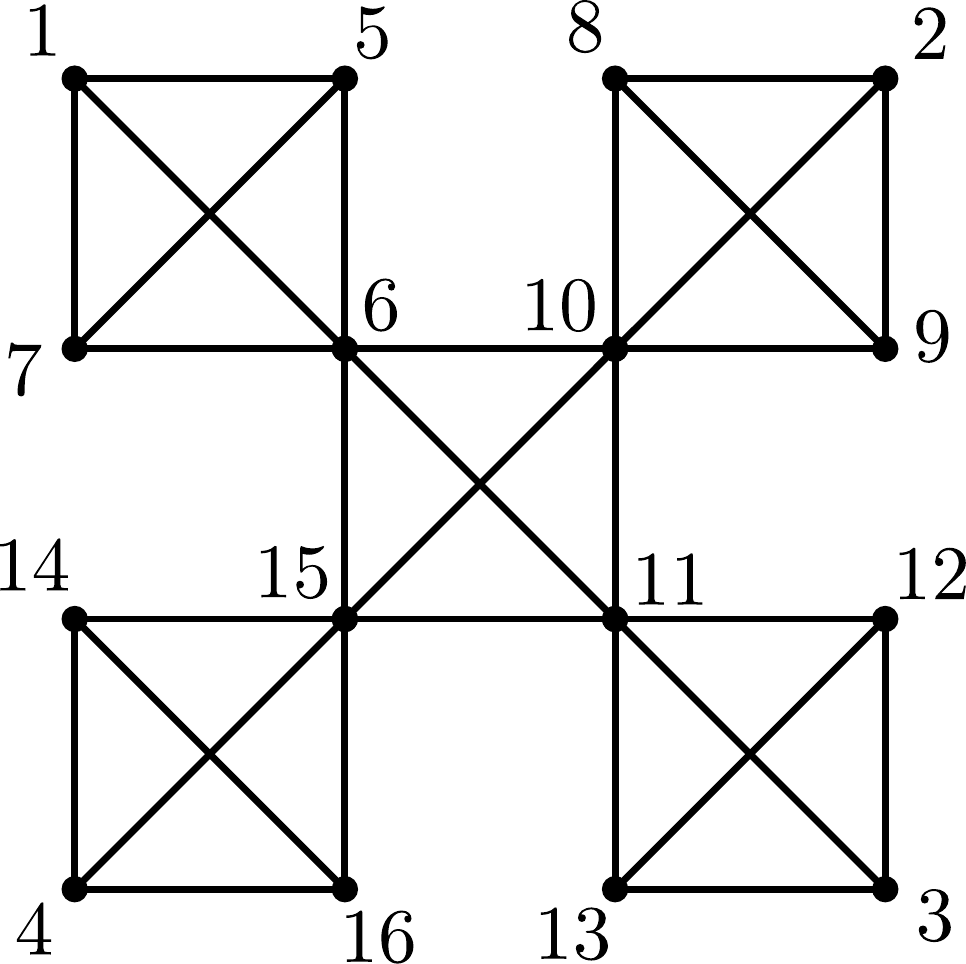}
  \caption{Labeling system for the vertices $V_1$ of the first level
  graph approximation to $\VS_2$.\label{fig:labeling}}
\end{figure}

In the case of $\VS_2$, we have $R(\lambda) = 36 \lambda^3 - 48
\lambda^2 + 15 \lambda$.  The forbidden eigenvalues are $0$, $1/2$,
$4/3$, and $(7\pm \sqrt{17})/12$. There is a 0-eigenvalue born at
level 0, and a 4/3 eigenvalue born at every level thereafter, and
continued eigenvalues are formed by successively choosing one of the
three inverse functions of $R$ (see Figure~\ref{fig:inverses}), so
long as this does not lead to a forbidden eigenvalue.
Using the labeling system described in Figure~\ref{fig:labeling}, the
matrix which allows us to continue eigenfunctions is given by
\begin{align} \label{matrix}
  -(X+\lambda M)^{-1} J = 
  \gamma
  \left(
  \begin{array}{cccccccccccc}
    a&b&a&c&c&d&d&c&c&c&d&c\\
    c&d&c&a&a&b&d&c&c&c&d&c\\
    c&d&c&c&c&d&b&a&a&c&d&c\\
    c&d&c&c&c&d&d&c&c&a&b&a
  \end{array}
  \right)^\top
\end{align}
where 
\begin{align}\label{abcd}
  \begin{aligned}
    a &= 9 - 42 \lambda + 36 \lambda^2, &
    c &= 1, \\
    b &= 6(1 -4 \lambda + 3 \lambda^2), &
    d &= 2 - 3\lambda,\\
    \gamma &= \frac 1{3(4 - 29 \lambda + 60 \lambda^2 - 36 \lambda^3)},
  \end{aligned}
\end{align}
(Note that only roots of $-3/\gamma = (1-2\lambda)f(\lambda)$ are
forbidden eigenvalues, so $\gamma$ is well-defined as long as
$\lambda$ is not forbidden.)

We denote the \emph{4/3-series} as those eigenvalues continued from a
4/3-eigenvalue, and the \emph{0-series} as those eigenvalues continued
from the 0-eigenvalue.  To find $\lambda_m$ from $\lambda_{m-1}$ we
have to invert $R$; in the case of $\VS_2$, there are three inverses,
shown in Figure~\ref{fig:inverses}.  Note that for the sequence $15^m
\lambda_i$ to converge to an eigenvalue $\lambda$ on $\VS$, we need
$\lambda_m$ to approach zero, so we must choose the smallest of the
three inverses all but finitely many times.

A proof in \cite{zhou} guarantees that spectral decimation produces
all possible eigenvalues and eigenfunctions (up to linear
combination), so this formula allows us to explicitly determine the
values of eigenfunctions at arbitrarily high graph approximations.  We
make several observations from numerical calculation of the
eigenfunctions (see Section~\ref{sec:diag}).  One is that the restrictions of
certain eigenfunctions to the diagonal (the segment in $\mathbb{R}^2$
between (0, 0) and (1, 1)) are periodic with period proportional to
$1/m$ and approximate sine functions; this suggests that higher Vicsek
sets $\VS_n$, as they converge to a cross, will have eigenfunctions
approaching the sine and cosine functions in the classical case.  We
will prove this fact in Section~\ref{sec:higher}.

Secondly, we observe that for the 0-series eigenfunctions, choosing
the smallest inverse function of $R$ first means that $\lambda_1 = 0$
so the eigenfunction is extended to be constant on $V_1$.  On each of
the five 1-cells, we start as before, with the eigenfunction having a
value of 1 one all boundary points; so the eigenfunction is
miniaturized into identical copies at each graph approximation, and
the eigenvalue is multiplied by 15.  The same thing happens for any
number of initial choices of the smallest inverse function.

We next describe the structure of the spectrum of the Neumann
Laplacian on $\VS_n$ in complete detail.  Let $\phi_1, \phi_2, \ldots,
\phi_{2n-1}$ denote the inverse functions of the polynomial $R$ in
Theorem 2.1, in increasing order.  We note that $\phi_j$ is an
increasing function when $j$ is odd and is a decreasing function when
$j$ is even.  We write $\rho_n = (4n-3)(2n-1)$ for the Laplacian
renormalization factor.  We write $0 = \lambda_0 < \lambda_1 <
\lambda_2 < \cdots $ for the distinct eigenvalues.  The spectral
decimation rules are summarized as follows:

\begin{enumerate}
\item[(i)] Each eigenvalue has the form
  \begin{align*} 
    &\lim_{m \to \infty} \rho_n^m \phi_{w_m} \circ \phi_{w_{m-1}} \circ
    \cdots \circ \phi_{w_1}(0),\\
    \qquad
    \textup{~or~}
    \qquad
    &\lim_{m \to \infty} \rho_n^{m+k} \phi_{w_m} \circ \phi_{w_{m-1}}
    \circ \cdots \circ \phi_{w_1} (4/3).
  \end{align*}
  where in the first case the eigenvalue is in the 0-series and in the
  second case it is in the 4/3-series and born on level $k$.

\item[(ii)] All but a finite number of the $w_m$ are equal to 1.

\item[(iii)] For the 0-series, the first $w_j$ with $w_j \neq 1$ must
  have $w_j$ an odd number; for the 4/3-series, $w_1$ must be an odd
  number but $w_1 \neq 2n-1$.

\item[(iv)]  The multiplicity of each 0-series eigenvalue is 1, while
  the multiplicity of each 4/3-series eigenvalue born on level $k$ is
  $2(4n-3)^k + 1$.  
\end{enumerate}
Condition (ii) is required in order that the limits in (i) exist.
Let $m_0$ denote the largest value of $m$ for which $w_m \neq 1$ (if
this never happens, let $m_0 = 0$.)  Then we can rewrite the limits in
(i) in terms of a single function $\psi_n$ defined by 
\begin{align*}
  \psi_n(t) = \lim_{m \to \infty} \rho_n^m \phi_1^{(m)} (t)
\end{align*}
(here $\phi_1^{(m)}$ denotes the $m$-fold composition of $\phi_1$).
This limit exists because the Taylor expansion of $R(t)$ about $t = 0$
is $\rho_n t + O(t^2)$, so the Taylor expansion of $\phi_1(t)$ about $t =
0$ is $\rho_n^{-1} t + O(t^2)$.  Then (i) says the eigenvalues are either 
\begin{align*}
  &\rho_n^{m_0} \psi_n(\phi_{w_{m_0}} \circ \phi_{w_{m_0 - 1}} \circ
  \cdots \circ \phi_{w_1}(0)),\\
  \text{~or~} 
  \qquad
  &\rho_n^{m_0+k} \psi_n(\phi_{w_{m_0}} \circ \phi_{w_{m_0}-1} \circ \cdots
  \phi_{w_1}(4/3)) .
\end{align*}
Condition (iii) spells out explicitly the rules for avoiding forbidden
eigenvalues.  We may explain the multiplicities in (iv) as follows.
To satisfy the 4/3-eigenvalue equation on level $k$ we may assign
initial values at the points in $V_k$ so that the sum of the values on
the four boundary points of every $k$-cell is 0. This gives a space of
dimension $\# \{ V_k\} - \# \{ k\text{-cells} \}$ and it is easy to
see that $\#(V_k) = 3 (4n-3)^k + 1$ and $\# \{ k\text{-cells} \} =
(4n-3)^k$. 

\begin{thm}
  Eigenvalues in the 0-series and 4/3-series alternate: $\lambda_j$ is
  0-series for $j$ even and 4/3-series for $j$ odd.  More precisely, the
  spectrum consists of an initial segment of length $2n$ followed by
  segments of length $4n-2$.  In each segment all the 4/3-series
  eigenvalues are born level 0 (hence have multiplicity 3) except the
  last one. 
\end{thm}

\begin{proof}
  Because $\phi_1$ is increasing, so is $\psi_n$, so applying $\psi_n$
  does not change the order of eigenvalues.  Thus the ordering of the
  eigenvalues can be inferred from the ordering and
  increasing/decreasing nature of the $\phi_j$'s.  The lowest
  eigenvalues have $m_0 = 0$, $\lambda_0 = \phi_n(0) = 0$, and
  $\lambda_1 = \psi_n(4/3)$.  After that come those with $m_0 = 1$, in
  the order
  \begin{multline*}
    \rho_n\psi_n( \phi_3(0)) < \rho_n\psi_n( \phi_3(4/3)) < \rho_n\psi_n(
    \phi_5(0)) < \rho_n\phi_n(\phi_5(4/3))
    < \cdots \\ < \rho_n\psi_n( \phi_{2n-1}(0)).
  \end{multline*}
  (here we have used the fact that $\phi_{2j-1}$ is increasing and $\max
  \phi_{2j-1} < \min \phi_{2j + 1}$ and $\phi_{2n + 1}(4/3)$ is a
  forbidden eigenvalue.) Together these give the first $2n-1$
  eigenvalues of the initial segment, with 0-series and 4/3-series born
  at level 0 alternating.  The segment is completed by $\lambda_{2n-1} =
  \rho_n\psi_n(4/3)$, a  4/3-series eigenvalue born on level 1.
  This is valid because $\phi_{2n-1} (0) < 4/3$ (indeed $\max
  \phi_{2n-1} < 1$) and $\rho_n\psi_n(4/3) < \rho_n^2\psi_n( \min
  \phi_2)$.

  Let $\Sigma_1$ denote the sequence 
  \begin{align*}
    \phi_1(4/3), \phi_3(0), \phi_3(4/3), \phi_5(0), \phi_5(4/3),
    \ldots, \phi_{2n-1}(0), 4/3 
  \end{align*}
  and let $\widetilde\Sigma_1$ denote the sequence in reverse order with
  the first and last terms omitted.  Then the initial segment of the
  spectrum has the form $\{ 0, \rho_n\psi_n( \Sigma_1) \}$ (note that
  $\psi_n(4/3) = \rho_n\psi_n((\phi_1(4/3)$).  Let $\Sigma_2$ denote the
  sequence
  \begin{multline*}
    \phi_1(\Sigma_1), \phi_2(\widetilde\Sigma_1), \phi_3(\Sigma_1),
    \phi_4(\widetilde\Sigma_1), \phi_5(\Sigma_1), \ldots,\\
    \phi_{2n-2}(\widetilde\Sigma_1), \phi_{2n-1}(\Sigma_1 \setminus \{4/3 \}), 4/3.  
  \end{multline*}
  Then $\{ 0, \rho_n^2\psi_n(\Sigma_2) \}$ is a larger initial segment
  of the spectrum.  Note that each of the segments 
  \begin{multline*}
    \{\rho_n^2\psi_n(
    \phi_4(\widetilde\Sigma_1)), \rho_n^2\psi_n( \phi_3(\Sigma_1)) \}, \{
    \psi_n (\rho_n^2 \phi_4(\widetilde\Sigma_1)), \psi_n (\rho_n^2
    \phi_5(\Sigma_1) \}, \ldots,\\ \{ \rho_n^2\psi_n( \phi_{2n-2}
    (\widetilde\Sigma_1)), \rho_n^2\psi_n( \phi_{2n-1}
    (\Sigma_1 \setminus \{ 4/3 \})), \rho_n^2\psi_n( 4/3) \} 
  \end{multline*}
  has length $4n-2$ and alternates 0-series and 4/3-series, where all
  except the last 4/3-series eigenvalues are born on level 0. 

  Inductively, we define $\Sigma_k$ to be the sequence
  \begin{align*}
    \phi_1(\Sigma_{k-1}), \phi_2(\widetilde\Sigma_{k-1}),
    \phi_3(\Sigma_{k-1}), \ldots, \phi_{2n-2}(\widetilde\Sigma_{k-1}),
    \phi_{2n-1}(\Sigma_{k-1} \setminus \{4/3 \}), 4/3.  
  \end{align*}
  Then $\{ 0, \rho_n^k\psi_n( \Sigma_k) \}$ is an initial segment of the
  spectrum, and after $\{0, \rho_n\psi_n( \Sigma_1) \}$ it breaks up into
  segments of length $4n-2$ with 0-series and 4/3-series alternating,
  and all but the last 4/3-series alternating, and all but the last
  $4/3$-series eigenvalues are born on level 0.  
\end{proof}

\subsection{Scaling Inner Products}
In order to find an orthonormal basis for eigenspaces, we have to
relate the graph inner product $\langle f, g \rangle_m$ to the inner
product on the next graph approximation, $\langle f, g \rangle_{m+1}$.
This is necessary because we need to compute the inner product
exactly, and we would like to be able to show that functions
orthogonal on one graph approximation will remain orthogonal when
spectrally decimated at higher levels.  We now prove, as
\cite{sampling} does for the Sierpinski Gasket,  a multiplicative
formula for $\langle f, g \rangle_{m + 1}$ in terms of $\langle f, g
\rangle_m$ and the current discrete eigenvalue $\lambda_m$.

\begin{thm}\label{thm:scaling}
  If $u$ and $v$ are eigenfunctions born on level $m'<m$, both with
  the same graph eigenvalues $\lambda_{m-1}$ and $\lambda_m$, then
  \begin{align*}
    \langle u, v \rangle_{m} =  N(m) \langle u, v \rangle_{m -1}
  \end{align*}
  where 
  \begin{align*}
    N(m) = \frac{20 -143 \lambda_m + 240 \lambda_m^2 -108 
    \lambda_m^3}{4 - 28 \lambda_m 
    + 60 \lambda_m^2 - 36 \lambda_m^3}.
  \end{align*}
  The product below converges, and in the limit gives the inner
  product on $\VS_2$ for $u$ and $v$ eigenfunctions born on level 0
  with the same eigenvalue:
  \begin{align*}
    \langle u, v \rangle = \langle u, v \rangle_0 \prod_{m=1}^\infty N(m).
  \end{align*}
\end{thm}

\begin{proof}
  On a graph approximation of the Vicsek Set, we call two points
  \emph{neighbors} if they are connected by an edge.  All points have
  either three or six neighbors.  We define \emph{junction points} to
  be those with six neighbors, and \emph{non-junction points} to be
  those with three neighbors.

  For simplicity we take $u = v$ as the general case is essentially
  the same.  

  The graph inner product of two functions on the graph approximation
  $V_m$  is defined as 
  \begin{align*}
    \langle u, v \rangle_m= \frac 14 \cdot 5^{-m} \sum_{|w| = m} \sum_i 
    u(F_w q_i)\, v(F_w q_i),
  \end{align*}
  where we need to multiply by $\frac 14$ so that $\langle 1, 1\rangle_m
  = 1$.  This makes the limit $\mu$ a probability measure.  Here each
  $w$ is a ``word,'' that is, a string of numbers corresponding to the
  five similarities that define $\VS_2$.  $F_w$ is the composition
  $F_{i_1} \circ F_{i_2} \circ \cdots \circ F_{i_m}$ where $i_j$ are
  the constituents of the word $w$.  

  At each graph approximation, these similarities map two distinct
  points to the junction points, and only one point to the boundary
  points, so we account for double-counting as follows:
  \begin{align*}
    \langle u, v \rangle_m 
    &= \frac 14\cdot 5^{-m} \bigg(2 \sum_{\rm junction} u(x) v(x) 
    + \sum_{\rm nonjunction} u(x)v(x)\bigg)\\
    &= \frac 14\cdot 5^{-m} \sum_x \frac{\deg x}3 \; u(x) v(x).
  \end{align*}

  Fix an $(m-1)$-cell $C$ and let $u_1$, $u_2$, $u_3$, $u_4$ be the values
  of $u$ on its boundary. Then the contribution to $\| u \|_{m-1}^2$
  due to $C$ is 
  \begin{align*}
    \frac 14\cdot 5^{-(m-1)} (u_1^2 + u_2^2 + u_3^2 + u_4^2).
  \end{align*}
  Now we can use spectral decimation (i.e.~\eqref{matrix} with
  $\lambda = \lambda_m$) to get the values of $u$ on the level $m$
  vertices in $C$. Letting $a$, $b$, $c$, $d$, and $\gamma$ as in
  \eqref{abcd} we see that the contribution of $C$ to 
  $\| u \|_m^2$ is 
  \begin{multline*}
    \frac 14\cdot 5^{-m} \Big[(u_1^2+u_2^2+u_3^2+u_4^2) 
    (1+2\gamma^2 (a^2 + b^2 + 3c^2 + 3d^2)) \\
    + 8\gamma^2\sum_{i< j} u_i u_j (ac + bd + c^2 + d^2)\Big].
  \end{multline*}
  Applying this to all $(m-1)$-cells we obtain
  \begin{multline}
    \label{ba}
      \|u\|_m^2 = \tfrac 15 \|u\|_{m-1}^2 
      (1 + 2\gamma^2(a^2 + b^2 + 3c^2 + 3d^2))\\
      + 2\cdot 5^{-m}\gamma^2(ac + bd + c^2 + d^2)\sum_{x\sim y} u(x)u(y),
  \end{multline}
  where the sum is over $V_{m-1}$.  To deal with the cross-terms in
  terms we apply the Gauss-Green formula,
  \begin{align*}
    E_{m-1}(u) &= \sum_{x \sim y} |u(x) - u(y)|^2 \\
    &= -12\cdot 5^{m-1}\langle u, \Delta_{m-1} u \rangle_{m-1} \\
    &= 12 \cdot 5^{m-1}\lambda_{m-1}\| u \|_{m-1}^2.
  \end{align*}
  Since
  \begin{align*}
    \sum_{x\sim y} |u(x) - u(y)|^2 = \sum_x (\deg x)\, u(x)^2 - 
    2 \sum_{x \sim y} u(x) u(y),
  \end{align*}
  this implies
  \begin{align*}
    \sum_{x\sim y} u(x)u(y) 
    = 6 \cdot 5^{m-1}(1 - \lambda_{m-1}) \|u\|_{m-1}^2.
  \end{align*}
  Combining this with \eqref{ba}, we see
  \begin{multline}
    \|u\|_m^2 = \frac 15\Big[
    (1 + 2\gamma^2(a^2 + b^2 + 3c^2 + 3d^2))\\
    + 12\cdot \gamma^2(1 - \lambda_{m-1})(ac + bd + c^2 + d^2)
\Big] \|u\|_{m-1}^2
  \end{multline}
  Simplifying using the values for $\gamma$, $a$, $b$, $c$, $d$, and
  $\lambda_{m-1}$ in terms of $\lambda_m$, we get the normalization
  formula
  \begin{align*}
    N(m) = \frac 15 \cdot \frac{20 -143 \lambda_m + 240 \lambda_m^2 -108 
    \lambda_m^3}{4 - 28 \lambda_m 
    + 60 \lambda_m^2 - 36 \lambda_m^3}.
  \end{align*}
  This allows us to compute the norm of a function the Vicsek set at
  any graph approximation, and, in the limit, on the Vicsek set
  itself:
  \begin{align*}
    \|u\|^2 = \|u\|_0^2 \cdot \prod _{m = 1}^\infty N(m).
  \end{align*}
\end{proof}

\subsection{Center values} \label{centervalues}
It is also useful to have a formula for the value of an eigenfunction
at the center $q_0$ of $\VS_2$. Using \eqref{matrix} to continue a
function $u$ on $V_0$ to $V_1$, we see that the values $u_6$, $u_{10}$,
$u_{11}$ and $u_{15}$ are related to the values of $u$ on $V_0$ by 
\begin{align*}
  u_6 + u_{10} + u_{11} + u_{15} 
  &= \gamma ( 3d+b ) (u_1+u_2+u_3+u_4).
\end{align*}
Substituting for $d$, $b$, and $\gamma$ we get
\begin{align*}
  u_6 + u_{10} + u_{11} + u_{15} 
  &= \frac{4-3\lambda_1}{4-21\lambda_1+18\lambda_1^2} (u_1+u_2+u_3+u_4).
\end{align*}
Continuing this process we get
\begin{align*}
  u(q_0) = \frac{u_1+u_2+u_3+u_4}4  \prod_{m=1}^\infty N'(m)
\end{align*}
where 
\begin{align*}
  N'(m) = \frac{4-3\lambda_m}{4-21\lambda_m+18\lambda_m^2}.
\end{align*}
In particular, since any 4/3-series eigenfunction satisfies
$u_1+u_2+u_3+u_4=0$, all 4/3-series eigenfunctions vanish at $q_0$.

\section{Spectral projections at boundary points}\label{sec:boundary}
We would like to be able to solve differential equations such as the
wave equation
\begin{align*}
  \frac{\partial^2u}{\partial t^2} = \Delta u 
\end{align*}
and the heat equation
\begin{align*}
  \frac{\partial u}{\partial t} = \Delta u 
\end{align*}
on the Vicsek set with Neumann boundary conditions and suitable
initial conditions.  These equations are solved in terms of an
orthonormal basis of eigenfunctions $u_j$: we have
\begin{align*}
  u(x, t) = \sum S(t, j) \left(\int f(y) u_j(y) d\mu(y)\right) u_j(x) 
\end{align*}
where $S$ depends on the equation and the expression in parentheses is
a Fourier coefficient.  Usually the sum and integral can be
interchanged to yield 
\begin{align*}
  u(x,t) = \int K_t(x,y) f(y) d\mu(y),
\end{align*}
where $f$ is defined by the initial conditions and
\begin{align*}
  K_t(x,y) = \sum S(j, t) u_j(x) u_j(y), 
\end{align*}
for instance
\begin{align*}
  \sum e^{-t \lambda_j} u_j(x) u_j(y) 
\end{align*}
for the heat equation and
\begin{align*}
  \sum \frac{\sin \sqrt{\lambda_j} t}
  {\sqrt{\lambda_j}} u_j(x) u_j(y) 
\end{align*}
for the wave equation.  

We can get a better understanding of projection kernels $K_t(x,y)$
when we restrict one of the variables to specific boundary points.
Suppose $y = q_i$, $i = 1, 2, 3, 4$, and note that $E_\lambda^0 = \{ u
\in E_\lambda: u(q_i) = 0 \}$ has codimension 1.  We can compute a
normalized function defined to be perpendicular to this space
$u_0^\lambda \in (E_\lambda^0)^\perp.$  In that case we can simplify
\begin{align*}
  P_\lambda (q_i, x) = \sum_{j=0}^N u_j^\lambda (q_i) u_j^\lambda(x) =
  u_0^\lambda(q_i) u_0^\lambda(x).
\end{align*}

If $\lambda$ is a 4/3-series eigenvalue born on some level $m_0 \ge
1$, there is an easy characterization of $u_0^\lambda$; spectral
decimation works ``in reverse,'' i.e.~$u_0^\lambda$ is an
eigenfunction of $\Delta_{m_0-1}$ with eigenvalue $R(4/3) = 20$.
We can then continue spectral decimation back to all 
levels $< m_0$ because we will never encounter a forbidden eigenvalue. 

\begin{thm}[cf. \cite{spectralops} Theorem 3.2, \cite{teplyaev}
  Theorem 3.6]
  Let $\lambda$ be a 4/3-series eigenvalue with $m_0 \ge 1$. Then
  \begin{align*}
    \Delta_{m_0-1} u_0^\lambda = 20u_0^\lambda,
  \end{align*}
  i.e.~$u_0^\lambda$ is an eigenfunction of $\Delta_{m_0-1}$ with
  eigenvalue $20$.
\end{thm}
\begin{proof}
  Fix a point $x \in V_{m_0-1}$. First assume $x$ is only part of a
  single $1$-cell in $\Gamma_{m_0}$, and let $y_1,y_2,y_3$ be the
  other boundary points of that cell. Then the function $v$ shown in
  Figure~\ref{fig:perp}a is a 4/3-series eigenfunction born on level
  $m$ (this is easy to see since the sum around any small square is
  0).
  \begin{figure}
    \centering
    \includegraphics[width=4.7in]{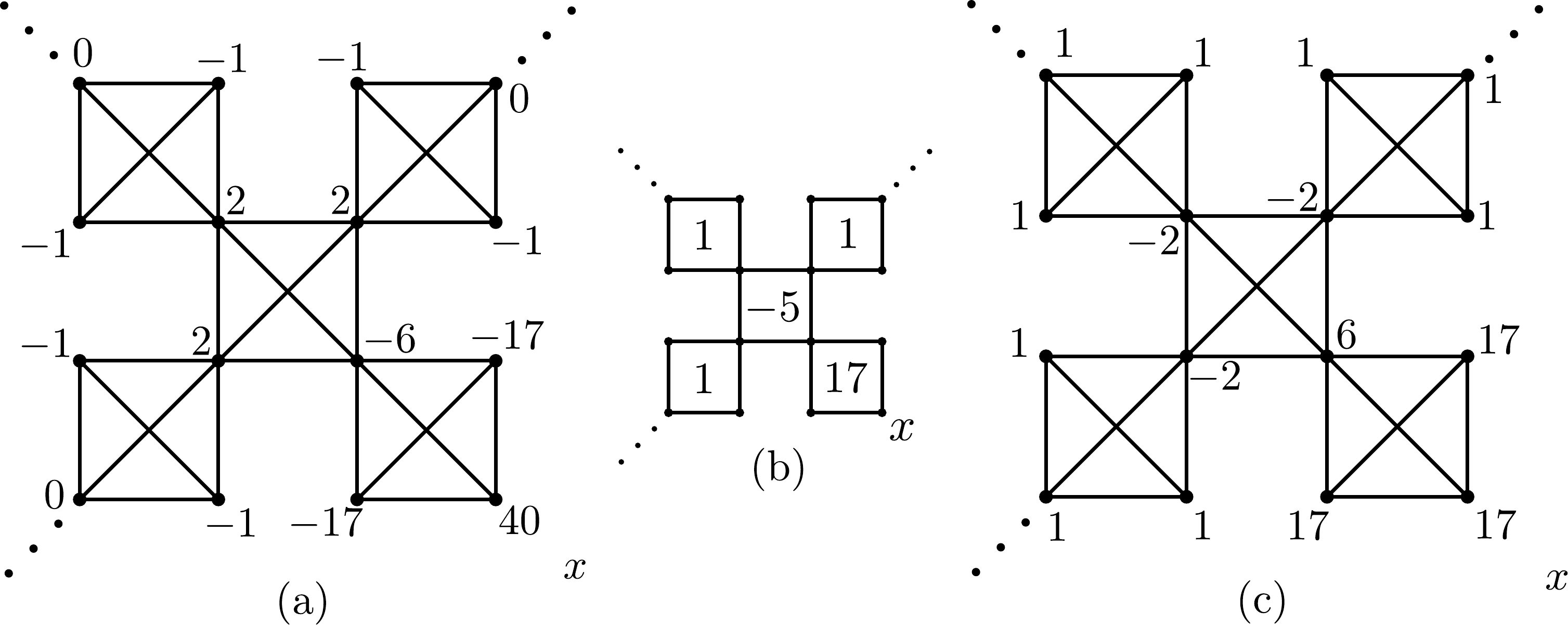}
    \caption{Functions on $\Gamma_{m_0}$, 0 outside of the 1-cell
    containing $x$: (a) $v$.
    (b) Coefficients in the linear combination defining $w$.
    (c) $w$.
    \label{fig:perp}}
  \end{figure}
  Now $v$ vanishes on $V_{m_0-1}\setminus \{x\}$, so in particular
  $v(q_i) = 0$, which forces $\langle u_0^\lambda, v \rangle
  = 0$ and hence (by Theorem~\ref{thm:scaling}) $\langle u_0^\lambda,
  v \rangle_{m_0} = 0$.

  Since $u_0^\lambda$ is a 4/3-series eigenfunction born on level
  $m_0$, we also know that the sum around any small square in
  $\Gamma_{m_0}$ is 0. Taking a linear combination of these equations,
  with coefficients given by Figure~\ref{fig:perp}b, and recalling
  that the inner product weights the center 4 vertices by 2, we see
  that $\langle u_0^\lambda, w \rangle_{m_0} = 0$ where $w$ is given
  by Figure~\ref{fig:perp}c.  Writing $\langle u,v+w\rangle_{m_0} = 0$
  we get
  \begin{align*}
    57u_0^\lambda(x) + u_0^\lambda(y_1) + u_0^\lambda(y_2) +
    u_0^\lambda(y_3) = 0 
  \end{align*}
  or equivalently
  \begin{align*}
    20u_0^\lambda(x) - \frac{u_0^\lambda(y_1) + u_0^\lambda(y_2) +
    u_0^\lambda(y_3)}3 = 20u_0^\lambda(x)
  \end{align*}
  i.e.~$\Delta_{m_0-1} u_0^\lambda(x) = 20 u_0^\lambda(x)$. 
  
  A similar argument works when $x$ is instead part of two 1-cells in
  $\Gamma_{m_0-1}$, with the function on in Figure~\ref{fig:perpb}a
  playing the role of $v$ and the one in Figure~\ref{fig:perpb}b
  playing the role of $w$.
  \begin{figure}
    \centering
    \includegraphics[width=4.7in]{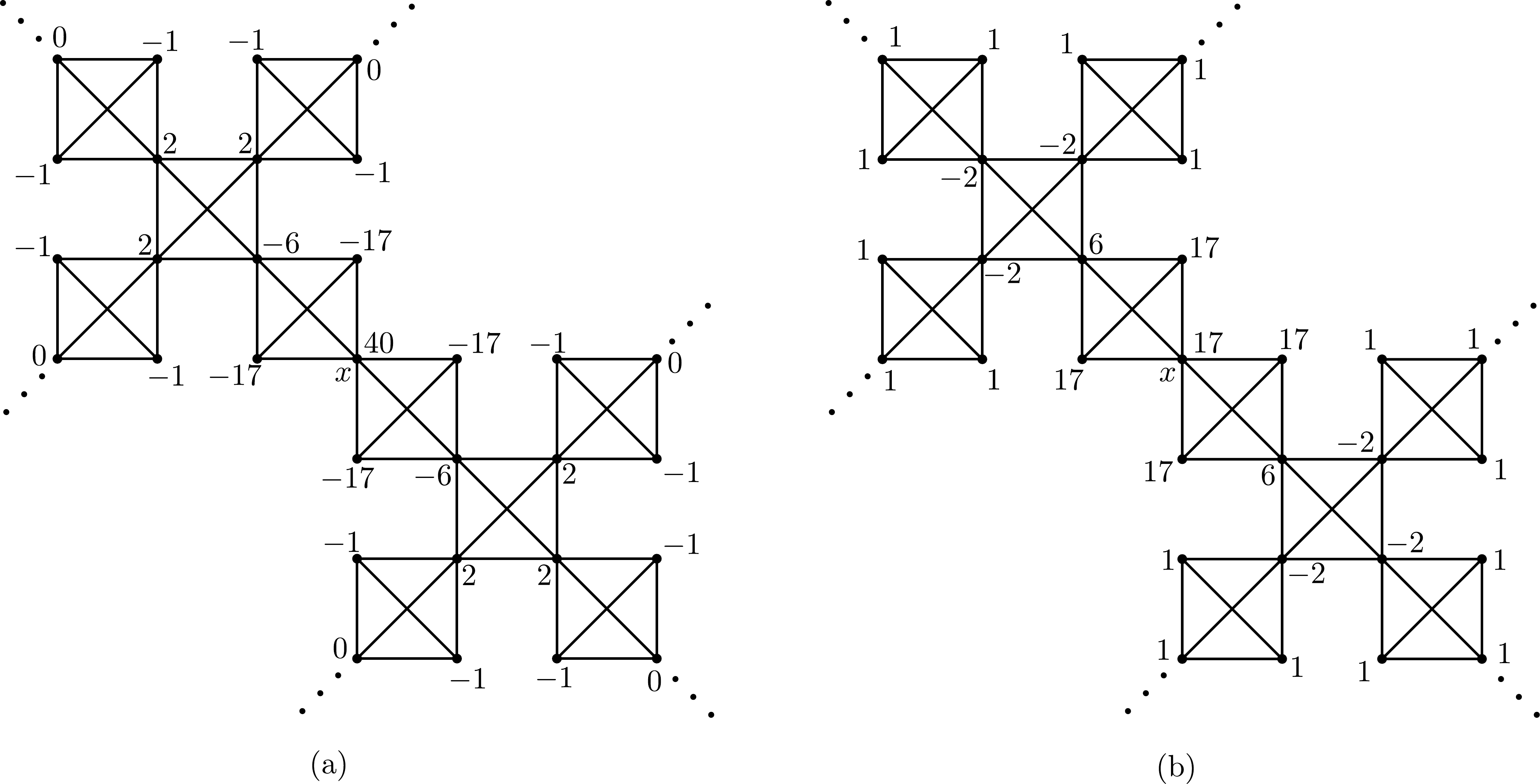}
    \caption{Functions on $\Gamma_{m_0}$, 0 outside of the two
    1-cells
    containing $x$: (a) $v$.  (b) $w$.
    \label{fig:perpb}}
   \end{figure}
\end{proof}
Another special case occurs when we fix $y = q_0$, where $q_0$ is the
center point of the Vicsek set.  At $q_0$, all the eigenfunctions
associated with 4/3-series eigenvalues are equal to zero (see Section
\ref{centervalues}).  This is a fortunate result because, in
calculating the projection kernel at $q_0$, all the terms from the
4/3-series contribute zero, so we only need to consider the
eigenfunctions associated with 0-series eigenvalues --- and these form
a one-dimensional vector space.

\section{Numerical data for eigenvalues and
eigenfunctions}\label{sec:data}

\begin{figure}
  \centering
  \includegraphics[width=5in]{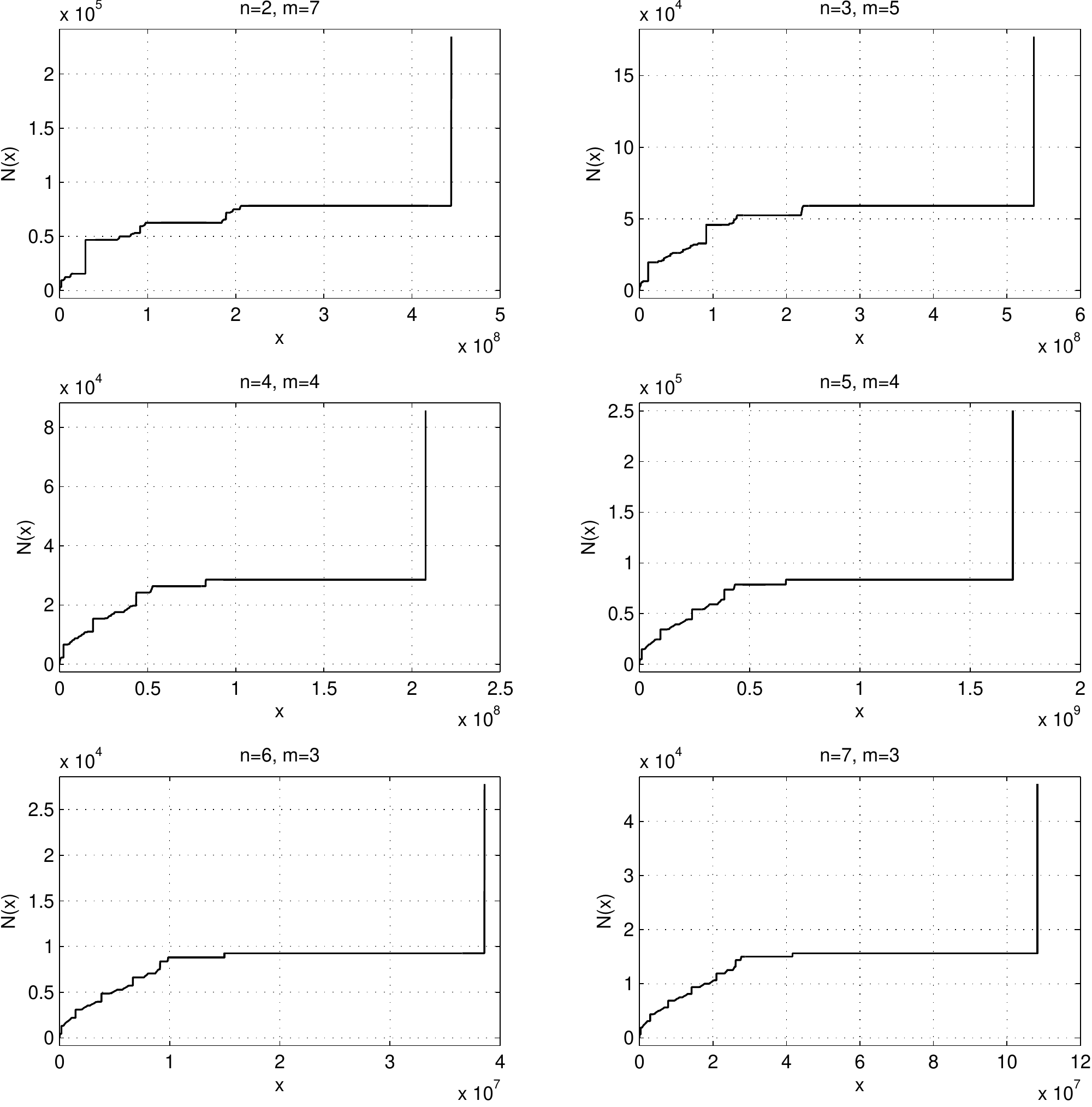}
  \caption{Eigenvalue counting function on $\VS_n$ for $n=2,\ldots,7$.
  Here $m$ refers to the level of the graph approximation.
  \label{fig:counting}}
\end{figure}

\begin{figure}
  \centering
  \includegraphics[width=5in]{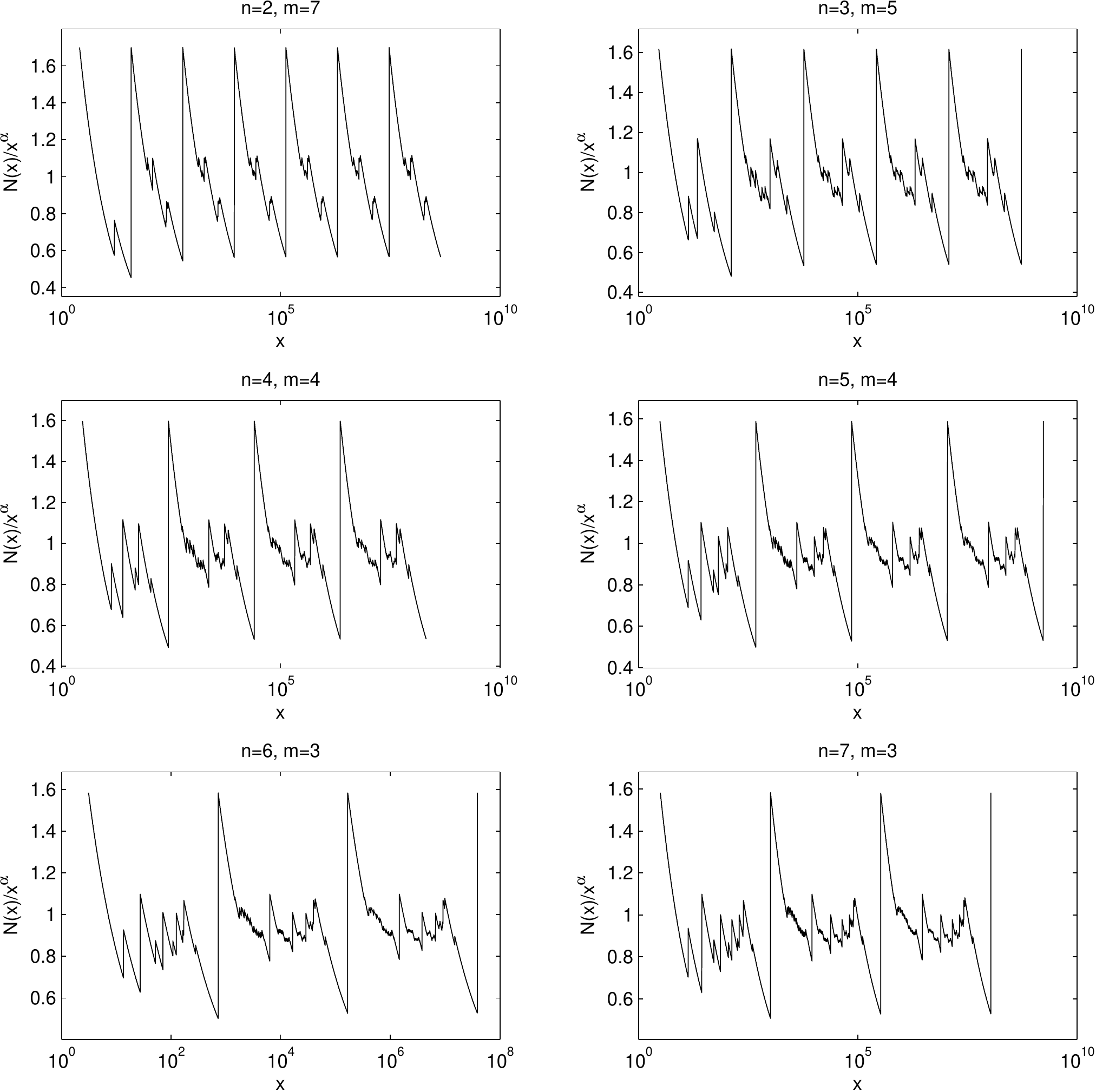}
  \caption{Weyl ratios on $\VS_n$ for $n=2,\ldots,7$. Here $m$ refers
  to the level of the graph approximation.  \label{fig:weyl}}
\end{figure}

Using our implementation of spectral decimation on $\VS_n$, we can
compute the eigenvalues of the graph Laplacians $\Delta_m$ on the
graph approximations $\Gamma_m$. By repeatedly applying the smallest
inverse of the spectral decimation function $R$, we can use these to
approximate the eigenvalues $\lambda_i$ of the standard Laplacian
$\Delta$. Figure~\ref{fig:counting} shows plots of the eigenvalue
counting function $N(x) = \#\{i : \lambda_i \le x\}$.
Since the eigenvalue counting function grows as $x^\alpha$ where
$\alpha = \log(4n-3)/\log((4n-3)(2n-1))$, it is also useful to look at
the Weyl ratio $N(x)/x^\alpha$, shown in Figure~\ref{fig:weyl}.
For each $n$, these functions are asymptotically periodic as a
function of $\log x$, as predicted in \cite{lapidus}. What is
rather striking and somewhat mysterious, there appears to be a
convergence as $n \to \infty$, after appropriate rescaling. We will
attempt to explain some of this in Section~\ref{sec:weyl}.

\begin{figure}[hp]
  \centering
  \includegraphics[width=5in]{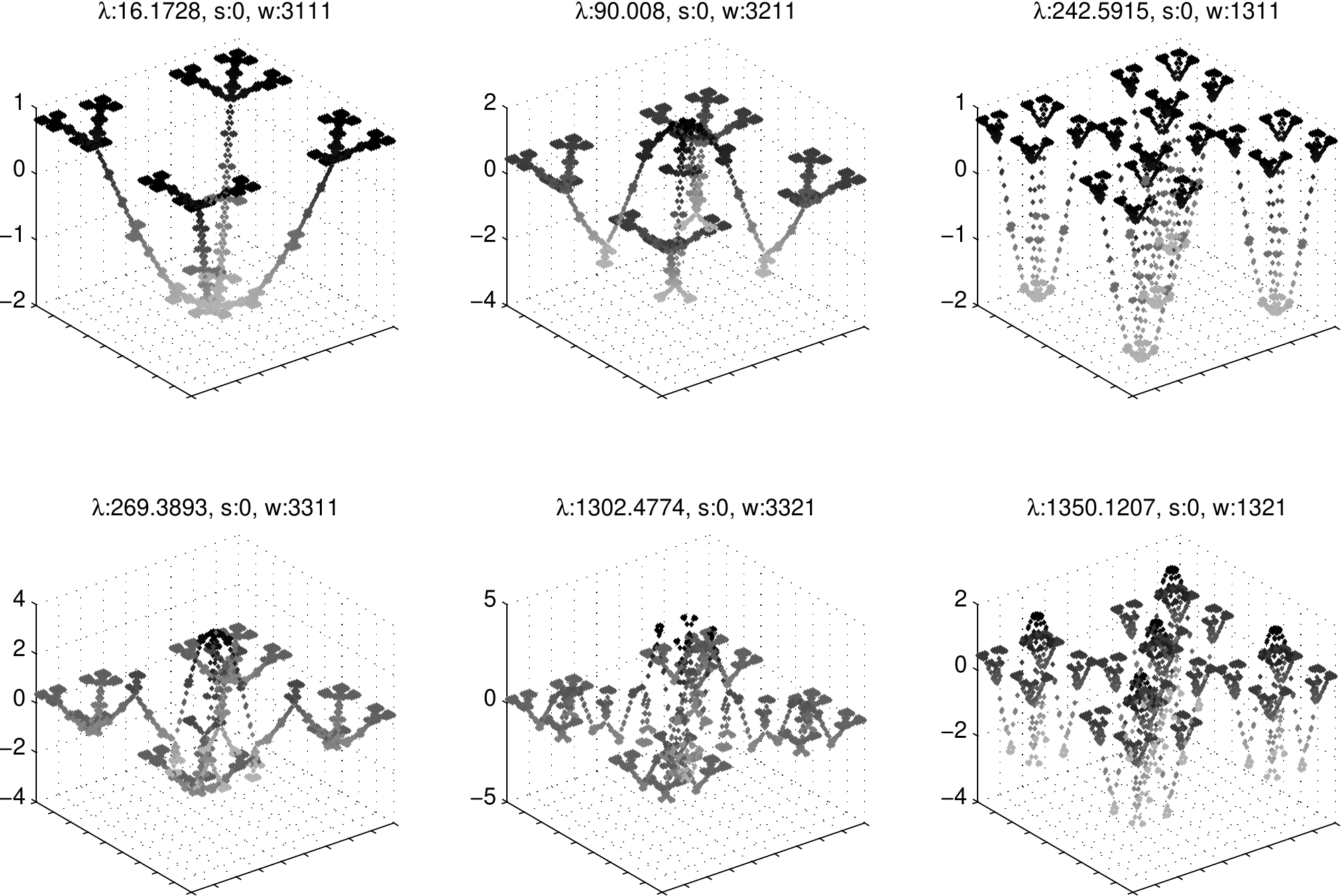}\vspace{8ex}
  \includegraphics[width=5in]{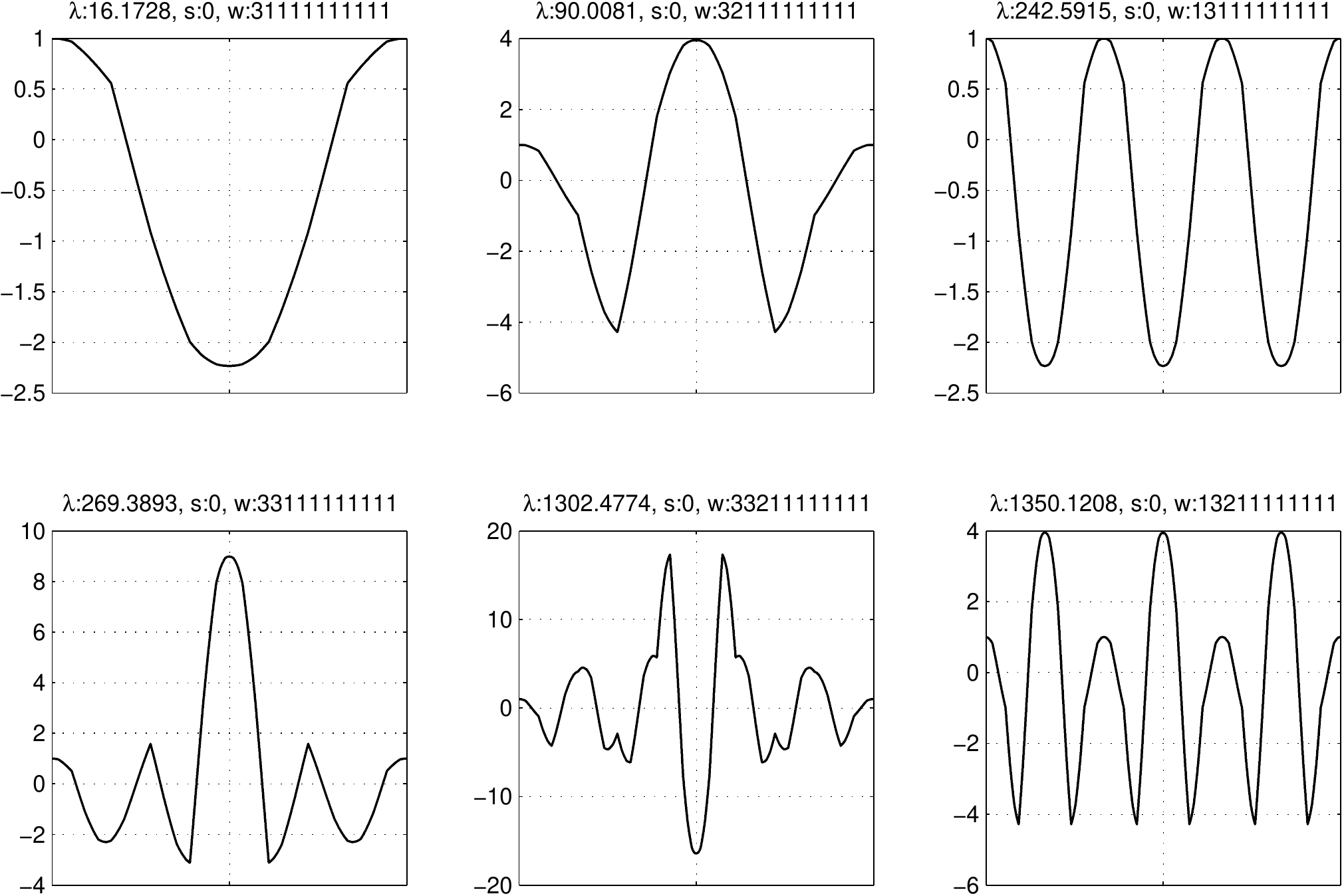}
  \caption{Some 0-series eigenfunctions on $\VS_2$ and their
  restrictions to the diagonal. \label{fig:0-series}}
\end{figure}

\begin{figure}[hp]
  \centering
  \includegraphics[width=5in]{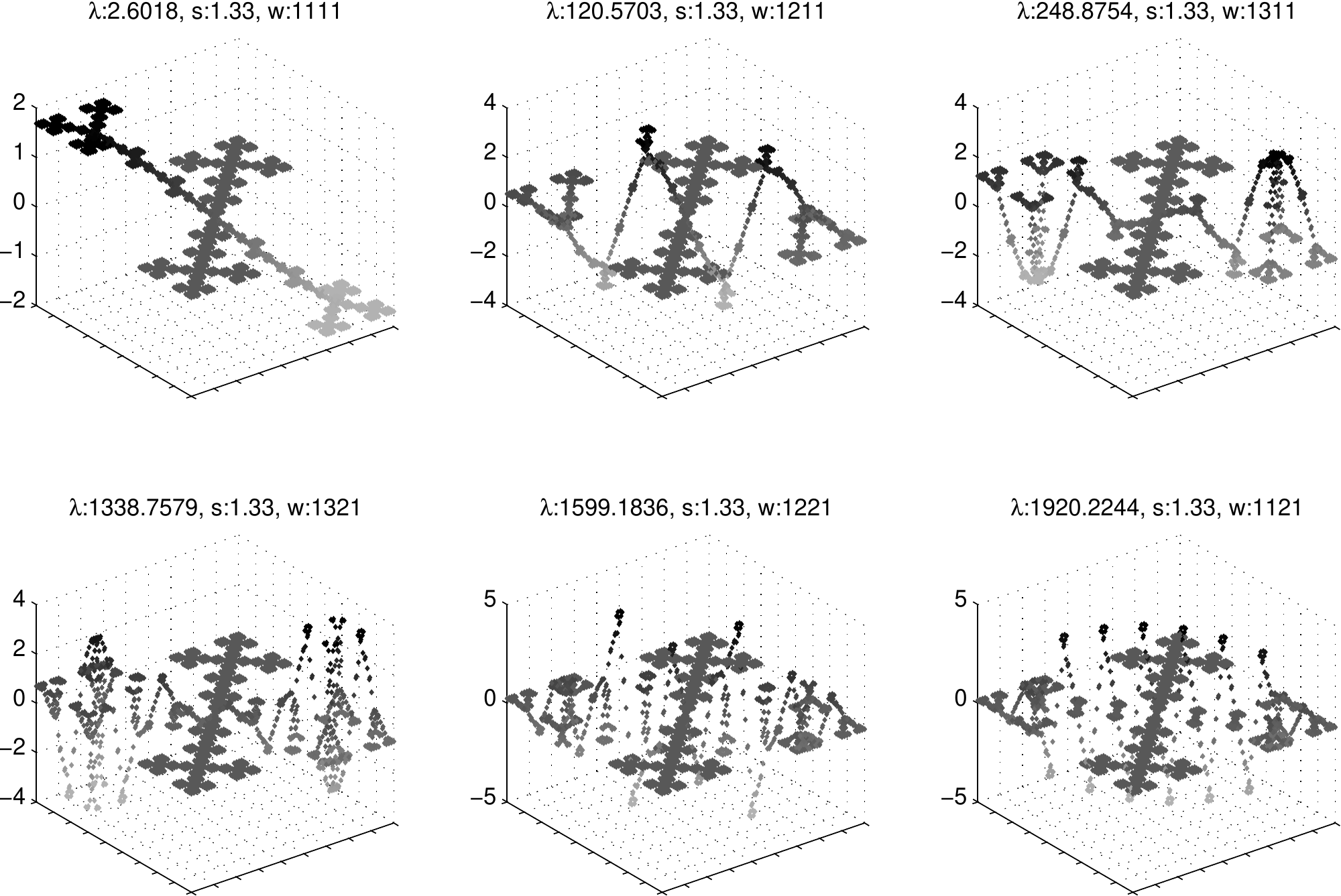}\vspace{8ex}
  \includegraphics[width=5in]{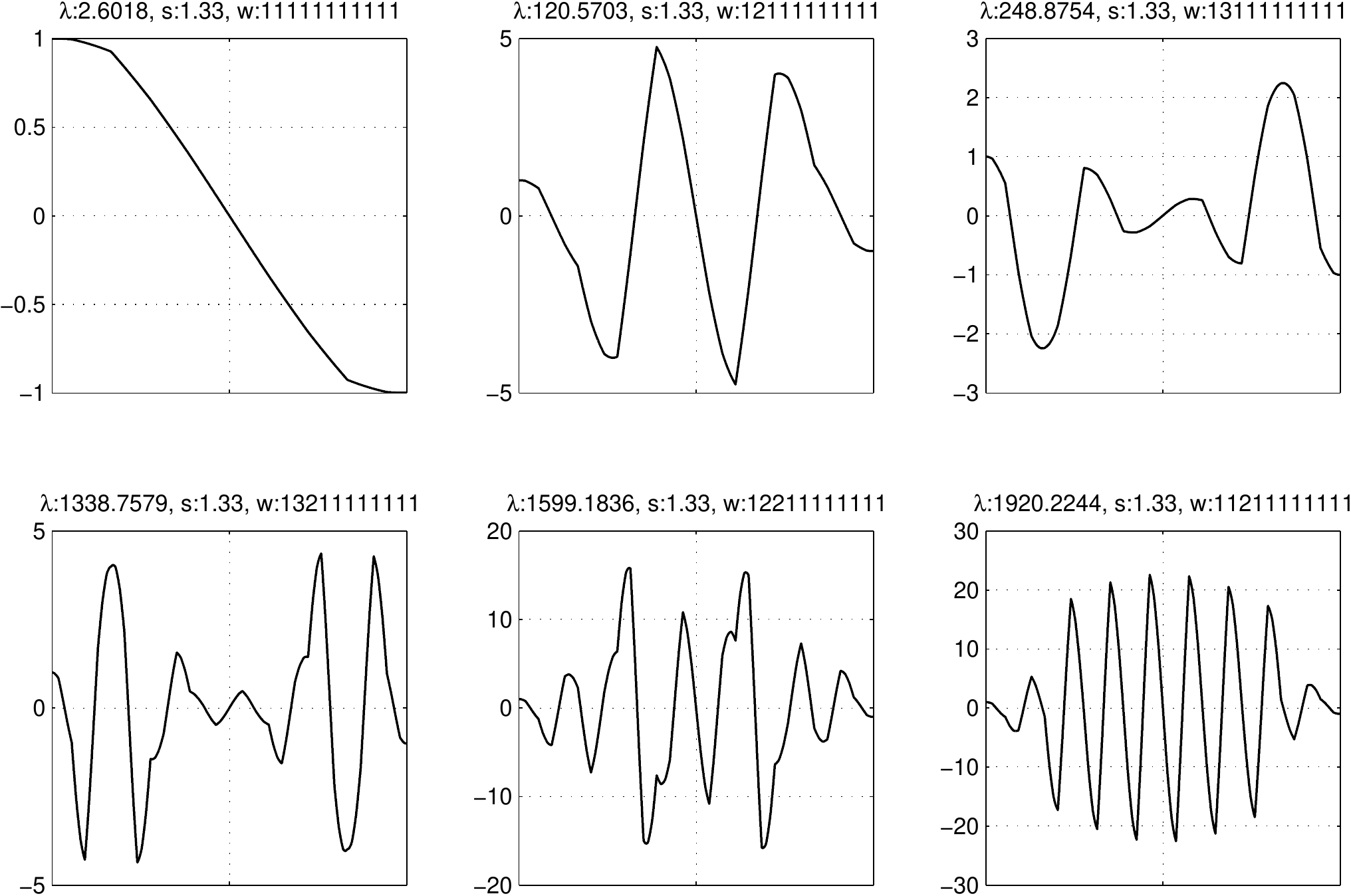}
  \caption{Some 4/3-series eigenfunctions on $\VS_2$ and their
  restrictions to the diagonal. \label{fig:4/3-series}}
\end{figure}

We can also compute eigenfunctions of the graph Laplacians. Figure
\ref{fig:0-series} shows 0-series eigenfunctions and their
restrictions to the diagonal, and Figure \ref{fig:4/3-series} shows
the same for some 4/3-series eigenfunctions. The eigenfunctions in the
diagonal plots have been continued with the lowest inverse several
times to increase the number of data points.
\begin{figure}[hp]
  \centering
  \includegraphics[width=5in]{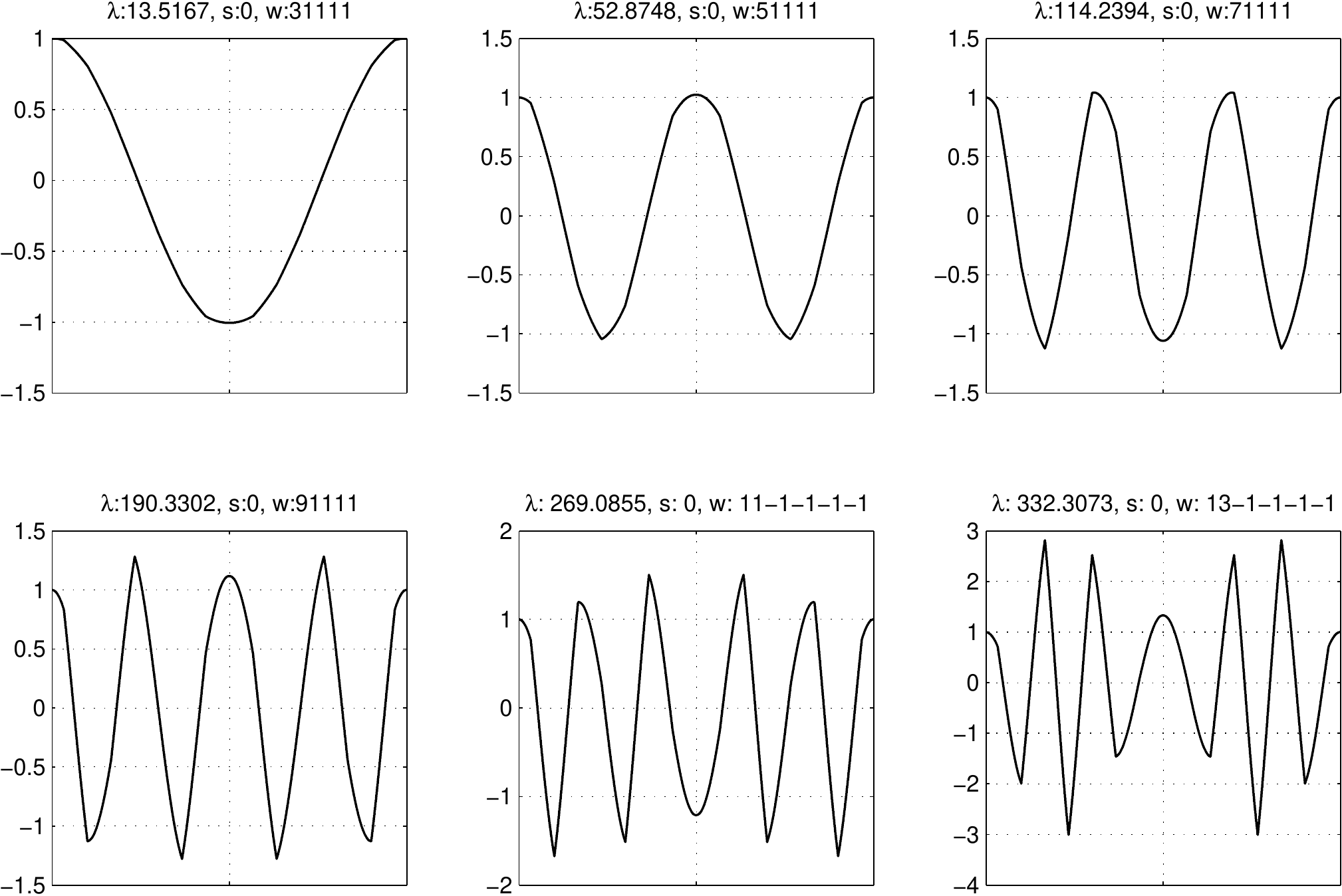}
  \caption{Some 0-series eigenfunctions on $\VS_8$, restricted to the
  diagonal. \label{fig:0diag8}}
\end{figure}
\begin{figure}[hp]
  \centering
  \includegraphics[width=5in]{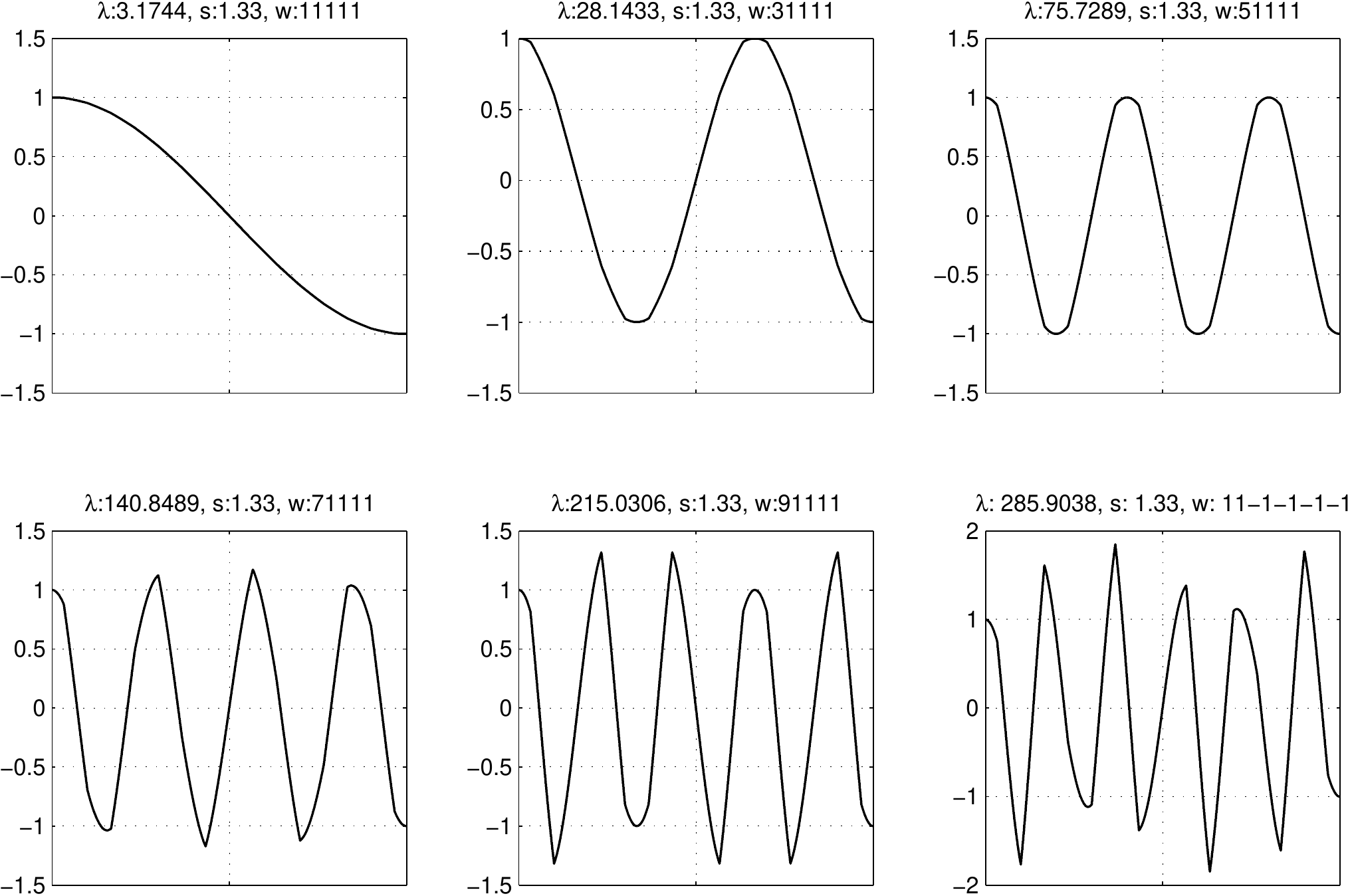}
  \caption{Some 4/3-series eigenfunctions on $\VS_8$, restricted to the
  diagonal. \label{fig:1diag8}}
\end{figure}
For $n > 2$, our implementation can only compute eigenfunctions
restricted to the diagonal. ./figures \ref{fig:0diag8} and
\ref{fig:1diag8} show these plots for $\VS_8$.
There is more data on the website \cite{website}.

We observe from the data a phenomenon known as miniaturization
\cite{outerapprox}.  Taking a 0-series eigenfunction on the $m$th
level approximation to $\VS_2$, if the function is continued by
spectral decimation to the $(m + 1)$th level of approximation, the new
eigenfunction is composed of 5 copies of the previous one; it is
``miniaturized.'' Thus, eigenfunctions of higher eigenvalue are
composed of copies of eigenfunctions of lower eigenvalue.

\section{Spectral operators}\label{sec:spectral}

We can apply the theory of spectral operators to finding solutions of
differential equations on the Vicsek set.  Two major equations are the
heat equation and the wave equation.

\subsection{Heat Kernel}
The heat equation for a function $u(x,t)$ with Neumann boundary
conditions states
\begin{align*}
  \left\{
  \begin{aligned}
    \frac{\partial u}{\partial t} &= \Delta_x u(x, t), \\
    \partial_n u(q_j,t) &= 0,\\
    u(x, 0) &= f(x).
  \end{aligned}
  \right.
\end{align*}
Formally this is solved by
\begin{align*}
  u(x, t) = e^{t \Delta} f(x)
\end{align*}
and since the Laplacian has a discrete spectrum with an orthonormal
basis $\{ u_j \}$ of eigenfunctions, $-\Delta u_j = \lambda_j u_j$,
the solution to the heat equation is
\begin{align*}
  u(x, t) = \sum_{j} e^{-t \lambda_j} 
  \left(\int f(y) u_j(y) d\mu(y) \right) u_j(x).
\end{align*}
Usually the sum and integral can be interchanged to yield
\begin{align*}
  u(x, t) = \int h(t, x, y) f(y) d\mu(y)
\end{align*}
where $h$ is defined to be 
\begin{align*}
  h(t, x, y) = \sum_j e^{-t \lambda_j} u_j(x) u_j(y)
\end{align*}
and called the heat kernel.

\begin{figure}
  \centering
  \includegraphics[width=5in]{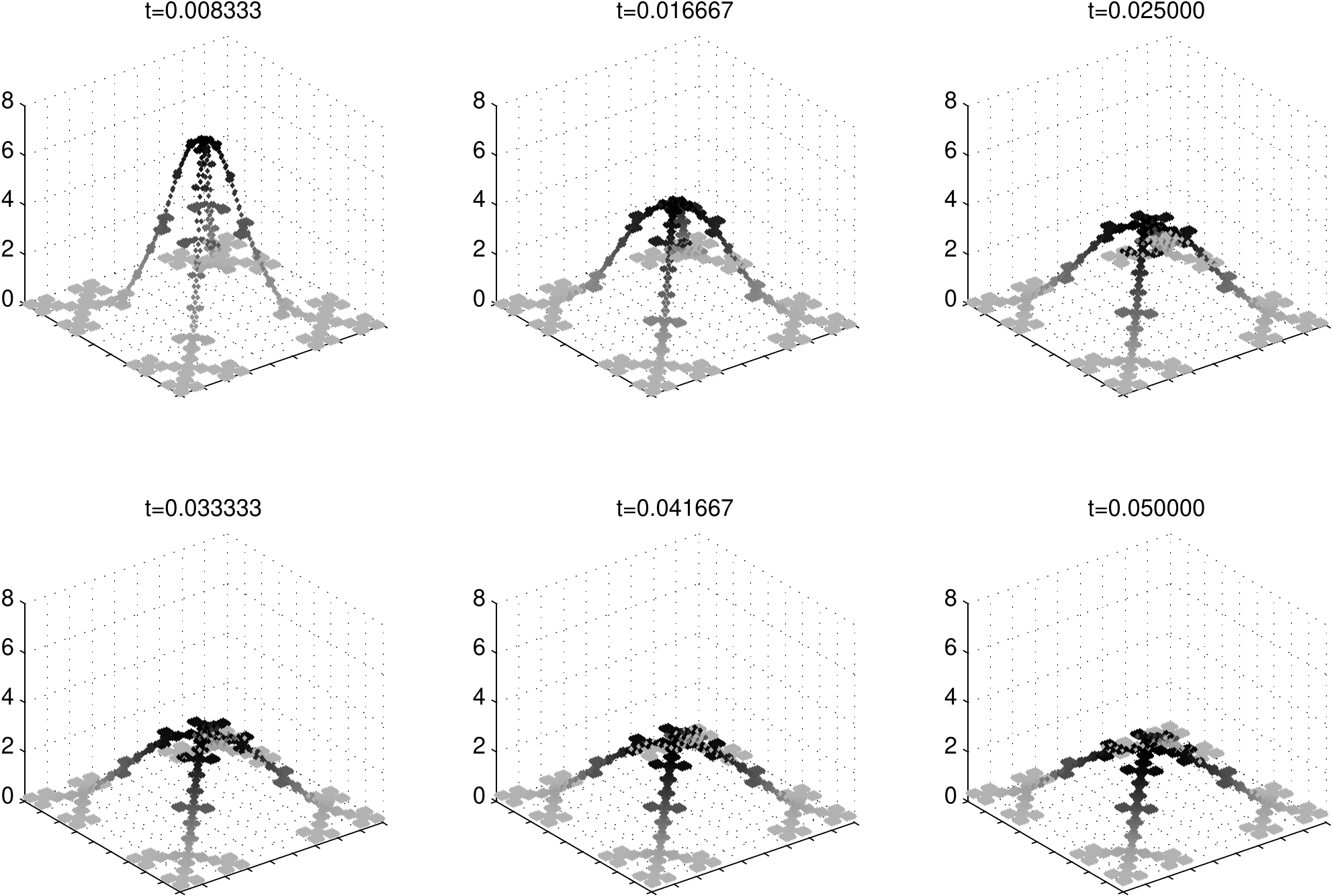}\vspace{10ex}
  \includegraphics[width=5in]{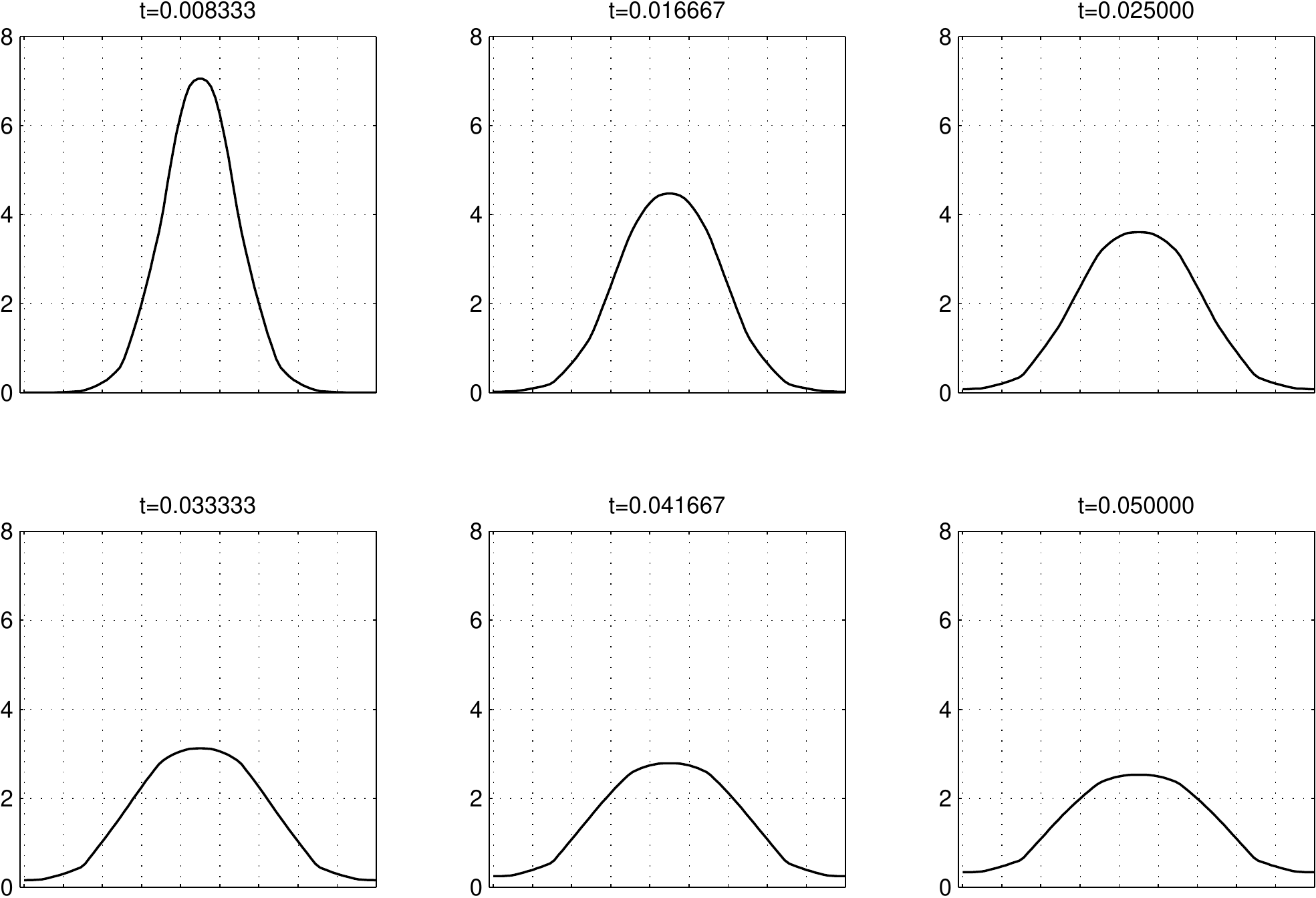}
  \caption{The heat kernel $h(t,q_0,x)$ on $\VS_2$ for various values
  of $t$, on the whole Vicsek set and restricted to the diagonal
  \label{fig:heat}}
\end{figure}

From the eigenvalues and eigenfunctions we can construct the heat
kernel on the standard Vicsek set. This is especially easy when one of
the arguments is the center point of the Vicsek set, since then we
only need to consider 0-series eigenfunctions. Plots of the heat
kernel on the $m=4$ approximating graph are shown in
Figure~\ref{fig:heat}.

\begin{figure}[h]
  \centering
  \includegraphics[width=5in]{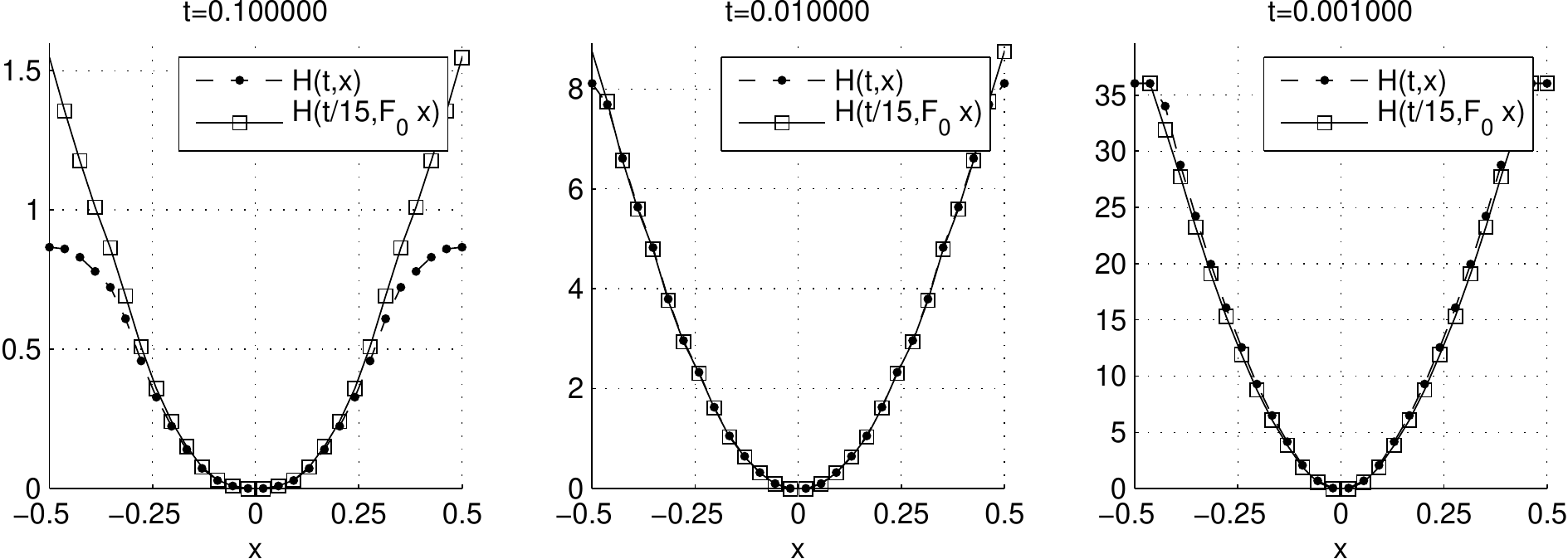}
  \caption{$H(t,x)$ and $H(t/15,F_0x)$ on the diagonal (using the
  $m=4$ graph approximation).
  \label{fig:heatscaling}}
\end{figure}

\begin{figure}
  \centering
  \includegraphics[width=2.5in]{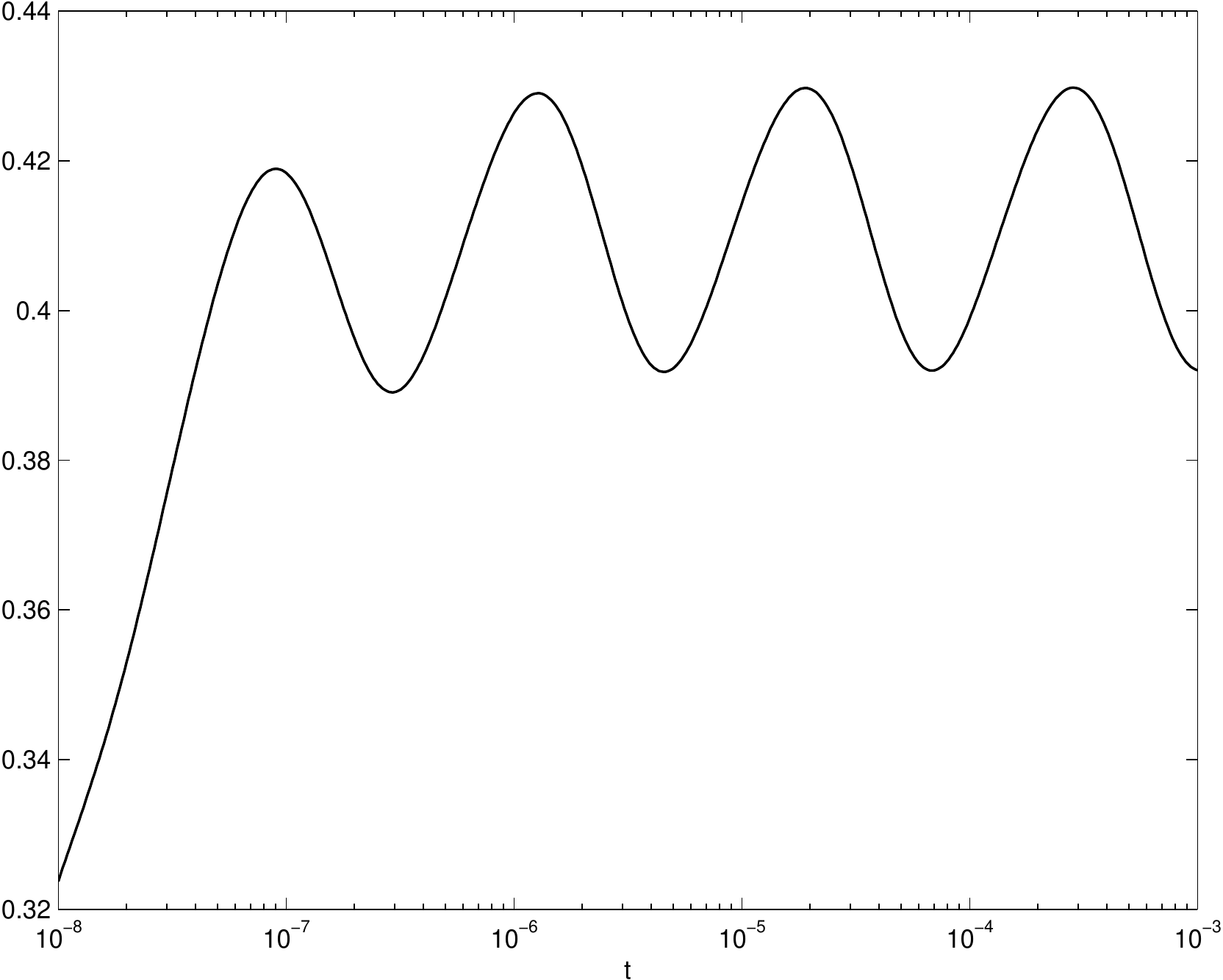}
  \caption{Plot of $t^\alpha h(t,q_0,q_0)$ versus $t$ on the $m=7$ graph
  approximation to $\VS_2$. \label{fig:centerheat}}
\end{figure}

Our data allows us to examine the behavior of the heat kernel
$h(t,q_0,x)$ in greater detail. Estimates for the heat kernel are
known, but they involve constants of unknown size. It is expected that
$h(t,q_0,x)$ should involve a factor of $t^{-\alpha}$ multiplying a
term that drops off exponentially as $x$ moves away from $q_0$. Since
that data in Figure~\ref{fig:centerheat} suggest that the
$t^{-\alpha}$ factor is modified by an oscillating factor, we look at
the ratio $h(t,q_0,x)/h(t,q_0,q_0) = H(t,x)$. Actually, it seems more
plausible that
\begin{align*}
  \frac{h(t,q_0,x)}{\sqrt{h(t,q_0,q_0)} \sqrt{h(t,x,x)}} 
\end{align*}
will be better behaved than $H(t,x)$, but since we don't know how to
compute $h(t,x,x)$ effectively, this isn't an option. Note that
$H(t,x)$ is normalized so that $H(t,q_0) = 1$. Also, if we ignore the
influence of the boundary, which is certainly very slight for small
$t$, we expect $H(t/15,F_0 x)$ should be very close to $H(t,x)$.
Figure~\ref{fig:heatscaling}
illustrates this invariance property.

\begin{figure}[h]
  \centering
  \includegraphics[width=5in]{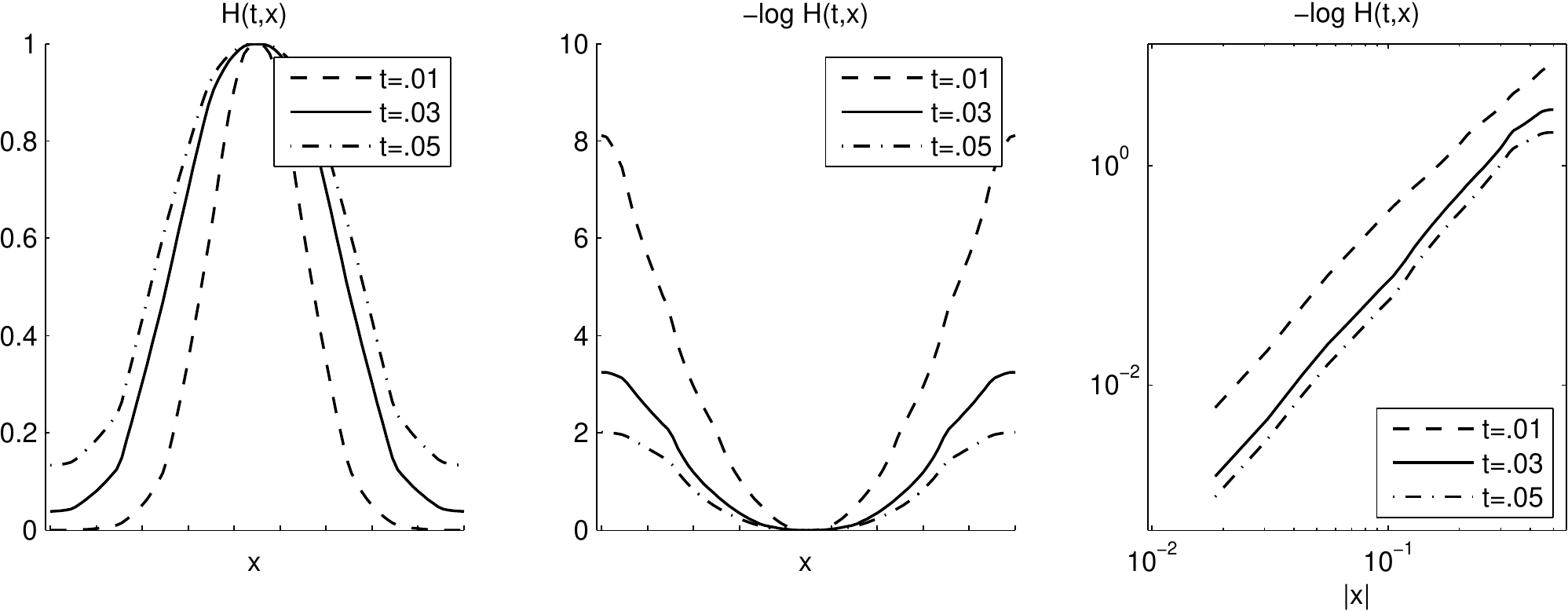}
  \caption{$H(t,x)$ and $-\log H(t,x)$ for several choices of $t$
  (using the $m=4$ graph approximation).
  \label{fig:logheat}}
\end{figure}

First we look at the behavior of $H(t,x)$ for $x$ restricted to the
diagonal. Figure~\ref{fig:logheat}
shows some typical graphs. We also look at $-\log H(t,x)$, again shown
in Figure~\ref{fig:logheat}.  Since $-\log H(t,x)$ vanishes at
$x=q_0$, we try to fit a power law $-\log H(t,x) \approx a|x|^b$ where
the constants $a$ and $b$ depend on $t$, and $|x|$ denotes the
distance to $q_0$. However, we find that the power $b$ varies
significantly as we vary the neighborhood of $q_0$ where we do the
fit. This leads us to doubt the power law model.
Figure~\ref{fig:logheat} shows a log-log plot of $-\log H(t,x)$ for
some choices of $t$.

\begin{figure}[h]
  \centering
  \includegraphics[width=5in]{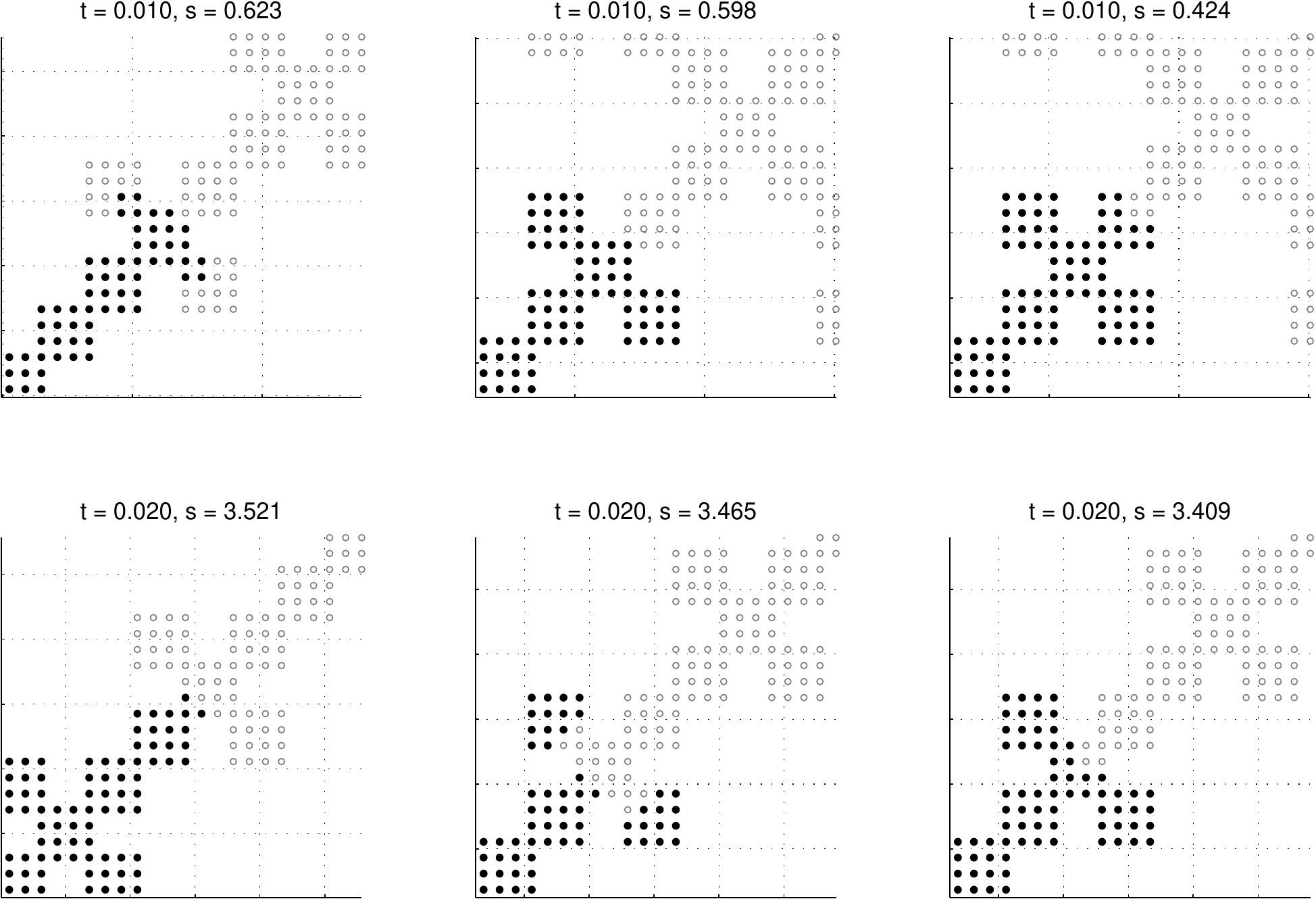}
  \caption{Some heatballs for $t=.01,.02$. \label{fig:heatballs}}
\end{figure}

There is no compelling reason to restrict $x$ to the diagonal in
studying the heat kernel. In a crude sense, the heat kernel
$h(t,q_0,x)$ should depend on the distance of $x$ to $q_0$ in the
resistance metric, which coincides with geodesic distance in $\VS_2$.
But in fact, this is not very accurate. What we want to look at are
what might be called ``heatballs'', sets of the form
\begin{align*}
  \{ x : h(t,q_0,x) \ge s \} 
\end{align*}
for different choices of $t$ and $s$. A naive guess would be that the
heatballs form a 1-parameter family of sets, at least if we stay
toward the center of $\VS_2$ where the influence of the boundary is
small. Again this is only valid in a crude sense.
Figure~\ref{fig:heatballs}
shows some examples of heatballs for two different choices of $t$ and
a variety of $s$-values. One observation is that heatballs tend to
spread further in directions perpendicular to the diagonal. Decreasing
the value of $s$ increases the size of the heatballs, so we may
imagine that the heatballs for fixed $t$ represent an ``invasion''
that spreads out from the center point $q_0$. By and large the
invasion follows an orderly patter, with a cells that lies on the
diagonal being invaded first at the point closest to $q_0$. However,
there are examples where the invasion jumps around, and this produces
examples of heatballs that are disconnected. Apparently, disconnected
heatballs also may occur in the setting of manifolds \cite{laurent}.
Of course, it is also possible to study invasions with $s$ fixed and
$t$ increasing.

\begin{figure}
  \centering
  \includegraphics[width=2.5in]{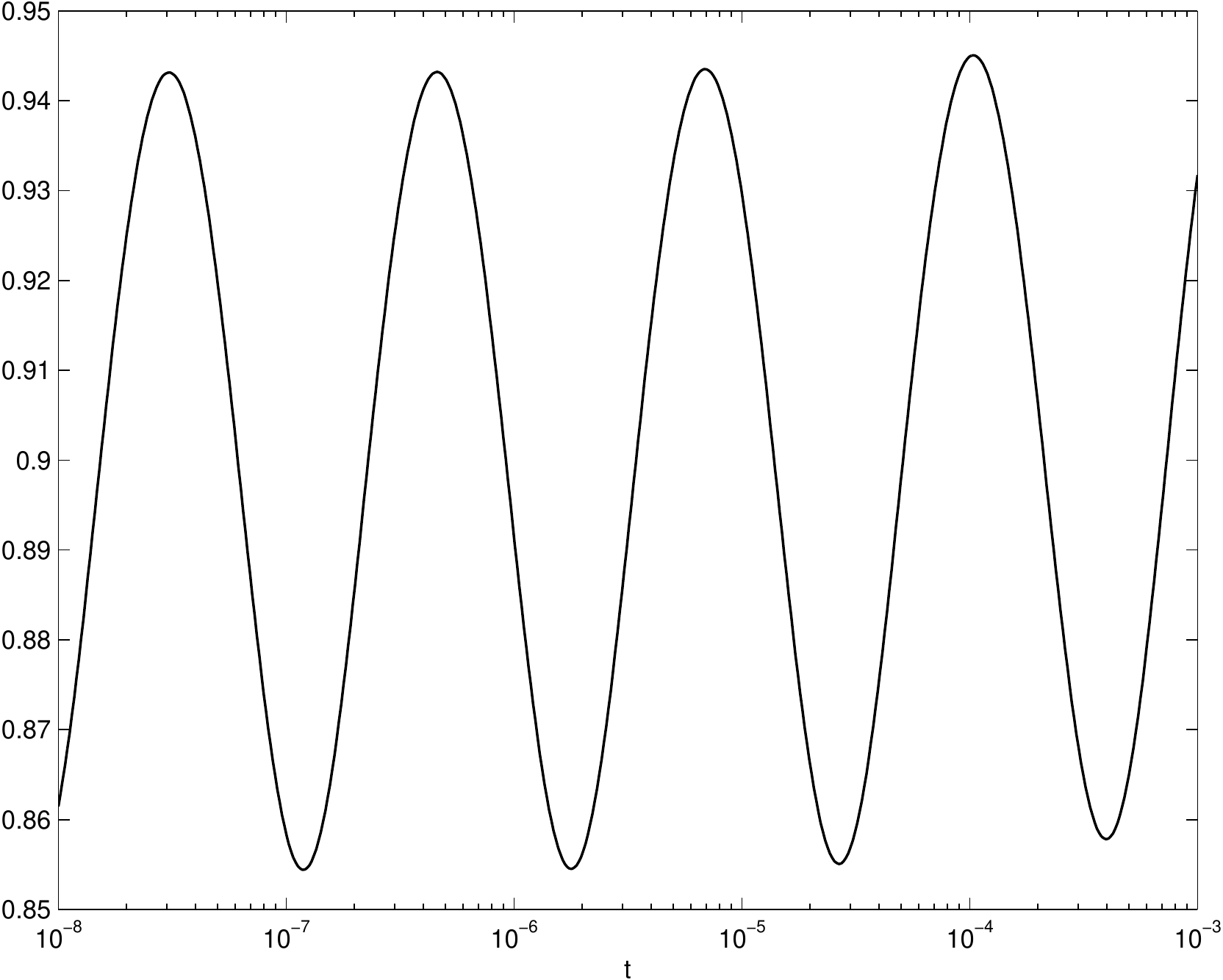}
  \caption{Plot of $t^\alpha \int h(t,x,x)\, d\mu(x)$ versus $t$ on
  the $m=7$ graph approximation to $\VS_2$.
  \label{fig:heatkerneltrace}}
\end{figure}

The trace of the heat kernel and its value at the center, when
multiplied by $t^\alpha$, are both periodic in $\log t$ (see
\cite{grabner}). This is shown in ./figures \ref{fig:heatkerneltrace}
and \ref{fig:centerheat} on the $m=7$ graph approximation. The
approximate sinusoidal behavior is explained for the trace in
\cite{spectralops}, and at the center in \cite{grabner}.  We note that
the approximate sines are out of phase: Fitting to $a+b\sin(c\log t +
d)$ we get $a = .90$, $b = .045$, $c = 2.33$, and $d = -2.2$ for the
trace of the heat kernel, and $a = .4110$, $b=.0191$, $c = 2.33$, and
$d=1.803$ for the heat kernel at the center.

\subsection{Wave Propagator}
The wave equation is given by 
\begin{align*}
  \frac{\partial^2 u}{\partial t^2} = \Delta u 
\end{align*}
If we impose Neumann boundary conditions and initial conditions
$u(x,0) = 0$, $\frac{\partial}{\partial t} u(x,0) = f(x)$, then the
solution is given by 
\begin{align*}
  u(x,t) = \int W(t,x,y) f(y) d\mu(y) 
\end{align*}
where the
{\em wave propagator} $W(t,x,y)$ is given by 
\begin{align*}
  W(t,x,y) = \sum \frac{\sin \sqrt{\lambda_j} t}{\sqrt{\lambda_j}}
  u_j(x) u_j(y).
\end{align*}

\begin{figure}
  \centering
  \includegraphics[width=5in]{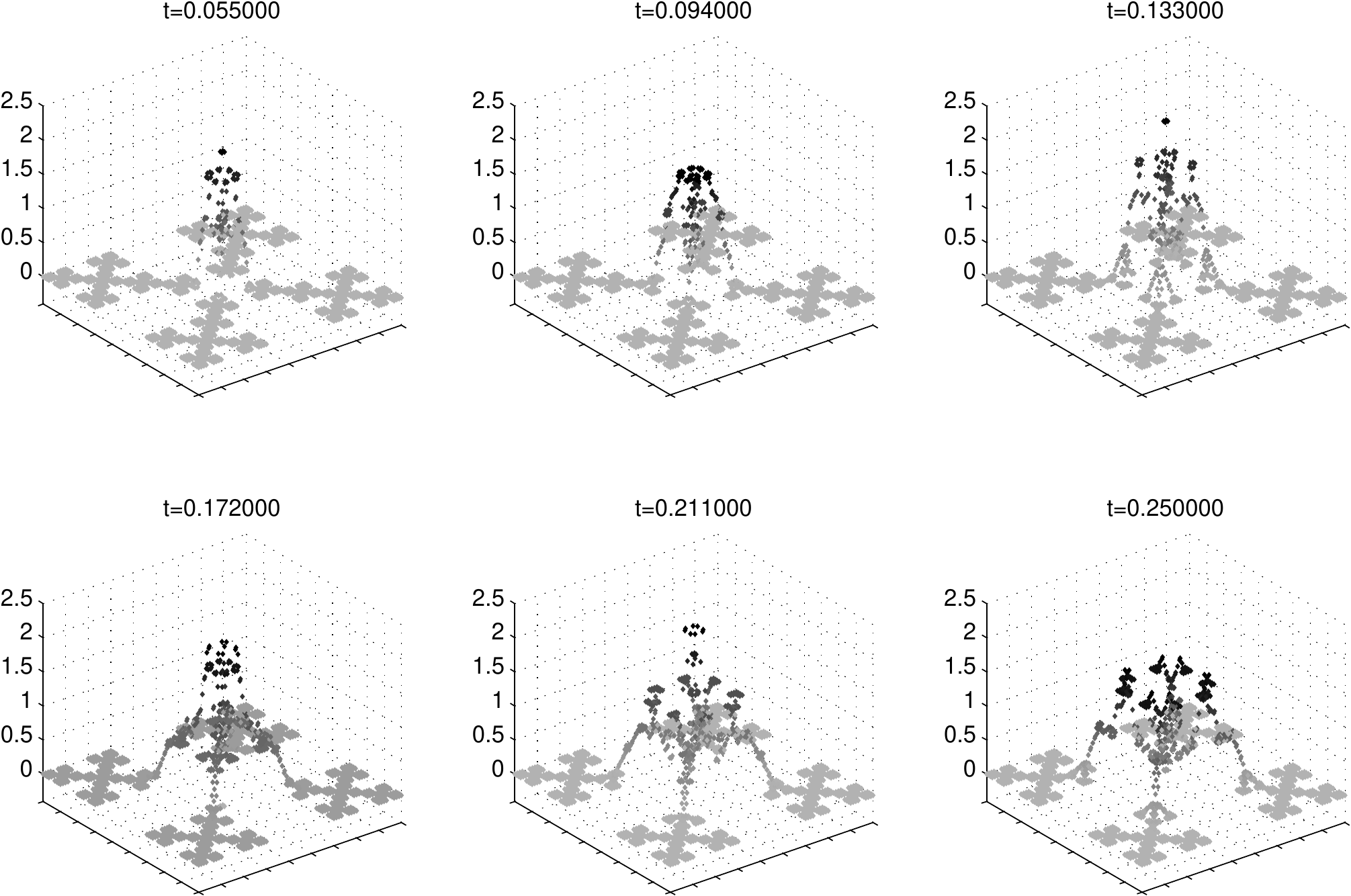}\vspace{10ex}
  \includegraphics[width=5in]{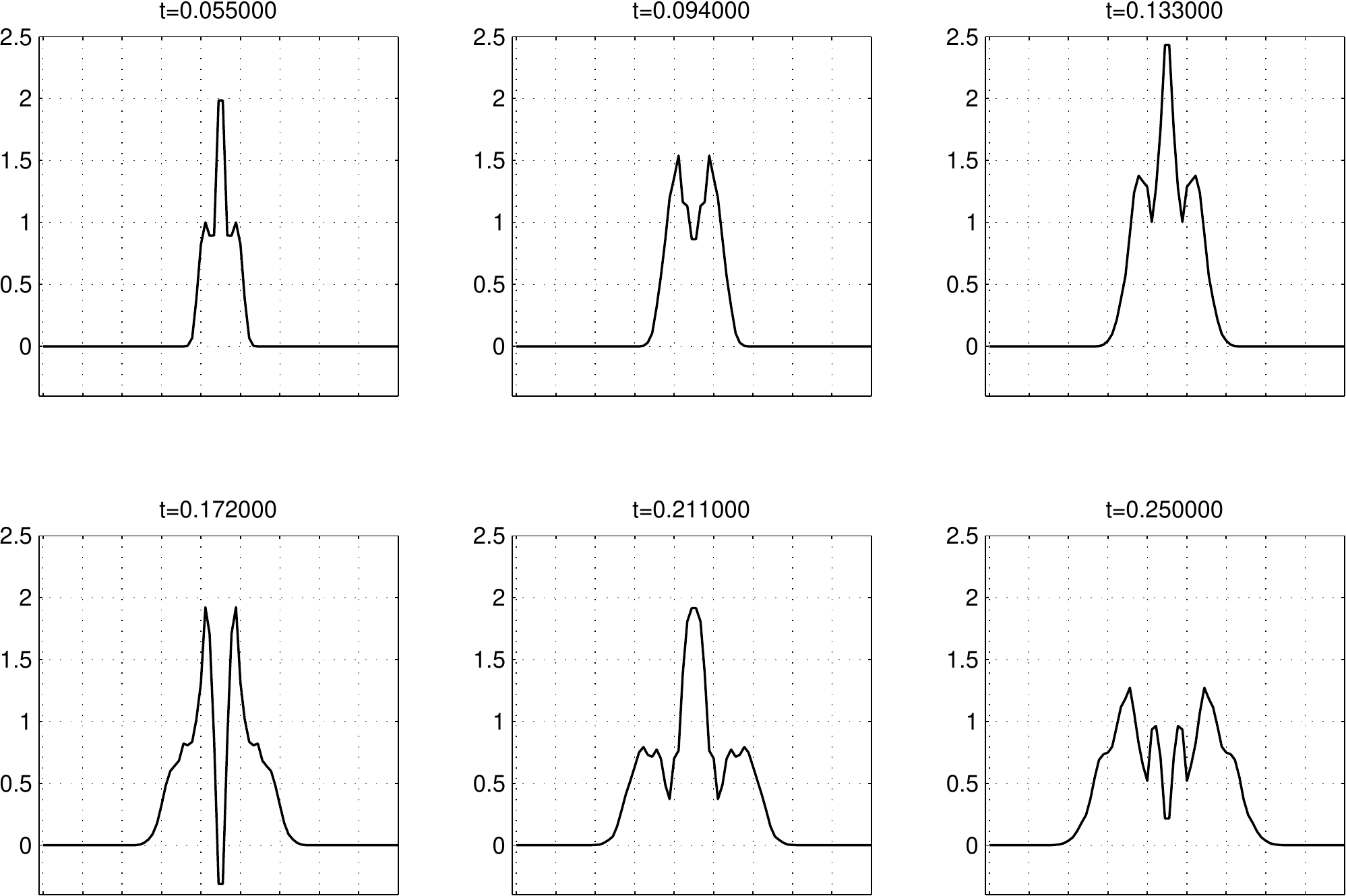}
  \caption{The wave propagator $w(t,q_0,x)$ on $\VS_2$ for various
  values of $t$, on the whole Vicsek set and restricted to the
  diagonal.  \label{fig:wave}}
\end{figure}

From the eigenvalues and eigenfunctions we can also construct the wave
propagator on the standard Vicsek set. As with the heat kernel, this
is easiest to compute when on of the arguments is the center point of
the Vicsek set, since then we only need to consider 0-series
eigenfunctions. This is shown in Figure~\ref{fig:wave}.

As already observed in the case of $\SG$ in \cite{dalrymple} the wave
propagator $W(t,q_0,x)$ is not supported in a small neighborhood of
$q_0$ for fixed $t$; in other words, waves propagate at infinite
speed. This is easily explained because the differential operators on
either side of the wave equation do not have the same order. However,
the amount of energy that propagates at high speed is relatively
small. So we can expect a weak substitute for finite propagation
speed. Attempts to understand this in \cite{dalrymple} and
\cite{coletta} were stymied by the complexity of the wave propagator
on $\SG$ (in \cite{spectralops}
it was shown that time integrals of the wave propagator are
computationally tamer on $\SG$, but this did not help with a weak
finite propagation speed).

\begin{figure}
  \centering
  \includegraphics[width=5in]{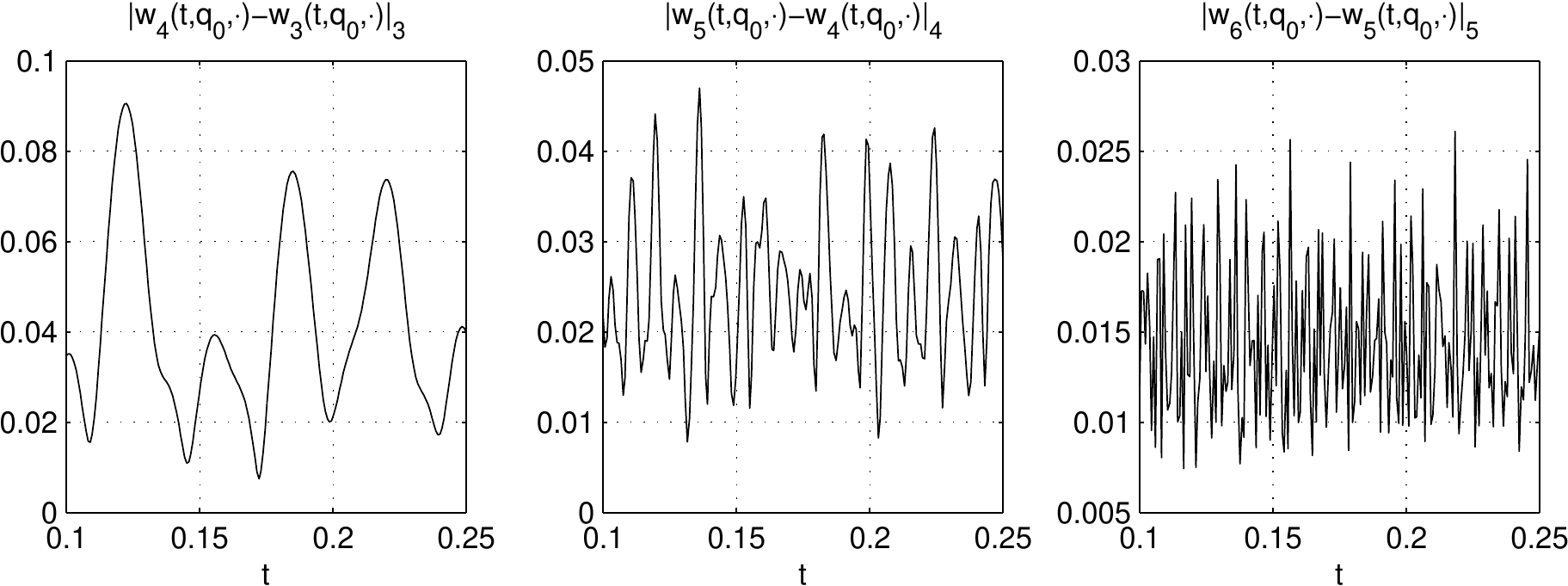}
  \caption{$L^2$ distances between $w_m(t,q_0,\cdot)$ and
  $w_{m-1}(t,q_0,\cdot)$. \label{fig:stability}}
\end{figure}

On $\VS_2$ the wave propagator at the center point may be effectively
computed. In particular, when we increase the level of approximation
the graph does not change appreciably: Figure~\ref{fig:stability}
shows $L^2$ distances between $w_m(t,q_0,\cdot)$ and
$w_{m-1}(t,q_0,\cdot)$, where $w_m$ is the level $m$ approximation to
the wave propagator.
In Figure~\ref{fig:wave} we display the graphs for some values of $t$.
Unlike the heat kernel, the wave propagator is not known to be
positive, and indeed we see time where negative values occur. We know
\begin{align*}
  \int W(t,q_0,x) d\mu(x) =t 
\end{align*}
so that positive values predominate, and it seems from the data that
\begin{align*}
  \int |W(t,q_0,x)| d\mu(x)
\end{align*}
is bounded by a multiple of $t$. Recall that in Euclidean space, the
singularity of the wave propagator worsens as the dimension increases.
Our data is more in line with the $n=1$ case.

\begin{figure}
  \centering
  \includegraphics[width=5in]{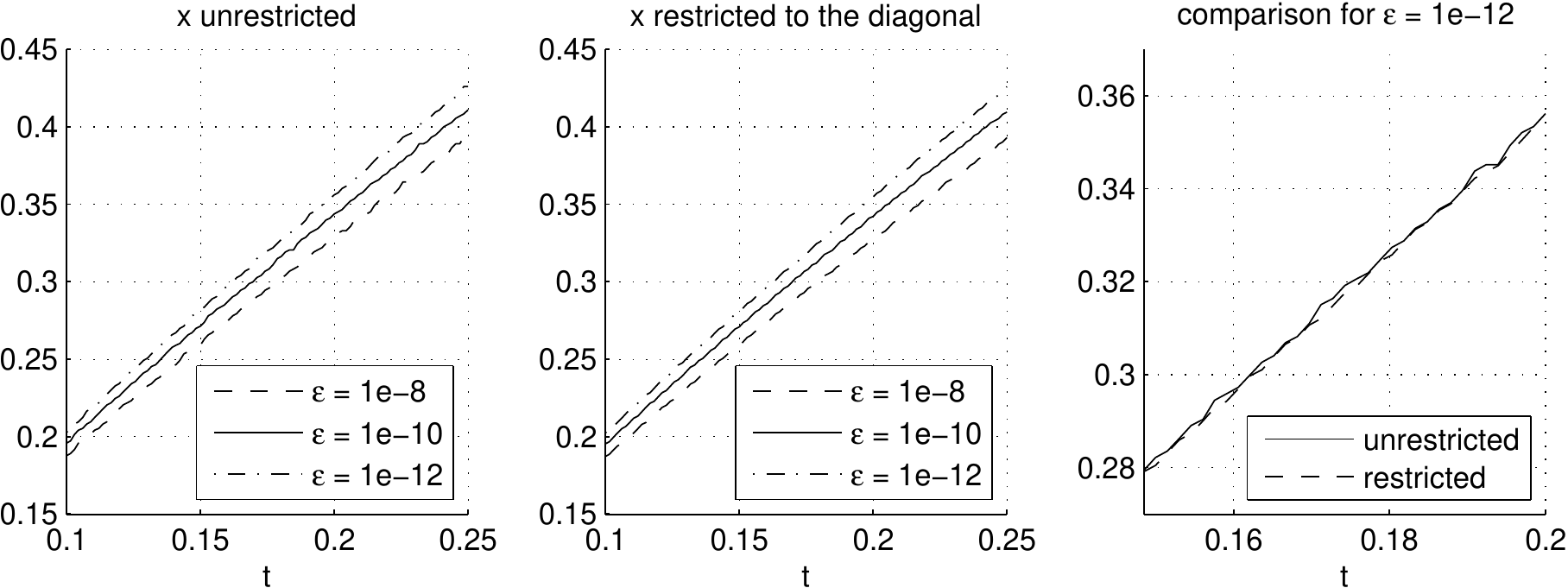}
  \caption{$\max\{|x| : |w_6(t,q_0,x)| \ge \varepsilon\}$ for several
  choices of $\varepsilon$, both in the case when $x$ is restricted to
  the diagonal and when $x$ varies over all $\VS_2$.\label{fig:width}}
\end{figure}

Our data strongly suggests an approximate finite propagation speed. We
can quantify this by choosing a small cutoff $\varepsilon$ and looking
for the maximum value of $|x|$ where $|w_m(t,q_0,x)| \ge \varepsilon$
for fixed $t$, and then letting $t$ vary. In Figure~\ref{fig:width}
we show plots of this function, both in the case when $x$ is
restricted to the diagonal and in the case where $x$ varies over all
$\VS_2$, for different choices of $\varepsilon$.
Notice that in both cases the slope of the function increases with
$\varepsilon$.

\subsection{Spectral Projections}
Another important class of spectral operators are the spectral
projections. Let $\Lambda$ be a subset (usually infinite) of the
spectrum, and define
\begin{align*}
  P_\Lambda f = \sum_{\lambda\in\Lambda} P_\lambda f 
\end{align*}
where
\begin{align*}
  P_\lambda f(x) = \sum \int u_j(x) u_j(y) f(y) d\mu(y) 
\end{align*}
for $\{u_j\}$ an orthonormal basis of the $\lambda$-eigenspace. Such
operators are always bounded on $L^2$ (with norm 1) and usually not
bounded on $L^1$ or $L^\infty$. A natural question to ask is under
what conditions is $P_\Lambda$ bounded on $L^p$ for $1< p < \infty$?
In the classical setting such results can be obtained from the
Marcinkiewicz multiplier theorem \cite{stein} and we expect that
analogous results should be valid in the fractal setting, perhaps
related to the transplantation theorems of \cite{duong} and
\cite{sikora}.  We note that the results of \cite{gaps} imply that we
can always ``segment'' such problems; we write $\Lambda_k = \Lambda
\cap [0,N_k]$ for a natural sequences of cutoffs $N_k$ that lie at the
beginning of spectral gaps (in our case we take spectral decimation
through level $k$). Then $P_\Lambda$ is bounded on $L^p$ if and only
if $P_{\Lambda_k}$ is uniformly (in $k$) bounded on $L^p$.

In \cite{spectralops} we looked at some spectral projection on $\SG$, but
it was difficult to arrive at meaningful predictions because of the
computational complexity of the data. Here we are able to examine one
example in detail: the case that $\Lambda$ consists of the 0-series
eigenvalues. Because these eigenvalues all have multiplicity one, it
is straightforward to compute kernels $K_k(x,y)$ of the segmented
projection operators $P_{\Lambda_k}$ for $k \le 5$ on $\VS_2$. We make
a few simple observations. The first is that
\begin{align*}
  \int K_k(x,y) d\mu(y) = 1
  \qquad
  \text{for every $x$}.
\end{align*}
This follows from the fact that the constant 1 is in the 0-series, and
every other 0-series eigenfunction is orthogonal to it. The second is
that
\begin{align*}
  K_k(x,y) 
  = K_k(\Phi(x),y) 
  = K_k(x,\Phi(y)) 
  = K_k(\Phi(x),\Phi(y)) 
\end{align*}
where $\Phi$ is any isometry of $\VS_2$. This is an immediate
consequence of the fact that each 0-series eigenfunction is invariant
under $\Phi$ (if $u$ is a 0-series eigenfunction then so is $u \circ
\Phi$, with the same eigenvalue, and the multiplicities are all one).
Incidentally, we remark that invariance under all isometries does not
characterize the 0-series spectrum; it easy to construct 4/3-series
eigenfunctions (on a higher level) that show this invariance.

We examine the behavior of $\int |K_k(x,y)| d\mu(y)$ as a function
$k$. Table~\ref{fig:abskernel}
\begin{table}
  \begin{tabular}{c|c}
    $k$  &  $\max_x \textstyle \int |K_k(x,y)| d\mu(y)$\\
    \hline
    1 & 1.4476 \\      
    2 & 1.7336 \\      
    3 & 2.9958 \\      
    4 & 4.7955 \\      
    5 & 7.6572 \\ 
  \end{tabular}
  \vspace{2ex}
  \caption{Maximum value of $\int |K_k(\cdot,y)| d\mu(y)$ for several
  $k$. \label{fig:abskernel}}
\end{table}
\begin{figure}
  \centering
  \includegraphics[width=5in]{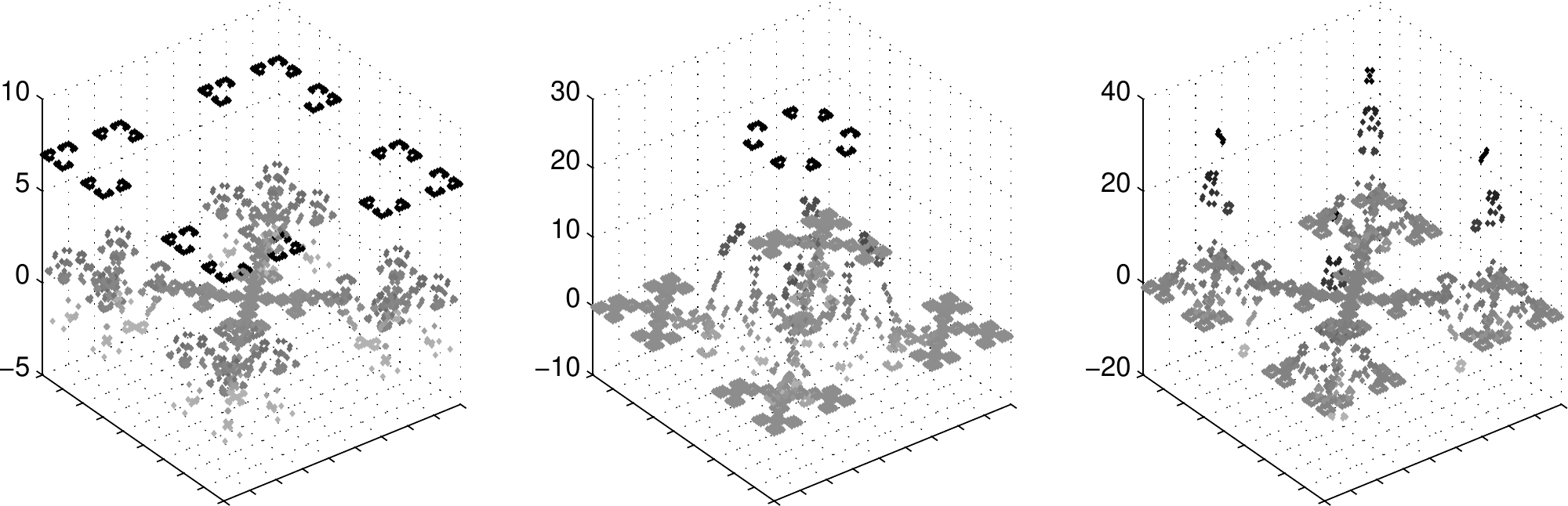}
  \caption{$K_4(x,\cdot)$ for $x=q_1, F_{0233} q_3, F_{2042}q_2$
  \label{fig:projection}}
\end{figure}
\begin{figure}
  \centering
  \includegraphics[width=3.5in]{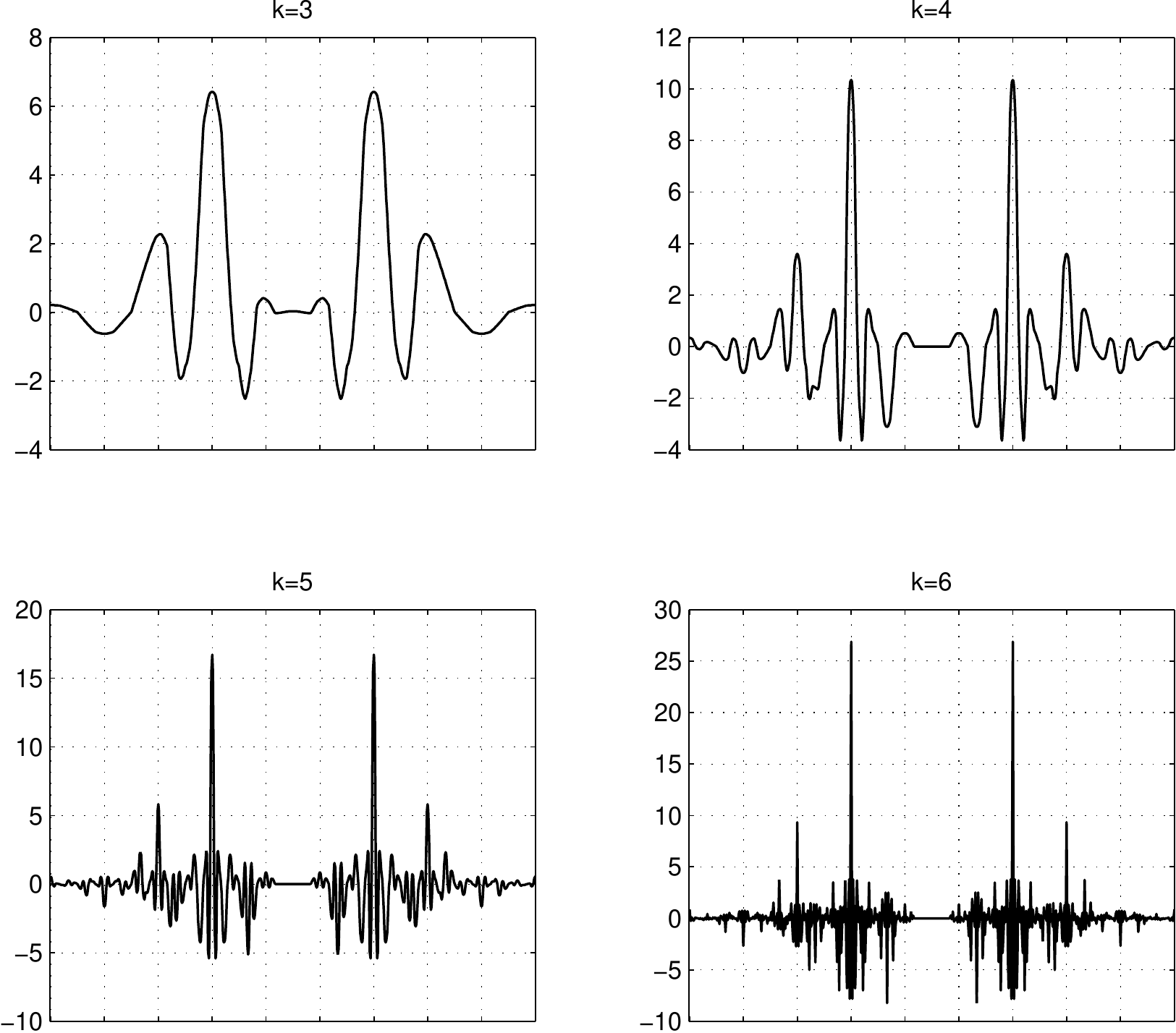}
  \caption{The restriction of $K_k(x,\cdot)$ to the diagonal when $x$
  is the junction point between two 1-cells, for $3\le k \le 6$.
  \label{fig:projectiondiag}}
\end{figure}
shows the maximum over $x$ for $k\le 5$. This is overwhelming evidence
that $\max_x \int |K_k(x,y)| d\mu(y) \to \infty$ as $k \to \infty$,
and this implies that $P_\Lambda$ is not bounded in $L^1$ or
$L^\infty$. Next we ask, for fixed $x$, what are the $y$ values where
$|K_k(x,y)|$ is large? Looking at the graphs of $K_k(x,\cdot)$ in
Figure~\ref{fig:projection}
we see evidence that the answer is the values of $y$ that are close to
$\Phi(x)$ for some isometry $\Phi$. Note that for some choices of $x$,
the set of all $\Phi(x)$ is finite, but for other choices it may be
infinite. (For example, if $x$ is a boundary point, then it is a dense
subset of a Cantor set that includes the intersection of $\VS_2$ with
the boundary of the unit square.)

In Figure~\ref{fig:projectiondiag} we show the restriction to the
diagonal of $K_k(x,\cdot)$ when $x$ is the junction point between two
1-cells, for $3 \le k \le 6$. The behavior is certainly more
complicated than the kernels in the standard Calderon-Zygmund theory.
On the other hand, the graphs to appear to be converging to some
limiting shape. It would be interesting to make this statement more
precise, and to investigate whether there is $L^p$ boundedness of
$P_\Lambda$ for some values of $p$ in $1<p<\infty$ other than $p=2$.

\section{Diagonals and the 0-series}\label{sec:diag}

We can write $L^2(\VS) = \mathcal{H}_0 \oplus \mathcal{H}_{4/3}$ where
$\mathcal{H}_0$ represents the eigenfunctions associated with the
0-series, and $\mathcal{H}_{4/3}$ represents those associated with
the 4/3-series.  These are orthogonal because the eigenvalues are
distinct.

\begin{thm}
Each 0-series eigenfunction of the Laplacian on the Vicsek Set is
determined by its restriction to the diagonal.
\begin{proof}
  If we look at the fractal Laplacian, we can view $\VS$ as the union
  of the diagonals with little trees attached, each tree a small copy
  of $1/4\, \VS$, one arm of the Vicsek set.   $u|_T$ satisfies
  $-\Delta u = \lambda u$, with $\partial_n u = 0$ at the outer
  boundary, and $u(q_0)$ is a specified value if $q_0$ is the center
  point, because the center point lies on the diagonal.  

  Let $v_\lambda$ denote the function on $1/4 \, \VS$ that satisfies
  $-\Delta v_\lambda = \lambda v_\lambda$, $\partial_n v_\lambda(q) =
  0$ if $q$ is a boundary point, and $v_\lambda(q_0) = 1$.  To show
  existence and uniqueness, we have to show that $-\Delta u = \lambda
  u$ on $1/4\, \VS$, $\partial_n u (\textup{boundary}) = 0$, and
  $u(q_0) = 0$ imply that $u$ must be identically zero.  Indeed, given
  such a function $u$, extend it by odd reflection across the center
  to the opposite arm of the Vicsek set, and set it identically zero
  on the other two arms.  Then we obtain a global eigenfunction
  satisfying $\sum u(q_j) = 0$ for the boundary points $q_j$, so it
  belongs to the 4/3-series.  But $\lambda$ is a 0-series eigenvalue,
  and by spectral decimation, there are no simultaneous 0-series and
  4/3-series eigenvalues;  the only way out is if $u = 0$.

  Now let $T$ denote any tree of level $m$ that attaches to the
  diagonal at $y$.  Then there exists $\psi_T: T \to 1/4\, \VS$ with
  $\psi_T(y) = q_0$ and $\psi_T(\textup{bdry}(T)) =
  \textup{bdry}(1/4\, \VS)$, and
  \begin{align*}
    \Delta (f \circ \psi_T) = (15)^m (\Delta f) \circ \psi_T.  
  \end{align*}
  This says that any tree can be put in one-to-one correspondence with
  an arm of the Vicsek set in such a way as to respect similarities.

  Let $u$ be our 0-series eigenfunction. Then $(u|_T) \circ
  \psi_T^{-1} = f$ satisfies
  \begin{align*}
    -\Delta f &= (15)^{-m} \lambda f
    \quad \text{on $1/4\, \VS$,}\\
    \partial_n f(\textup{bdry}) &= 0,\\ 
    f(q_0) &= u(y).
  \end{align*}
  Since $(15)^{-m} \lambda$ is not a 4/3-series eigenvalue (if it
  were, then so would $\lambda$ be) we have $f = u(y) v_{15^{-m}
  \lambda} \circ \psi_T$.
  Hence
  \begin{align*}
    u|_T = u(y) v_{15^{-m} \lambda} \circ \psi_T.
  \end{align*}
  So $\lambda$ and $u|_{\rm diagonal}$ determine $u$ according to the
  above equation.
\end{proof}
\end{thm}

We would like to go further and say that any function in
$\mathcal{H}_0$ is determined by its restriction to the diagonal, and
aside from symmetry there are essentially no other conditions on 
the restrictions to the diagonal of $\mathcal{H}_0$ functions. We
begin with the analogous statement on the discrete approximations. 

Let $D_m$ denote the intersection of $V_m$ with one arm of the
diagonal. Note that $\# D_0 = 1$, $\# D_1 = 2$ and $\# D_m = 3\#
D_{m-1} -1$. Let $Z_m$ denote the span of the 0-series eigenfunctions
through level $m$. We may consider elements of $Z_m$ as functions
either on $V_m$ or $\VS_2$. Note that $\dim Z_0 = 1$, $\dim Z_1 = 2$
and $\dim Z_m = 3(\dim Z_{m-1})-1$, so $\dim Z_m = \# D_m = \frac 12
(3^m+1)$. Thus it is plausible that every function on $D_m$ is the
restriction of a function in $Z_m$, and every function in $Z_m$ is
uniquely determined by its restriction to $D_m$. In fact, these
statements are equivalent. We conjecture a little more.

Let $x_0,x_1,x_2,\ldots$ denote the points in $D_m$ moving from the
outside inward. Let $D_m^{(k)} = \{x_0,\ldots,x_k\}$ and let
$V_m^{(k)}$ denote all points in $V_m$ that attach to one of the
midpoints of the intervals in $D_m^{(k+1)}$.
\begin{conj}
  If $f \in Z_m$ vanishes on $D_m^{(k)}$, then it vanishes on
  $V_m^{(k)}$. In particular, $f$ is determined on $V_m$ by its values
  on $D_m$.
\end{conj}
This conjecture implies that there is a formula of the form 
\begin{align*}
  f(z) = \sum c_k^{(m)}(z) f(x_k) 
\end{align*}
for all $f \in Z_m$ for each $z \in V_m$ and
suitable coefficients $c_k^{(m)}$. In fact $f(z) = c_k^{(m)}(z)$
defines the function in $Z_m$ that takes on the values $f(x_j) =
\delta_{jk}$ on $D_m$. The conjecture implies $c_k^{(m)}(z) = 0$ if $z
\in D_m^{(k)}$. 
\begin{figure}
  \centering
  \includegraphics[width=5.0in]{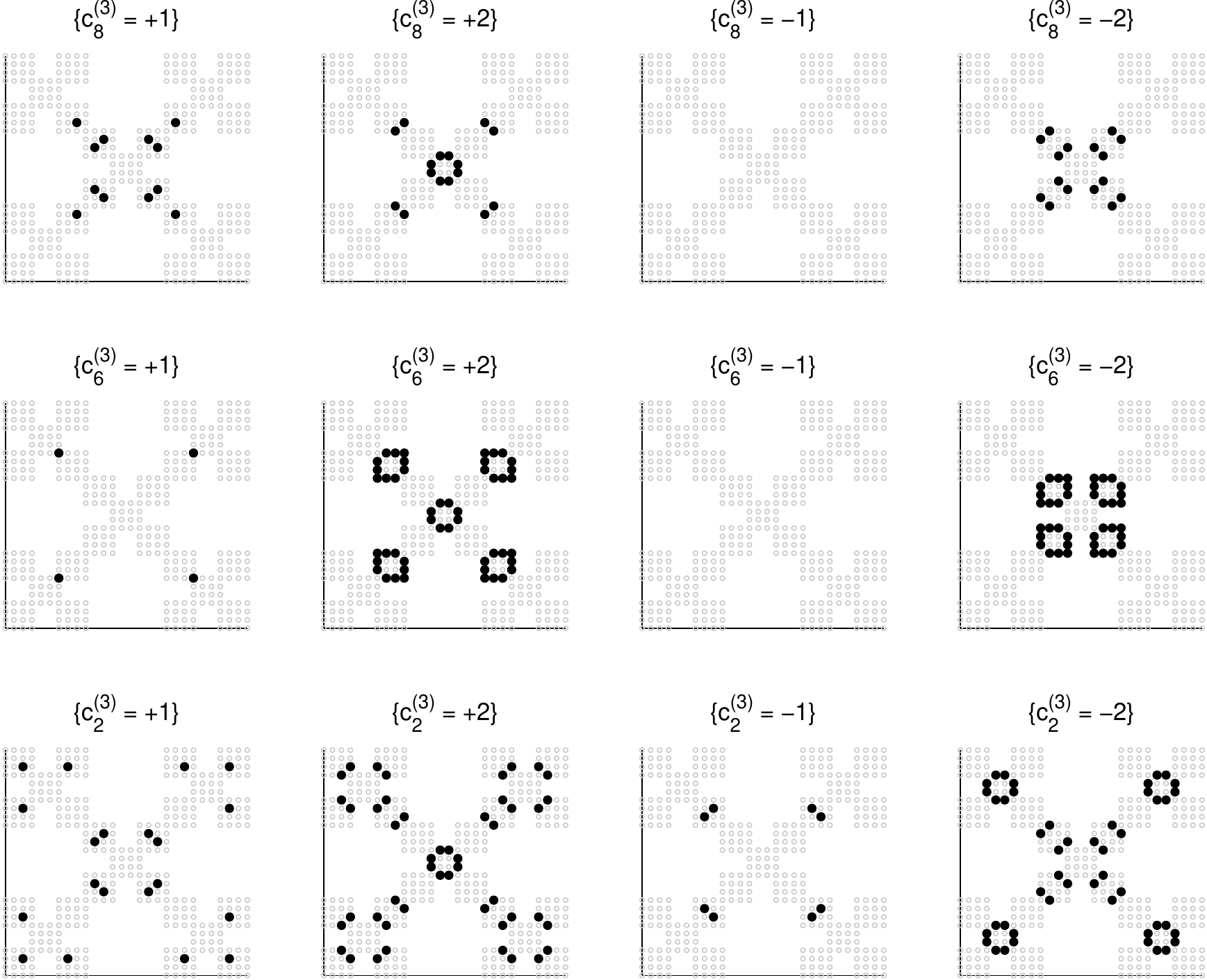}
  \caption{The values of $c_k^{(3)}$ for three choices of $k$.
  \label{fig:restriction}}
\end{figure}
Experimental evidence indicates that $c_k^{(m)}$ can take on only 4
nonzero values, namely $\pm 1$ and $\pm 2$. Some examples of
$c_k^{(3)}$ are shown in Figure~\ref{fig:restriction} and many more
are available (for $m=4$) on our website \cite{website}.

To pass from the discrete to the continuous version we consider
function $\mathcal{H}_0 \cap C$ (here $C$ denotes the continuous
functions on $\VS_2$). Such functions have well-defined restrictions
to $D$ (one arm of the diagonal). To show that the restriction $f|_D$
of such a function determines $f$, it suffices to show that it
determines $f|_{V_m}$ for all $m$, since $\cup_m V_m$ is dense in
$\VS_2$ and $f$ is continuous. Let $f_m$ denote the projection of $f$
onto $Z_m$. By the results of \cite{gaps} we know $f_m$ converges to
$f$ uniformly. If the conjecture is valid then $f_{m'}(z) = \sum
c_k^{(m)} f_{m'}(x_k)$ for $z \in V_m$ and $m' \ge m$ (since
$f_{m'}|_{V_m} \in Z_m$), so passing to the limit $f(z) = \sum
c_k^{(m)} f(x_k)$ for $z \in V_m$. Despite the fact that this is a
finite sum for each $m$, it is a rather peculiar formula. The
coefficients oscillate rapidly but do not go to zero as $m$ increases.
It does not seem likely that we can make any sense out of it if we do
not assume that $f$ is continuous. It seems unlikely that the
existence of a continuous restriction to $D$ for a function in
$\mathcal{H}_0$ implies that it is continuous on $\VS_2$. A more
plausible conjecture is that if the restriction to $D$ is H\"older
continuous of some order then the function is H\"older continuous of
the same order on $\VS_2$. Another reasonable conjecture is that the
restrictions of $\mathcal{H}_0 \cap C$ to $D$ form a dense subset of
the continuous functions on $D$. A less likely conjecture is that the
restrictions give all continuous functions on $D$.

\section{Ratio Gaps}\label{sec:gaps}

\begin{figure}
  \centering
  \includegraphics[width=3.5in]{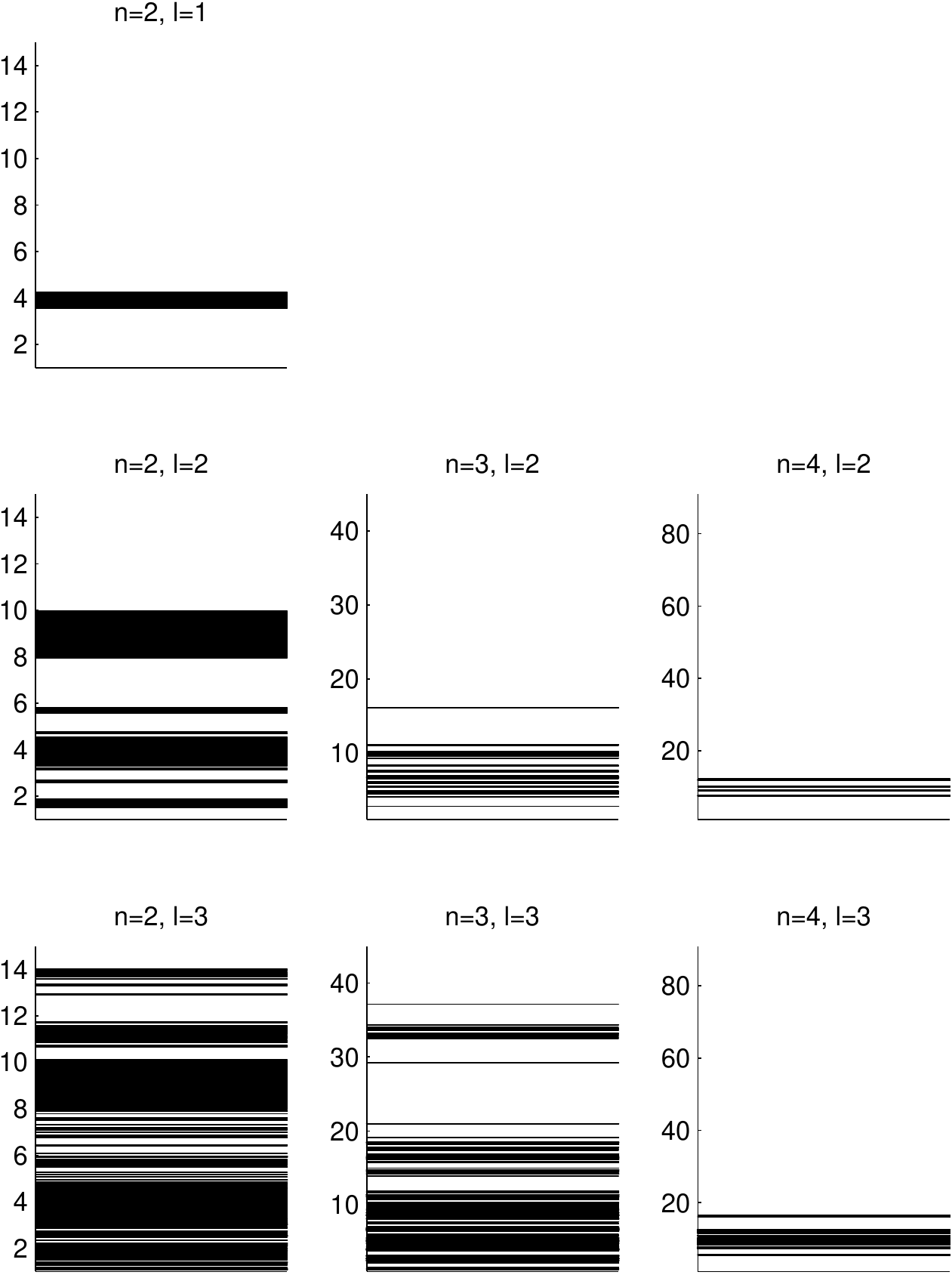}
  \caption{Ratio gaps for $n=2,3,4$ and $\ell = |w| = 1,2,3$. The gaps are
  the black regions. No gaps were found for $n=3,4$ and $\ell = 1$.
  \label{fig:gaps}}
\end{figure}
In \cite{products} it was shown that on $\SG$ there exist gaps in the
ratios of eigenvalues. As a consequence, it is possible to define
operators of the form $\Delta' - a\Delta''$ on the product of two
copies of $\SG$ ($\Delta'$ and $\Delta''$ denote the Laplacian on each
copy of $\SG$) where $a$ lies in a gap, and these operators
paradoxically behave in some ways like elliptic operators, despite the
fact that the coefficient $-a$ has the wrong sign. These operators
were called quasielliptic in \cite{products}. There are no analogous
operators in classical PDE theory. Thus it is of great interest to
know whether similar operators exist for products of fractals other
than $\SG$. In fact \cite{ratiogaps} shows that this is the case for
$\VS_2$ and $\VS_3$. Also \cite{drenning} investigates this question
for a variant of the $\SG$ type fractal. The method used in
\cite{drenning}, which we follow here, yields a computer-assisted
proof. The idea is that the method introduced in \cite{products} leads
to a large number of tedious calculations, and these are best left to
the computer. In our method there is a parameter $\ell$ that may be
chosen at will. Increasing $\ell$ will do a better job finding gaps,
at the cost of increasing the number of computations.

Let $\lambda_m$ be a graph eigenvalue born on level $m_0$. Then
\begin{align*} 
  \lambda_m
  = \phi_{w_m} \circ \phi_{w_{m-1}} \circ \cdots \circ
  \phi_{w_{m_0+1}}(\lambda_{m_0})
\end{align*}
where $\lambda_{m_0} \in \{0,4/3\}$. Let $p$ be the fixed point of
$\phi_{2n-1}$ and $q = \phi_{2n-1}(4/3)$. If $\lambda_{m_0} = 0$ we
have $\lambda_m \le p$, while if $\lambda_{m_0} = 4/3$ then we
either have $m=m_0$ and $\lambda_m = 4/3$ or $\lambda_m \le q$.
Since $q \ge p$, we can simply write
\begin{align}
  \label{graph}
  \lambda_m \in [0,q] \cup \{4/3\}.
\end{align}

Fix $\ell > 0$. Any fractal eigenvalue $\lambda$ is of the form
\begin{align*} 
  \lambda &= \rho^{m_0} \lim_{m \to \infty} \rho^m \phi_{v_m} \circ
  \phi_{v_{m-1}} \circ \cdots \circ \phi_{v_{m_0+1}} (\lambda_{m_0})
\end{align*}
where all but finitely many of the $v_j = 1$. Thus there must be a
word $w$ of length $\ell$ and some graph eigenvalue $\lambda_m$ so
that
\begin{align*} 
  \lambda 
  &= \rho^{m_0} \lim_{k \to \infty} 
  \rho^k\phi_1^k \circ  \phi_w \circ \phi_{v_m} \circ \cdots
  \circ \phi_{v_{m_0+1}} (\lambda_{m_0})\\
  &= \rho^{m_0+\ell+m} \lim_{k \to \infty} 
  \rho^k\phi_1^k (\phi_w  (\lambda_m))\\
  &= \rho^{m_0+\ell+m} 
  \psi_n(\phi_w  (\lambda_m)),
\end{align*}
($\phi_w = \phi_{w_\ell} \circ \cdots \circ \phi_{w_1}$).  Combining
this with \eqref{graph} we see that every fractal eigenvalue $\lambda$
can be written as
\begin{align}
  \label{lambda}
  \lambda = \rho^r \psi_n(x),
  \qquad
  x \in \bigcup_{|w|=\ell} \phi_w([0,q] \cup \{4/3\}).
\end{align}
for some integer $r$.

Consider the contribution of a word $w$ to the eigenvalues described
by \eqref{lambda}. If $w$ ends in a $1$, then as long as $m > m_0$ can
rewrite
\begin{align*}
   \phi_w \circ \phi_{v_m} \circ \cdots
  \circ \phi_{v_{m_0+1}} (\lambda_{m_0})
  = 
  \phi_1 \circ  \phi_{w'} \circ \phi_{v_{m-1}} \circ \cdots
  \circ \phi_{v_{m_0+1}} (\lambda_{m_0}),
\end{align*}
for some other word $w'$ of length $\ell$ (with one less 1 at the
end), while $m = m_0$ means
\begin{align*}
  \phi_w \circ \phi_{v_m} \circ \cdots
  \circ \phi_{v_{m_0+1}} (\lambda_{m_0})
  =
  \phi_w  (\lambda_{m_0}) \in \left\{\phi_w(0), \phi_w(4/3)\right\}.
\end{align*}
Thus \eqref{lambda} is still valid with $\phi_w([0,q] \cup \{4/3\})$
replaced by $\{\phi_w(0),\phi_w(4/3)\}$ in \eqref{lambda} for every word
$w$ ending in 1. Furthermore we can discard $\phi_w(\lambda_{m_0})$
if it is forbidden. 

So far we've found finitely many intervals $[a_i,b_i]$ (allowing $a_i
= b_i$) so that each eigenvalue $\lambda$ must satisfy
\begin{align*}
  \lambda \in \rho^r \psi_n([a_i,b_i])
  = \rho^r [\psi_n(a_i),\psi_n(b_i)] 
\end{align*}
for some $i$ and $r$. Therefore any ratio of eigenvalues $\lambda/\mu$
must satisfy
\begin{align*}
  \frac\lambda\mu \in \rho^r 
  \left[
  \frac{\psi_n(a_i)}{\psi_n(b_j)},
  \frac{\psi_n(b_i)}{\psi_n(a_j)}
  \right]
  \equiv \rho^r[R_{ij}, S_{ij}]
\end{align*}
for some $r$, $i$, and $j$.  Since $\rho\lambda$ is an eigenvalue if
$\lambda$ is, we can restrict our attention to ratios $\lambda/\mu \in
[1,\rho]$ and hence to the finite number of intervals $\rho^r [R_{ij},
S_{ij}]$ which intersect $[1,\rho]$. The gaps in the union of these
intervals are then guaranteed to be ratio gaps. 

Figure~\ref{fig:gaps} shows the ratio gaps that are proved to exist by
this method for $n=2,3,4$ using values of $\ell = 1,2,3$. For all of
these $n$ there are ratio gaps containing $\sqrt{\rho_n}$, given in
Table~\ref{fig:rhogaps}.
\begin{table}[h]
  \centering
  \begin{tabular}{cc|cc}
    $n$ & $\sqrt{\rho_n}$ & $\ell$ & ratio gap\\
    \hline
    2 & 3.8730 & 1 & $[3.5370,4.2409]$\\
    && 2 & $[3.2948,4.5526]$\\
    && 3 & $[3.2948,4.5526]$\\
    \hline
    3 & 6.7082 & 1 & no gap\\
    && 2 & $[6.6952,6.7212]$\\
    && 3 & $[6.6950,6.7214]$\\
    \hline
    4 & 9.5394 & 1 & no gap\\
    && 2 & no gap\\
    && 3 & $[9.5357,9.5431]$
  \end{tabular}
  \vspace{2ex} \caption{Ratio gaps containing $\sqrt{\rho_n}$ for
  $n=2,3,4$ and $\ell=1,2,3$. \label{fig:rhogaps} }
\end{table}
We see clearly that the number and size of the ratio gaps increases
with $\ell$. However, we have not been able to confirm the existence
of ratio gaps for $n\ge 5$. For $n=5$ none are revealed for $\ell \le
2$ and our MATLAB implementation (see \cite{website}) runs into memory
problems for $\ell \ge 3$. For $\ell \ge 3$ we can, however, use a
modified algorithm which searches only for ratio gaps containing a
particular point. These searches have failed to find ratio gaps
containing $\sqrt{\rho_5} \approx 12.3693$. It is not clear if these
failed searches should be interpreted as experimental evidence for the
nonexistence of ratio gaps, or just as evidence that we need to
consider higher values of $\ell$ to find ratio gaps.

\section{Eigenvalue clusters}\label{sec:clusters}
We say the spectrum of a Laplacian exhibits spectral clustering if the
following holds: for every integer $n$ and $\varepsilon > 0$ there exists
an interval $I$ of length $\varepsilon$ that contains $n$ distinct
eigenvalues.

This, for example, says you can find a million distinct eigenvalues
within a millionth of each other.  The eigenvalues will have to be
very large, so it becomes computationally challenging to find such
tight and large clusters.  Clustering does not occur on the Sierpinski
gasket $\SG$.  Experimental evidence suggests that it does occur on
the pentagasket \cite{pentagasket} and on the Julia sets \cite{flock}.
The following lemma allows us to prove it holds on $\VS_2$.
\begin{lem}
Suppose spectral decimation holds with spectral renormalization factor
$\rho$ and spectral renormalization function $R(\lambda)$.  Suppose $R$
has a fixed point $t$ ($R(t) = t$) such that $|R'(t)| > \rho.$  Then
spectral clustering occurs.
\begin{proof}
  Let $\phi_1, \phi_2, \ldots \phi_N$ be the inverses of $R(\lambda)$
  in increasing order.  Then $\phi_1(0) = 0$ and $\phi_1'(0) = 1/\rho$.
  There exists $k$ such that $\phi_k(t) = t$ and $|\phi_k'(t)| = b <
  \rho^{-1}$, by the assumption.  Choose $m$ large enough that $\Delta_m$
  has at least $n$ distinct eigenvalues $\lambda_1, \ldots, \lambda_n$.
  Then $\Delta_{m + j}$ has $n$ distinct eigenvalues $\phi_k^{(j)}
  (\lambda_1), \ldots, \phi_k^{(j)} (\lambda_m)$ and these give rise to
  $n$ distinct eigenvalues 
  \begin{align*}
    \lim_{l \to \infty} \rho^{m + j + l} \phi_1^{(l)} \phi_k^{(j)}
    (\lambda_p), 
    \qquad 1 \le p \le n. 
  \end{align*}
  Write $g(\lambda) = \lim_{l \to \infty} \rho^l \phi_1^{(l)} (\lambda)$.
  Then $g$ is a fixed function with bounded derivative $|g'(x)| \le M$ in
  the relevant interval of length $\varepsilon$ where we want to find
  $n$ distinct eigenvalues.

  By taking $j_0$ large enough, we can make all the values
  $\phi_k^{(j_0)} (\lambda_p)$ close enough to $t$ so that
  $|\phi_k'(x)| \le a \le \rho^{-1}$ for all $\phi_k^{j_0}
  (\lambda_p)$.  This means that $\{ \phi_k^{(j_0 + j_1)} (\lambda_p)
  \}$ belongs to an interval of length no more than $c a^{j_i}$ where
  $c$ is the length for $j_1 = 0$.  Then
  \begin{align*}
    \{ \rho^{m + j_0 + j_1} g(  \phi_k^{j_0 + j_1} (\lambda_p)) \} 
  \end{align*}
  belongs to an interval of length at most $c M \rho^{m + j_0}
  (a\rho)^{j_1}$.  Since $a\rho < 1$, this can be made $\le \varepsilon$ by
  taking $j_1$ large enough.  Thus we can find $n$ distinct eigenvalues
  in an interval of length no more than $\varepsilon$.
\end{proof}
\end{lem}
 
On $\VS_2$, $\rho = 15$ and $R(\lambda) = 36 \lambda^3 - 48 \lambda^2 +
15 \lambda$.  So $R(t) = t$ means $2t (18t^2 - 24t + 7) = 0$ with
solutions $0 , \frac{4 \pm \sqrt{3}}{6}$.  We are interested in the
largest $t$, $\frac{4 + \sqrt{2}}{6}$ which is the fixed point of
$\phi_3$.  
\begin{align*}
  R'(\lambda) 
  = 3 \cdot 36 \lambda^2 - 2 \cdot 48 \lambda + 15 
  = 15 + 12 \lambda(9 \lambda - 8)
\end{align*}
so to show $R'(t) > 15$ we need $t > 8/9$.  But $(4 + \sqrt{2})/6 =
0.902...$ so this is true.  We thus have clustering in $\VS_2$.
Computing the largest fixed point $t$ of $R$ on $\VS_n$ for
$n=3,\ldots,9$ we also get $R'(t) > \rho$ (see
Table~\ref{fig:clustering}) and hence that spectral clustering
occurs. Because the ratio $R'(t)/\rho$ increases rapidly with $n$ we
conjecture that spectral clustering occurs for all $n$.
\begin{table}[h]
  \centering
  \begin{tabular}{c|crr@{.}l}
    \multicolumn{1}{c|}{$n$} &
    \multicolumn{1}{c}{$t$} & 
    \multicolumn{1}{c}{$\rho$} & 
    \multicolumn{2}{c}{$R'(t)$} \\\hline
    2   & 0.9024 & 15  & 1&$6314\times 10^1$ \\ 
    3   & 0.8905 & 45  & 1&$3999\times 10^2$ \\ 
    4   & 0.8891 & 91  & 1&$2355\times 10^3$ \\ 
    5   & 0.8889 & 153 & 1&$1079\times 10^4$ \\ 
    6   & 0.8889 & 231 & 9&$9655\times 10^4$ \\ 
    7   & 0.8889 & 325 & 8&$9682\times 10^5$ \\ 
    8   & 0.8889 & 435 & 8&$0713\times 10^6$ \\ 
    9   & 0.8889 & 561 & 7&$2641\times 10^7$
  \end{tabular}
  \vspace{2ex}
  \caption{The largest fixed point $t$ of the spectral decimation
  function $R$ on $\VS_n$ satisfies $R'(t) > \rho$ for
  $n=2,\ldots,9$. \label{fig:clustering}}
\end{table}

\section{Green's Function on $\VS_n$}\label{sec:green}
 
The Green's function $G$ for the Laplacian is a function satisfying
\begin{align*}
  \left\{
  \begin{aligned}
    -\Delta G(x, y) &= \delta(x, y)\\
    G(q_j, y)  &= 0 \quad \text{if $q_j \in \mathop{\rm bdry}
    \VS_n$}.
  \end{aligned}
  \right.
\end{align*}
where $\delta$ is the Dirac delta function.  Then $u(x) = \int G(x, y)
f(y) d\mu(y)$ solves
\begin{align*}
  \left\{
  \begin{aligned}
    -\Delta u &= f, \\
    u|_{\mathrm{bdry}} &= 0.
  \end{aligned}
  \right.
\end{align*}
As a function of $x$, $G$ should be harmonic in the complement of $y$.
Suppose $y$ lies in the upper right arm of  $\VS_n$.  The boundary
points are labeled $q_1, q_2, q_3, q_4$, with $q_1$ corresponding to
the arm where $y$ is.  Let $z$ be the projection of $y$ onto the
diagonal of $\VS_n.$  (In the case that $y$ is on the diagonal already,
$z = y$.)

Now $G(q_j, y)= 0$ for $j = 1, 2, 3, 4$.  Define $G(q_0, y) = a$,
(where $q_0$ is the center point),  $G(z, y) = b$, $G(y, y) = c$.  The
values $a, b, c$ determine $G$ because $G(x, y)$ is linear on the arms
$(q_0, q_2), (q_0, q_3), (q_0, q_4)$, on $(q_0, z)$ and $(z, q_1)$,
and along the unique path joining $z$ to $y$.  It is constant on every
component of the complement of these 6 sets.

To determine the constants $a, b, c$ we have 3 equations that express
$-\Delta_x G(q_0, y) = 0$, $-\Delta_x G(z, y) = 0$, and $-\Delta_x G(y,
y) = \delta_y$. The first two equations say the sum of the 4
derivatives at $q_0$ (resp. $z$) vanish.  The last says the derivative
at $y$ is 1.

For simplicity assume the length of each arm is 1.  (This involves
rescaling by a factor of $\sqrt{2}/2$ for a unit square.)  Let $d(z,
q_0) = s$ (so $d(z, q_1) = 1-s$) and $d(x, y) = t$, measured along the
path.  Then
\begin{align*}
  3a + \frac{a-b}s &= 0 \qquad \text{(at $q_0$)}\\
  \frac{b-a}s + \frac b{1-s} + \frac{b-c}t &= 0 \qquad \text{(at $z$)}\\
  \frac{a-b}t &= 1  \qquad \text{(at $y$)}
\end{align*}
Note that the third equations says $c = b + t$ and the second equation says
\begin{align*}
  \frac{b-a}s + \frac b{1-s} = 1.
\end{align*}
(If $z = y$ then $t = 0$, $c = b$, and the second equation says
$-\Delta_x G(y, y) = 1$ as required.) Solving two equations for two
unknowns yields
\begin{align*}
  a = \frac{ 1-s }4,
  \qquad
  b = \frac{(1-s)(3s + 1)}4.
\end{align*}
In other words, 
\begin{align*}
  G(q_0, y) &= \frac{1 - d(z, q_0)}{4},\\
  G(z, y) &= \frac{(1 - d(z, q_0)) (3d(z, q_0) + 1)}{4},\\
  G(y, y) &= \frac{(1- d(z, q_0))(3 d(z, q_0) + 1)}{4} + d(z, y).
\end{align*}
Denote by $d'(x, q_0)$ the distance from $q_0$ to the point on one of
the main diagonals where $x$ attaches.  So $d'(y, q_0) = d(z, q_0)$, for
example.  Then
\begin{align*}
  G(x, y) = \frac{(1 - d'(x, q_0))(1 - d'(y, q_0))}{4} 
\end{align*}
if $x$ and $y$ lie on different arms.  If $x$ and $y$ lie on the same
arm and $d'(x, q_0) < d'(y, q_0)$ then 
\begin{align*}
  G(x, y) = \frac{(1 - d'(y, q_0))(3 d'(x, q_0) + 1)}{4},
\end{align*}
while if 
$d'(y, q_0) < d'(x, q_0)$ then
\begin{align*}
  G(x, y) = \frac{(1 - d'(x, q_0))(3 d'(y, q_0) + 1)}{4}.
\end{align*}
If $d(y, q_0) = d(x, q_0)$ then
\begin{align*}
  G(x, y) = \frac{(1 - d'(y, q_0))(3 d'(y, q_0 + 1))}{4} + d(z, w(x, y)),
\end{align*}
where $w(x, y)$ is the last point on the intersection of the paths
from $z$ to $y$ and $z$ to $x.$

\section{Higher Vicsek Sets}\label{sec:higher} 
It is clear that $\VS_n$ converges to a cross.  The eigenfunctions of
the Laplacian on the cross are well understood: the restriction to
either diagonal is an eigenfunction on the unit interval, while at the
center point the function is required to be continuous and to have the
sum of its normal derivatives equal to zero. Thus any eigenfunction is
either $\cos \pi kx$ on each diagonal or $a_j \sin \pi (k + 1/2) x$ on
each half diagonal, with $\sum a_j = 0$, for some integer $k$. We call
the first type symmetric, and the second nonsymmetric.  The symmetric
eigenvalues (obtained by taking second derivatives) are $\pi^2 k^2$,
and the nonsymmetric eigenvalues are $\pi^2 (k + 1/2)^2$.  

We claim that the symmetric spectrum is the limit of the spectrum of
the 0-series on $\VS_n$ as $n \to \infty$ (these are symmetric
eigenfunctions), and the nonsymmetric spectrum is the limit of the
spectrum of the 4/3-series born on level 0.  (the 4/3 series born on
levels $\ge 1$ does not contribute to the limit because the
eigenvalues go to infinity.)  We also claim that the limits of the
symmetric eigenfunctions are cosines, and the limits of the
nonsymmetric eigenfunctions are sines.

To understand the behavior of the eigenvalues as $n \to \infty$ we can
restrict attention to the initial segment consisting of the 0-series
eigenvalues $\rho_n \psi_n( \phi_{2j-1}(0))$ and the 4/3-series
eigenvalues $\rho_n \psi_n ( \phi_{2j-1}(4/3))$.  From  \cite{zhou} we
know
\begin{align*}
  r(\lambda) =   \lambda g_n(\lambda) h_n(\lambda),
  \qquad \text{and} \qquad
  3r(\lambda)-4 = f_n(\lambda) l_n(\lambda),
\end{align*}
so $\phi_k(0)$ are the zeroes of $g_n$ and $h_n$ and $\phi_k(4/3)$ are
the zeroes of $f_n$ and $l_n$.  The zeroes of $g_n$ and $f_n$ are
forbidden eigenvalues, and correspond to even values of $k$.  Thus
$\phi_3(0), \phi_5(0), \ldots, \phi_{2n-1}(0)$ are the zeroes of $h_n$
($\phi_1(0) = 0$, of course), and $\phi_1(4/3), \phi_3(4/3), \ldots,
\phi_{2n-3}(4/3)$ are the zeroes of $l_n$.  The exact values
\begin{align*}
  \phi_{2j-1}(4/3) 
  = \frac{1 + \cos \frac{2\pi(n-j)}{2n-1}}{3} 
  = \frac{2 \sin^2 \frac \pi 2 ( \frac{2j-1}{2n-1})}{3}
\end{align*}
are computed in \cite{zhou}.  If $j$ is small compared to $n$, which
will always happen if we fix $j$  and let $n \to \infty$, then
$\phi_{2j-1}(4/3) \approx \pi^2/6 (\frac{2j-1}{2n-1})^2$.  Since
$\psi_n(t) \sim t$ for $t$ near 0 we have
\begin{align*}
  \rho_n \psi_n(\phi_{2j-1}(4/3)) \sim \frac{4n-3}{2n-1} 
  \frac {\pi^2}6 (2j-1)^2 \to \frac{4\pi^2}3 (j-1/2)^2
\end{align*}
as $n \to \infty$.  

There is no exact computation of the zeroes of $h_n$, but the zeroes
of $g_n$ are known, so
\begin{align*}
  \phi_{2j}(0) = \frac{1 + \cos (\frac{2n-1-2j}{2n-1})\pi}{3} =
  \frac{2 \sin^2 \frac \pi 2 (\frac{2j}{2n-1})}{3}
\end{align*}
and we have interlacing of zeroes of $g_n$ and $h_n$, so
\begin{align*}
  \phi_{2j-2}(0) < \phi_{2j-1}(0) < \phi_{2j}(0).
\end{align*}
This implies
\begin{align*}
  \frac{4n-3}{2n-1} \frac{\pi^2}{6} (2j-2)^2 \le \rho_n \psi_n
  (\phi_{2j-1}(0)) \le \frac{4n-3}{2n-1} \frac{\pi^2}{6} (2j)^2.
\end{align*}
If we assume that the lower bound is the asymptotically correct value,
then we obtain the expected value $\frac{4 \pi^2}{3} (j-1)^2$ for the
limit. We will show below that this is indeed correct.

\begin{figure}
  \centering
  \includegraphics[width=2.9in]{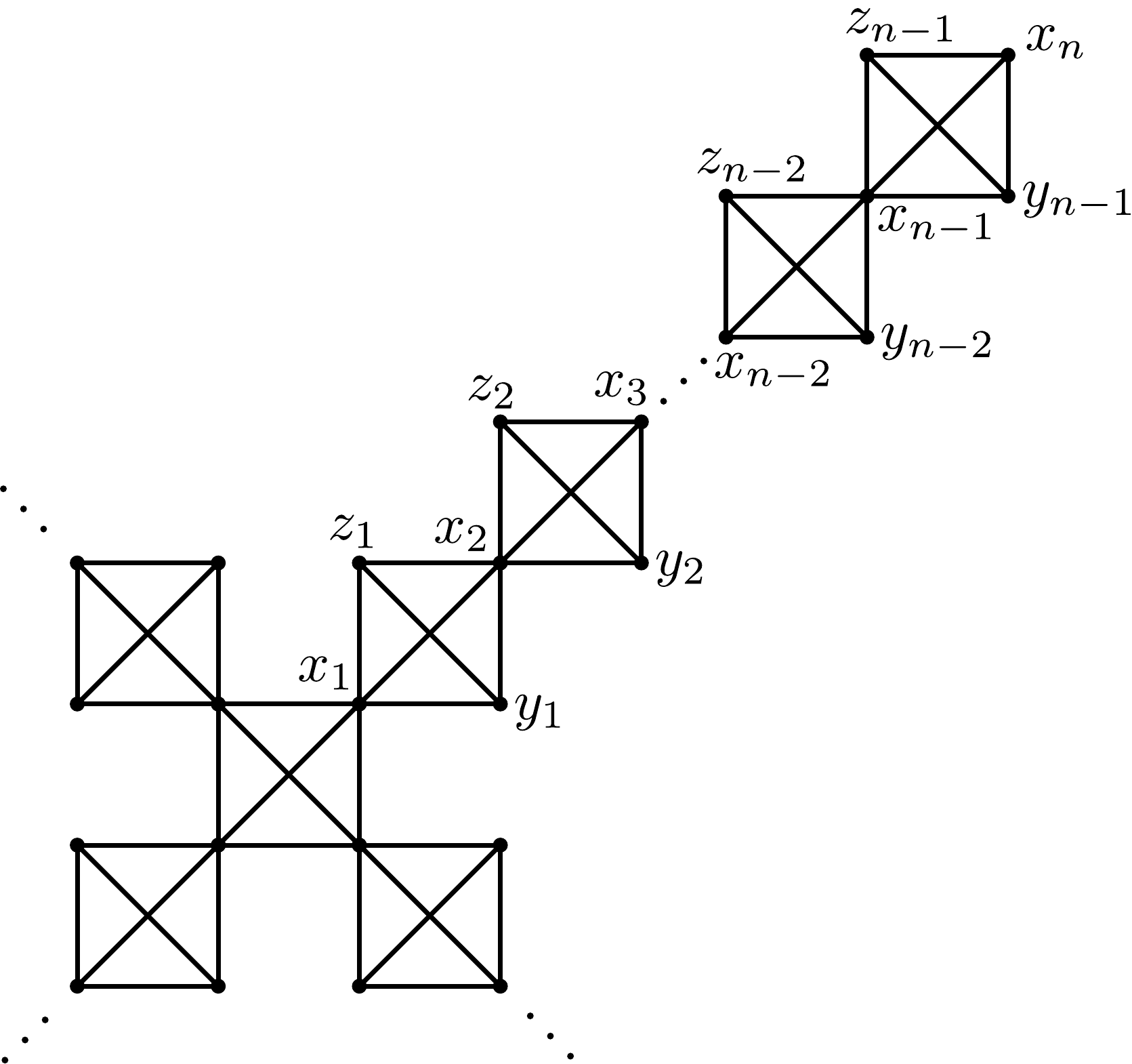}
  \caption{Vertices on one arm of $V_1$. \label{fig:armlabels}}
\end{figure}

We can also understand why the $\VS_n$ eigenfunctions, restricted to
the cross, converge to the eigenfunctions of the cross.  To see this,
we look at the graph eigenvalue equation on $V_1$.  Note that $V_1$
consists of four arms of $n-1$ squares joined at a central square.  We
label the diagonal vertices of one arm $x_1, x_2, \ldots, x_n$ and the
below and above diagonal vertices $y_1, \ldots, y_{n-1}$ and $z_1,
\ldots, z_{n-1}$ (see Figure~\ref{fig:armlabels}).  
By symmetry we will have $u(y_j) = u(z_j)$ for every eigenfunction.
The eigenvalue equation (with eigenvalue $\lambda_1$) at $y_j$ says 
\begin{align*}
  (1 - \lambda_1) u(y_j) = \tfrac 13 (u (x_j) + u(x_{j+1}) +
  u(y_j))
\end{align*}
So we obtain
\begin{align*}
  u(z_j) = u(y_j) = \frac{u(x_j) + u(x_{j+1})}{2 - 3 \lambda_1}.
\end{align*}
For $2 \le j \le n-1$ the eigenvalue equation at $x_j$ is 
\begin{align*}
  (1- \lambda_1) u(x_j) &= \frac 16 \Big( u (x_{j-1}) + u(x_{j+1}) +
  2 u(y_j) + 2 u(y_{j+1}) \Big)\\
  &= \frac 16 \left( u(x_{j-1}) + u(x_{j+1}) +
  \frac{2 u(x_{j-1}) + 4 u(x_j)+ 2 u(x_{j+1})}{2 - 3
  \lambda_1}\right).
\end{align*}
We can simplify this equation to
\begin{align*}
  (1-3\lambda_1) u(x_j) = \tfrac 12(u(x_{j-1}) + u(x_{j+1})).
\end{align*}
Note that this is exactly the eigenvalue equation (for eigenvalue
$3\lambda_1$) on the interior of the linear graph $x_1, \ldots, x_n$.
Similarly, at the endpoint $x_n$ the eigenvalue equation is 
\begin{align*}
  (1 - \lambda_1) u(x_n) &= \frac 13 \Big(u(x_{n-1}) + 2 u(y_{n-1})\Big)\\
  &=\frac 13\left(u(x_{n-1}) + \frac{2 u(x_{n-1}) + 2 u(x_n)}{2 - 3
  \lambda_1}\right),
\end{align*}
which simplifies to 
\begin{align*}
  (1 - 3 \lambda_1) u(x_n) = u(x_{n-1}),
\end{align*}
and this is the correct eigenvalue equation (for eigenvalue $3
\lambda_1$) with Neumann conditions at that endpoint.

The equation at the endpoint $x_1$ will depend on whether we are
looking at the 0-series or the 4/3-series.  For the 0-series the
values along all four arms will be identical, so the eigenvalue
equation is
\begin{align*}
  (1-\lambda_1)(u(x_1)) &= \frac 16 \Big(3 u(x_1) + u(x_2) + 2 u(y_1)\Big)\\
  &= \frac 16\left(3 u(x_1) + u(x_2) + \frac{2 u(x_1) + 2 u(x_2)}{2 - 3
  \lambda_1}\right),
\end{align*}
which simplifies to
\begin{align*}
  \left(1 - 18 \frac{1-\lambda_1}{4-3\lambda_1}\lambda_1\right) u(x_1)
  = u(x_2).
\end{align*}
For the 4/3-series the sum of the values on all four arms will be
zero, so the eigenvalue equation is 
\begin{align*}
  (1 - \lambda_1) u(x_1) &= \frac 16 \Big(-u(x_1) + u(x_2) + 2 u(y_1)\Big)\\ &= \frac 16\left(-u(x_1) + u(x_2) + \frac{2 u(x_1) + 2
  u(x_2)}{2 - 3 \lambda_1}\right),
\end{align*}
which simplifies to 
\begin{align*}
  (3 - 6 \lambda_1)u(x_1) = u(x_2).
\end{align*}

\begin{figure}
  \centering
  \includegraphics[width=5in]{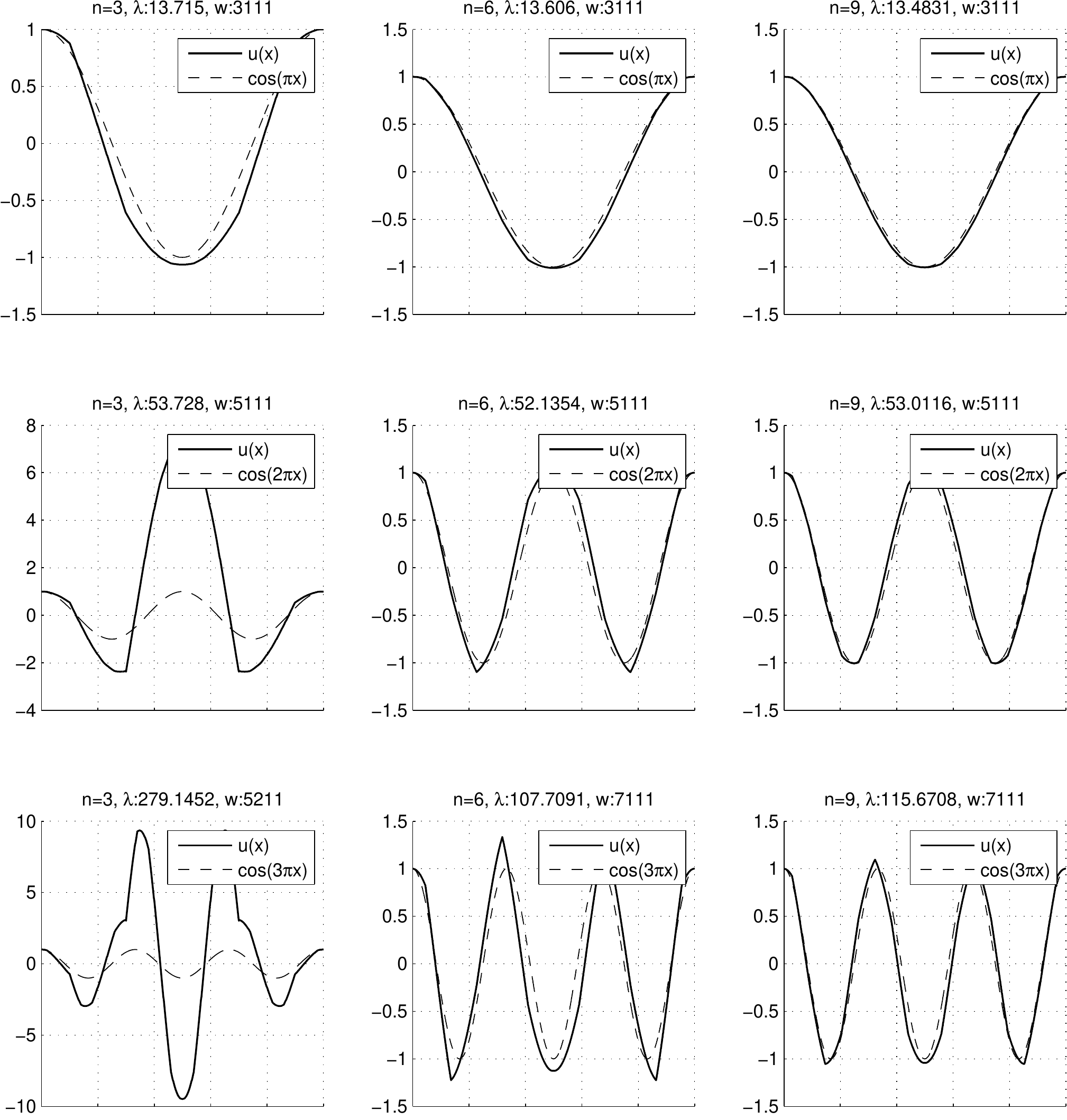}
  \caption{0-series eigenfunctions on $\VS_n$ for $n=3,6,9$ plotted
  against the corresponding symmetric eigenfunctions on the cross.
  \label{fig:convergence}}
\end{figure}

These should be compared with the eigenvalue equation for the
eigenfunction $\tilde u$ with eigenvalue $3 \lambda_1$ on two copies
of the linear graph with even and odd symmetries, namely
\begin{align*}
  (1 - 3 \widetilde\lambda_1) \tilde u(x_1) = \tfrac 12 (\tilde u(x_2) \pm \tilde u(x_1)).
\end{align*}
Note that we get the identical equation in the odd case, but in the
even case we get
\begin{align*}
  (1 - 6 \widetilde\lambda_1 ) \tilde u(x_1) = \tilde u(x_2),
\end{align*}
so there is a significant distinction. In the case of the 4/3-series,
we can therefore identify the restriction of the eigenfunctions to the
diagonal with
\begin{align*}
  \tilde u(x_k) = \sin \pi (j-1/2) \left(\frac{2k-1}{2n-1}\right), 
  \qquad 1 \le j \le n-1.  
\end{align*}
Figure~\ref{fig:convergence} shows some 0-series eigenfunctions plotted
against the symmetric eigenfunctions on the cross for $n=3,6,9$.
It appears that $u$ closely approximates 
\begin{align*}
  \tilde u(x_k) = \cos \pi j \left(\frac{2k-1}{2n-1}\right).
\end{align*}

We now sketch a proof that the eigenvalues $3\lambda_1$ approach
\begin{align*}
  3\widetilde\lambda_1
  = 1 - \cos \frac{2\pi j}{2n-1}
  = 2 \sin^2 \frac{\pi j}{2n-1}
\end{align*}
and the eigenvectors $u(x_k)$ approach $\tilde u(x_k)$ as $n \to
\infty$. Here we fix the value of $j$, and we require the appropriate
error estimate since both $3\lambda_1$ and $3\widetilde\lambda_1$ tend to
zero. The idea is to use standard perturbation theory, using the fact
that the two eigenvalue equations differ only at the single point
$x_1$, and the fact that the eigenvector $\tilde u(x_k)$ is fairly
uniformly distributed, so the value $\tilde u(x_1)$ is relatively
small.

Let $E$ denote the symmetric $n\times n$ matrix of tridiagonal form, with
\begin{align*}
  E_{kk} = 
  \begin{cases}
    1 & 2 \le k \le n-1 \\
    \frac 12 & k = 1 \text{~or~} n
  \end{cases}
\end{align*}
and $E_{k(k+1)} = E_{k(k-1)} = -1/2$. Let $\widetilde G$ denote the
diagonal matrix with
\begin{align*}
  \widetilde G_{kk} = 
  \begin{cases}
    1 & 1 \le k \le n- 1 \\
    \frac 12 & k = n
  \end{cases}
\end{align*}
and let $G$ denote the diagonal matrix with
\begin{align*}
  G_{kk} = 
  \begin{cases}
    1 & 2 \le k \le n-1 \\
    \frac 12 & k = n \\
    \frac{3\lambda_1-3}{3\lambda_1-4} & k = 1.
  \end{cases}
\end{align*}
Then the two eigenvalue equations may be written
\begin{align*}
  Eu &= 3\lambda_1 Gu,\\
  E\tilde u &= 3\widetilde \lambda_1 \widetilde G\tilde u.
\end{align*}
Note that the first equation is not a linear generalized eigenvalue
equation because $G$ depends on $\lambda_1$, but this does not really
matter in our argument.

The gist of the argument is that $\widetilde G - G$ is a matrix with only
one non-zero entry ($\widetilde G_{11} - G_{11}$) and we can bound this
entry since $\lambda_1$ is bounded away from $4/3$ for the 0-series;
and also we know $\tilde u$ exactly, hence $|\tilde u(x_1)| \le 1$
while $\langle \widetilde G \tilde u, \tilde u \rangle = n/2$. This yields
the estimate
\begin{align*}
  \frac{\langle (\widetilde G - G) \tilde u, \tilde u \rangle}
  {\langle\widetilde G \tilde u, \tilde u \rangle}
  = O\left( \frac 1n \right).
\end{align*}
With a little more work, we can get the estimate
\begin{align*}
  \frac{\langle (\widetilde G - G) \tilde u, u \rangle}
  {\langle\widetilde G \tilde u, u \rangle}
  = O\left( \frac 1n \right)
\end{align*}
for the first $N$ eigenfunctions ($N$ is fixed as $n \to \infty$). But
this is exactly what we need to estimate $\lambda_1 - \widetilde
\lambda_1$. If we take the inner product of the first eigenvalue
equation with $\tilde u$ and the second eigenvalue equation with $u$,
then using the symmetry of all the matrices we obtain
\begin{align*}
  \lambda_1 \langle G \tilde u, u  \rangle 
  = \widetilde \lambda_1  \langle \widetilde G \tilde u, u \rangle
\end{align*}
hence
\begin{align*}
  \lambda_1 - \widetilde \lambda_1 = \lambda_1  
  \frac{\langle (\widetilde G - G) \tilde u, u \rangle}
  {\langle\widetilde G \tilde u, u \rangle}
  = O\left( \frac 1{n^3} \right)
\end{align*}
for the first $N$ eigenvalues, since we know $\widetilde\lambda_1 =
O(1/n^2)$.  With a little more work we can show that $u - \tilde u =
O(1/n)$ when $u$ is properly normalized.
\begin{figure}
  \centering
  \includegraphics[width=2.9in]{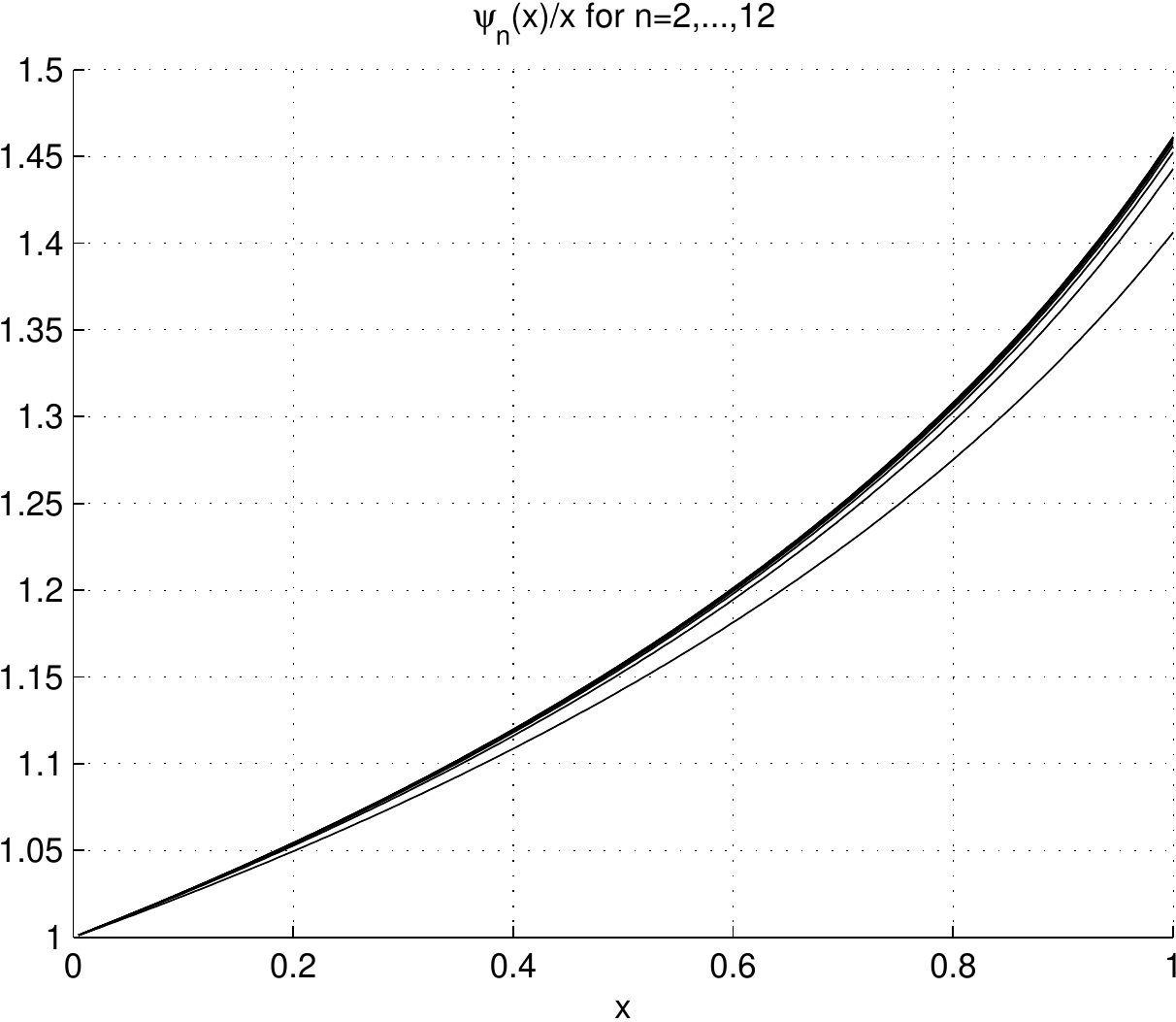}
  \caption{$\psi_n(x)/x$ on $[0,1]$ for $n=2,\ldots,12$
  ($\psi_{n+1} \ge \psi_n$). \label{fig:psi}}
\end{figure}

So far we have dealt with the level 1 eigenvalues $\lambda_1$. The
actual eigenvalues $\lambda$ on $\VS_n$ are given by $\lambda =
\psi_n(\lambda_1)$ for the lowest segment of the spectrum (this will
include the first $N$ eigenvalues once $n$ is large enough).
Figure~\ref{fig:psi} gives experimental evidence for the estimate $t
\le \psi_n(t) \le t + ct^2$ on $0 \le t \le 1$ for a constant $c$
independent of $n$. This shows that $\lambda-\widetilde\lambda_1 =
O(1/n^3)$ as $n\to \infty$ for the first $N$ eigenvalues.

\section{Weyl Ratio}\label{sec:weyl}

We now describe in more detail the Weyl ratio $W_n(t) =
N_n(t)/t^{\alpha_n}$ on $\VS_n$, where 
\begin{align*}
  \alpha_n =
  \frac{\log(4n-3)}{\log \rho_n}, 
  \qquad \text{and} \qquad
  N_n(t) = \sum_{\lambda_j \le t} m(\lambda_j) 
\end{align*}
is the counting function for the number (counting multiplicity) of
eigenvalues.  According to a general theorem of Kigami and Lapidus
\cite{lapidus},
$w_n(t) = \lim_{k \to \infty} W_n(\rho_n^k t)$ exists.  In order to
compare $w_n$ for different values of $n$, we normalize by
$\tilde w_n(s) = w_n(\lambda_1 \rho_n^s)$ so that $\tilde w_n$ is a
periodic function of period 1 with $\tilde w_n(0) = w_n(\lambda_1)$.

From the data it appears that $\tilde w_n$ is converging to a limit as
$n \to \infty$, but this limit has nothing to do with the Weyl ratio
on the cross, which tends to a constant. While we cannot supply a
complete explanation of this phenomenon, we can make a few
observations about the behavior of $\tilde w_n(s)$ for some values of
$s$. Because of high multiplicities the functions $w_n$ and $\tilde w_k$
have jump discontinuities.  We write $w_n(t^-) = \lim_{s \to t^-}
w_n(s)$ and similarly for $\tilde w_k(s^-)$.  First we note that it
is possible to compute $w_n(\lambda_i)$ for small values of $i$.

\begin{lem}
For $i \le j\le n-1$ we have 
\begin{align*}
  w_n(\lambda_{2j-1}^-) = \frac{4j-3}{(\lambda_{2j-1})^{\alpha_n}},
  \qquad
  w_n(\lambda_{2j-1}) = \frac{4j-1}{(\lambda_{2j-1})^{\alpha_n}},
  \qquad
  w(\lambda_{2j}) = \frac{4j}{(\lambda_{2j})^{\alpha_j}}.
\end{align*}
\end{lem}

\begin{proof}
A simple induction argument shows that
\begin{align*}
  N(\rho_n^k \lambda_{2j-1}) &= (4j-1)(4n-3)^k + 1,\\
  m(\rho_n^k \lambda_{2j-1}) &= 2(4n-3)^k + 1,
\end{align*}
since $\rho_n^k \lambda_{2j-1}$ is a 4/3-series eigenvalue born on level
$k$.  Thus
\begin{align*}
  W(\rho_n^k \lambda_{2j-1}) 
  = \frac{(4j-1)(4n-3)^k + 1}{(\rho_n^k \lambda_{2j-1})^{\alpha_n}} 
  = \frac{(4j-1)(4n-3)^k+1}{(4n-3)^k (\lambda_{2j-1})^{\alpha_n}}
\end{align*}
and the computation of $w_n(\lambda_{2j-1})$ follows by taking the
limit.  Similarly we obtain the result for $w_n(\lambda_{2j-1}^-)$.
We can also show by induction that $N(\rho_n^k \lambda_{2j}) = 4j(4n-3)^k
+1$, and the result for $w_n(\lambda_{2j})$ follows.
\end{proof}

In particular, we have $\tilde w_n(0^-) = 1/\lambda_1^{\alpha_n}$ and
$\tilde w_n(0) = 3/\lambda_1^{\alpha_n}$.  As $n \to \infty$ we have
$\lambda_1 \to \pi^2/3$ and $\alpha_n \to 1/2$ so $\lim_{n \to \infty}
\tilde w_n(0^-) = \sqrt 3/\pi$ and $\lim_{n \to \infty} \tilde w_n(0)
= 3 \sqrt 3/\pi$.  It is more difficult to get information about
limiting behavior of $\tilde w_n(s)$ for other values of $s$ because
we would have to simultaneously let $j \to \infty$ as $n \to \infty$.
Although we know $\lambda_l \to \frac{4}{3} (\frac{l \pi}{2})^2$ for
fixed $l$ as $n \to \infty$, the convergence is not uniform in $l$.

Because the spectrum has large gaps on either side of $\rho_n^k
\lambda_1$, we can say more about the behavior of $\tilde w_n(s)$ for
$s$ near 0.  In fact the eigenvalue just below $\rho_n^k \lambda_1$ is
$\rho_n^k \psi_n(\phi_{2n-1}^{(k)}(0))$, and the eigenvalue just above it
is $\rho_n^{k+1} \psi_n(\phi_2 \phi_{2n-1}^{(k)}(0))$.  So for
\begin{align*}
  \rho_n^k \psi_n( \phi_{2n-1}^{(k)} (0)) < t < \rho_n^k \lambda_1, 
\end{align*}
the value of $N_n(t)$ is $(4n-3)^k$ so $W_n(t) =
(4n-3)^k/t^{\alpha_n}$.  In taking the limit as $k \to \infty$ we note
that $\phi_{2n-1}^{(k)}(0) \to p_n$, the fixed point of $\phi_{2n-1}$,
so $w_n(t) = 1/t^{\alpha_n}$ for $\psi_n(p_n) \le t <
\lambda_1$ or equivalently 
\begin{align*}
  \tilde w_n(s) = \frac{1}{\lambda_1^{\alpha_n}(\rho_n^{\alpha_n})^s} 
  \qquad \text{for~} -\left(\frac{\log
  \lambda_1 - \log \psi_n(p_n)}{\log \rho_n}\right) \le s < 0.
\end{align*}

Similarly $w_n(t) = 3/t^{\alpha_n}$ for $\lambda_1 \le t \le \rho_n
\psi_n(\phi_2(p_n))$, or equivalently 
\begin{align*}
  \tilde w_n(s) = \frac{3}{\lambda_1^{\alpha_n}(\rho_n^{\alpha_n})^s} 
  \qquad
  \text{for~} 0 \le s \le \left(\frac{\log \rho_n + \log \psi_n(\phi_2(p_n)) -
  \log \lambda_1}{\log \rho_n}\right).
\end{align*}

\end{document}